\documentclass[numbers,webpdf,imaiai]{ima-authoring-template}%
\usepackage{algorithm}
\usepackage{amsbsy,amsmath,latexsym,amsfonts, epsfig, color, amssymb, bm, caption}
\usepackage{bbm}
\usepackage{cases}
\usepackage{epsf,slidesec,epic,eepic}
\usepackage{epsfig}
\usepackage{epstopdf}

\usepackage{graphicx}
\usepackage{longtable}
\usepackage{makecell,booktabs}
\usepackage{subfig}
\DeclareMathAlphabet{\mathcal}{OMS}{cmsy}{m}{n}

\theoremstyle{thmstyletwo}%

\newtheorem{theorem}{Theorem}
\newtheorem{lemma}{Lemma}
\newtheorem{definition}{Definition}
\newtheorem{example}{Example}
\newtheorem{proposition}{Proposition}
\newtheorem{corollary}{Corollary}

\newtheorem{remark}{Remark}

\numberwithin{equation}{section}

\newcommand{\tp}{^\top}

\begin{document}

\DOI{DOI HERE}
\copyrightyear{2021}
\vol{00}
\pubyear{2021}
\access{Advance Access Publication Date: Day Month Year}
\appnotes{Paper}
\copyrightstatement{Published by Oxford University Press on behalf of the Institute of Mathematics and its Applications. All rights reserved.}
\firstpage{1}


\title[Qian et al.]{Error bound and exact penalty method for optimization problems with 
nonnegative orthogonal constraint}

\author{Yitian Qian
\address{\orgdiv{School of Mathematics}, \orgname{South China University of Technology}, \orgaddress{\state{Guangzhou}, \country{China}}}}
\author{Shaohua Pan
\address{\orgdiv{School of Mathematics}, \orgname{South China University of Technology}, \orgaddress{\state{Guangzhou}, \country{China}}}}
\author{Lianghai Xiao*
\address{\orgdiv{School of Mathematics}, \orgname{South China University of Technology}, \orgaddress{\state{Guangzhou}, \country{China}}}}

\authormark{Qian et al.}

\corresp[*]{Corresponding author: \href{email:lxiao@scut.edu.cn}{lxiao@scut.edu.cn}}

\received{9}{11}{2021}
\revised{4}{5}{2022}


\abstract{This paper is concerned with a class of optimization problems with 
 the nonnegative orthogonal constraint, in which the objective function is 
 $L$-smooth on an open set containing the Stiefel manifold ${\rm St}(n,r)$. 
 We derive a locally Lipschitzian error bound for the feasible points without 
 zero rows when $n>r>1$, and when $n>r=1$ or $n=r$ achieve a global Lipschitzian 
 error bound. Then, we show that the penalty problem induced by the elementwise 
 $\ell_1$-norm distance to the nonnegative cone is a global exact penalty, 
 and so is the one induced by its Moreau envelope under a lower second-order 
 calmness of the objective function. 
 A practical penalty algorithm is developed by solving approximately a series of 
 smooth penalty problems with a retraction-based nonmonotone line-search 
 proximal gradient method, and any cluster point of the generated sequence 
 is shown to be a stationary point of the original problem.
 Numerical comparisons with the ALM \citep{Wen13} and the exact penalty method 
 \citep{JiangM22} indicate that our penalty method has an advantage 
 in terms of the quality of solutions despite taking a little more time.}
 


\maketitle

 \section{Introduction}\label{sec1}

  Let $\mathbb{R}^{n\times r}\,(n\ge r)$ represent the vector space of all $n\times r$ 
  real matrices, equipped with the trace inner product $\langle \cdot,\cdot\rangle$ 
  (i.e., $\langle X,Y\rangle ={\rm tr}(X\tp Y)$ for any $X,Y\in\mathbb{R}^{n\times r}$) 
  and its induced Frobenius norm $\|\cdot\|_F$, and let $\mathbb{R}_{+}^{n\times r}$ 
  denote the cone consisting of all nonnegative matrices in $\mathbb{R}^{n\times r}$. 
  Denote by ${\rm St}(n,r)$ the Stiefel manifold in $\mathbb{R}^{n\times r}$, 
  	i.e., ${\rm St}(n,r):=\big\{X\in\mathbb{R}^{n\times r}\,|\, X\tp X=I_{r}\big\}$, 
  and write $\mathcal{S}_{+}^{n,r}:=\mathbb{R}_{+}^{n\times r}\cap{\rm St}(n,r)$. 
  We are interested in the following optimization problem with the nonnegative orthogonal constraint
 \begin{equation}\label{orth-prob}
  \min_{X\in\mathcal{S}_{+}^{n,r}}\!f(X),
 \end{equation}
 where $f\!:\mathcal{O}\to\mathbb{R}$ is a continuously differentiable ($\mathcal{C}^1$)  
 function on an open set $\mathcal{O}$ containing ${\rm St}(n,r)$. Due to the compactness 
 of $\mathcal{S}_{+}^{n,r}$, the global minimizer set of \eqref{orth-prob}, 
 denoted by $\mathcal{X}^*$, is nonempty. In addition, the compactness of ${\rm St}(n,r)$ and the continuous differentiability of $f$ on $\mathcal{O}$ implies that the function $f$ and its gradient $\nabla\!f$ are Lipschitz continuous on ${\rm St}(n,r)$. We denote their Lipschitz constants on ${\rm St}(n,r)$ by $L_{\!f}$ and $L_{\nabla\!f}$, respectively. 

 Problem \eqref{orth-prob} frequently arises from machine learning 
 and data sciences such as the nonnegative principal component analysis \cite{Montanari16,Zass07}, the nonnegative Laplacian embedding \cite{Luo09}, 
 the discriminative nonnegative spectral clustering \cite{YangY12}, and the orthogonal  
 nonnegative matrix factorization \cite{Ding06}. For more applications, the reader is 
 referred to \cite[Section 1.1]{JiangM22}. Note that the nonnegative constraint 
 $X\!\in\!\mathbb{R}_{+}^{n\times r}$ introduces a sparse structure to problem 
 \eqref{orth-prob}, i.e., for every $X\in\mathcal{S}_{+}^{n,r}$, 
 its each row has at most one nonzero entry. When $n=r$, the column orthogonality 
 makes the feasible set of \eqref{orth-prob} become a permutation matrix set.
 Consequently, problem \eqref{orth-prob} with $n=r$ covers some combinatorial
 optimization problems, and a typical example is the quadratic assignment problem \cite{Burkard97}:
 \begin{equation}\label{QAP}
  \min_{X\in\mathbb{R}^{n\times n}}\Big\{\langle A, XBX\tp \rangle\ \ {\rm s.t.}\ \
  X\tp X=I_{n},\ X\ge 0\Big\},
 \end{equation}
 where $A\in\mathbb{R}^{n\times n}$ and $B\!\in\mathbb{R}^{n\times n}$ are the given matrices.
 In addition, the discrete constraint $X\!\in\mathcal{S}_{+}^{n,n}$ is
 found to have many applications in graph matching of image science 
 (see \cite{Jiang17,Yu18}).
 \subsection{Related works}\label{sec1.1} 
  
  It is well known that error bound for optimization problems plays a significant role
  in establishing their global or local exact penalty. Because the constraint system 
  $\{X\in\mathbb{R}^{n\times r}\ |\ X\in\mathcal{S}_{+}^{n,r}\}$ can be expressed as 
  $\{X\in\mathbb{R}^{n\times r}\ |\ X\tp X-I=0,X_{ij}\ge 0\ \ {\rm for}\ i=1,\ldots,n, j=1,\ldots,r\}$, by invoking \citep[Theorem 2.2]{Luo94}, 
  for each compact set $\Omega\subset\mathbb{R}^{n\times r}$,
  there exist constants $\tau>0$ and $\sigma>0$ such that for all $X\in\Omega$, 
  \begin{equation*}
    {\rm dist}(X,\mathcal{S}_{+}^{n,r})\le\tau\big[\vartheta(X)+\|X\tp X-I\|_F\big]^{\sigma},
  \end{equation*}  
  where ${\rm dist}(X,\Xi)$ for a closed set $\Xi\subset\mathbb{R}^{n\times r}$ 
  denotes the distance of $X$ to $\Xi$ on the Frobenius norm, and $\vartheta(X):=\langle E,\max(0,-X)\rangle$ with $E$ being the matrix of all ones 
  is precisely the elementwise $\ell_1$-norm distance from $X$ to $\mathbb{R}_{+}^{n\times r}$. From Lemma \ref{Bound-relation} in Appendix A,
  ${\rm dist}(X,{\rm St}(n,r))\le\|X\tp X-I\|_F
  \le(1+\|X\|){\rm dist}(X,{\rm St}(n,r))$ for any $X\in\mathbb{R}^{n\times r}$. The above error bound is equivalent to saying that for each compact set 
  $\Omega\subset\mathbb{R}^{n\times r}$, there exist $\tau>0$ and $\sigma>0$ such that for all $X\in\Omega$,
  \begin{equation}\label{ebound0}
  	{\rm dist}(X,\mathcal{S}_{+}^{n,r})
  	\le\tau\big[\vartheta(X)+{\rm dist}(X,{\rm St}(n,r))\big]^{\sigma}.
  \end{equation}
  To the best of our knowledge, the explicit exponent $\sigma$ is unknown for the constraint 
  system $X\in\mathcal{S}_{+}^{n,r}$. Recently, Jiang et al. \citep{JiangM22} derived an error bound for this constraint system by using
  the characterization 
  $\mathcal{S}_{+}^{n,r}=\big\{X\in\mathcal{O}\mathcal{B}_{+}^{n,r}\ |\ \|XV\|_F=1\big\}$, 
  where $\mathcal{O}\mathcal{B}_{+}^{n,r}$ is the nonnegative multiple spherical set 
  and $V\in\mathbb{R}^{n\times r}$ is a constant matrix with $\|V\|_F=1$ and $VV\tp$ 
  having positive entries, but now it is unclear whether the error bound in \citep{JiangM22} implies the one in \eqref{ebound0} for an explicit $\sigma$. 
  
  Problem \eqref{orth-prob} is a special nonlinear program, and in the past several decades, there were extensive research on the exact penalty method for nonlinear programming. We find that most of existing results are about local exact penalty. For example, Han and Magasarian \citep{Han79} 
  established the exactness of the popular $\ell_1$-penalty function in a local sense, 
  Bonnans and Shapiro \citep[Section 3.4.2]{Bonnans03} achieved a local exact 
  penalty result for a general conic constrained optimization problem 
  under a local error bound condition, and Burke \citep{Burke91} showed that 
  the local exact penalization of a general nonlinear programming problem is equivalent to 
  the calmness of the problem itself. Although a global exact penalty result 
  was partly achieved in \citep[Theorem 4.1]{Han79}, the required condition 
  is unverifiable. For problem \eqref{orth-prob} with $n=r$, by noting that 
  $\mathcal{S}_{+}^{n,n}=\{X\in\mathbb{R}_{+}^{n\times n}\ |\ Xe=X\tp e =e\}\cap
  \big\{X\in\mathbb{R}^{n\times n}\ |\ \|X\|_0=n\big\}$, 
  Jiang et al. \citep{JiangLW16} showed that the penalty problem induced by 
  the zero-norm constraint is a global exact penalty and then achieved 
  the $\ell_p$-norm global exact penalty by the relation between 
  the $\ell_p\,(0<p<1)$-norm and the zero-norm. In addition, for problem \eqref{orth-prob}, 
  Jiang et al. \citep{JiangM22} also proposed a class of exact penalty by 
  the derived error bound, but the global optimal solution of their penalty problem 
  is not necessarily globally optimal to \eqref{orth-prob} though its sign belongs to 
  the optimal sign set of \eqref{orth-prob}, so this class of penalty problems is not strictly a global exact penalty of \eqref{orth-prob}. Thus, for \eqref{orth-prob}, it is still an open question whether the penalty problem induced by the elementwise $\ell_1$-norm distance to $\mathbb{R}_{+}^{n\times r}$ is a global exact penalty.
 \subsection{Our contributions}\label{sec1.2} 
  
  In this work, we investigate the local Lipschitzian error bound for the constraint system $X\in\mathcal{S}_{+}^{n,r}$ and demonstrate that the penalty problems induced by the elementwise $\ell_1$-norm distance to $\mathbb{R}_{+}^{n\times r}$ and its Moreau envelope (see section \ref{sec2} for its definition) are global exact penalty for \eqref{orth-prob}, i.e., the penalty problems induced by them with the penalty parameter greater than a certain threshold have the same global optimal solution set as problem \eqref{orth-prob} does. 
  
  One contribution is to derive a local Lipschitzian error bound 
  for those feasible points without zero rows when $n>r>1$ and a global Lipschitzian error bound for those feasible points with $n>r=1$ or $n=r$, which means that for such feasible points of \eqref{orth-prob}, there is a neighborhood such that the error bound \eqref{ebound0} holds with $\sigma=1$. Based on this, we also demonstrate that the following penalty problem 
  \begin{equation}\label{Nsmooth-pen}
  \min_{X\in{\rm St}(n,r)}\,\Phi_{\!\rho}(X):=f(X)+\rho\vartheta(X)
  \end{equation}
  is a global exact penalty for \eqref{orth-prob} if each global minimizer
  has no zero rows when $n>r>1$. Because $\vartheta(X)=0$ if and only if $e_{\gamma}\vartheta(X)=0$ 
  or $\|\max(0,-X)\|_F^2=0$, where $e_{\gamma}\vartheta$ denotes the Moreau envelope of $\vartheta$ associated to constant $\gamma>0$, we also consider the following smooth penalty problems   
  \begin{equation}\label{smooth-pen}
  \min_{X\in{\rm St}(n,r)}\,\Theta_{\rho,\gamma}(X):=f(X)+\rho e_{\gamma}\vartheta(X)
  \end{equation}
  and
 \begin{equation}\label{Quad-Pen}
 \min_{X\in{\rm St}(n,r)}\Theta_{\rho,0}(X):=f(X)+\rho\|\max(0,-X)\|_F^2.
 \end{equation}  
  Under a lower second-order calmness of $f$ relative to a neighborhood of 
  $\mathcal{X}^*$ (see equation \eqref{growth}), we show that problems \eqref{smooth-pen}-\eqref{Quad-Pen} are a global exact penalty for \eqref{orth-prob} if each global minimizer has no zero rows when $n>r>1$. The lower second-order calmness of $f$ 
  relative to a neighborhood of $\mathcal{X}^*$ is very weak and does not necessarily 
  imply its local Lipschitz continuity relative to this neighborhood. 
  The assumption that each global minimizer has no zero rows when $n>r>1$ is a little strong, but it cannot be removed (see Example \ref{example1}). Even so, our global exact penalty results are still applicable to those problems with $n=r$. Different from the exact penalty problems in \cite{JiangLW16,JiangM22}, our penalty problems \eqref{Nsmooth-pen}-\eqref{Quad-Pen} involve the manifold constraint $X\in{\rm St}(n,r)$ and their computation requires more 
  cost, but when solving them with a manifold algorithm, 
  the iterate sequence can satisfy the orthogonal constraint very well, 
  which can lead to better solutions especially for those difficult 
  problems (see Section \ref{sec6.3}).

  The other contribution is to develop a practical penalty algorithm by seeking the approximate stationary points to a series of smooth penalty problems 
  \eqref{smooth-pen} or \eqref{Quad-Pen} with increasing $\rho$, and verify that 
  any cluster point of the generated sequence is a stationary point of \eqref{orth-prob}. 
  Although the global exact penalty of \eqref{smooth-pen} and \eqref{Quad-Pen} requires that each global minimizer of \eqref{orth-prob} has no zero rows when $n>r>1$, this convergence result does not depends on this requirement. Unlike the penalty algorithm in \citep{JiangM22}, this asymptotic convergence removes the no zero row assumption on the cluster point. In addition, our penalty algorithm will return a local minimizer of \eqref{orth-prob} provided that a local minimizer $X^{l+1}$ of some penalty problem \eqref{EP-subprob} associated to $\rho_{l}>\widehat{\rho}$ lies in a neighborhood of a point of $\mathbb{R}_{+}^{n\times r}$ which is locally optimal to \eqref{smooth-pen} or \eqref{Quad-Pen} associated to $\widehat{\rho}$.
  This result does not require the local optimality of every iterate. 
    
  To seek an approximate stationary point to a single penalty problem \eqref{smooth-pen} or \eqref{Quad-Pen}, we propose a nonmonotone line-search proximal gradient method (PGM) on manifold, and show that if $f$ is definable in an o-minimal structure $\mathscr{O}$ 
  over $\mathbb{R}$ \cite{van1998tame}, the whole iterate sequence is convergent under a  restriction on the penalty objective value sequence; see Theorem \ref{theorem1-ManPG1}. 
  For a fixed $\rho>0$, due to the smoothness of $\Theta_{\rho,\gamma}$ with $\gamma\ge 0$, the existing algorithms proposed in \cite{Wen13,Huang15} for manifold optimization can be directly applied to \eqref{smooth-pen}, but only the subsequential convergence was obtained 
  in \cite{Wen13} and the convergence of whole sequence was obtained in \cite{Huang15} 
  under the strong retraction convexity of the cost function.

  For the proposed penalty method based on $\Theta_{\rho,\gamma}$ (SEPPG+, for short) and  $\Theta_{\rho,0}$ (SEPPG0, for short), we conduct numerical experiments on the 134 QAPLIB instances, graph pair matching on CMU house image dataset, projection problems onto $\mathcal{S}_{+}^{n,k}$, and orthogonal nonnegative matrix factorization. Numerical comparisons with the ALM \cite{Wen13} and the exact penalty method EP4Orth+ \citep{JiangM22} demonstrate that SEPPG+ and SEPPG0 can yield better solutions than both ALM and EP4Orth+ do for those difficult problems such as QAPs. For graph pair matching problems, SEPPG+, SEPPG0 and ALM yield the matching accuracy and objective values 
  comparable to those yielded by FGM-D in \cite{Zhou15}, which are much better 
  than those yielded by EP4Orth+. Here, FGM-D is a path-following algorithm 
  with a heuristic strategy designed by the convex and concave relaxations 
  for \eqref{orth-prob}. For projection problems onto $\mathcal{S}_{+}^{n,k}$ 
  and orthogonal nonnegative matrix factorization, the solutions yielded by 
  SEPPG+ and SEPPG0 have comparable quality with those yielded by EP4Orth+ and ALM.

  It is worth pointing out that a similar penalty method can be developed by 
  solving a series of nonsmooth exact penalty problems \eqref{Nsmooth-pen} with 
  a nonmonotone line-search retraction-based PGM as in \cite{Chen20} or 
  the Riemannian PGM as in \cite{Huang21}, but the nonsmoothness of $\vartheta$ 
  makes the asymptotic convergence analysis require the exact stationary point of 
  every penalty problem, which is impossible in practical computation. In addition, 
  our testing on the \textbf{134} instances in QAPLIB indicates that 
  the penalty algorithm based on \eqref{Nsmooth-pen} has worse performance than 
  the one by \eqref{smooth-pen}. We leave the development of exact penalty 
  methods by \eqref{Nsmooth-pen} for a future research topic.
   
 \subsection{Notation}\label{sec1.3}
 
 Throughout this paper, $M$ represents the Stiefel manifold ${\rm St}(n,r)$, 
 an embedded submanifold of the Euclidean space $\mathbb{R}^{n\times r}$, 
 $\mathcal{E}$ denotes a Euclidean space endowed with the Euclidean norm 
 $\|\cdot\|$, $(\mathcal{M},\langle\cdot,\cdot\rangle_{g})$ denotes a finite 
 dimensional Riemannian manifold, $\mathcal{C}^1(\mathcal{M})$ 
 represents the collection of all $\mathcal{C}^1$ functions on $\mathcal{M}$, 
 and $\mathbb{S}^r$ means the vector space consisting of all $r\times r$ 
 real symmetric matrices. For each $k\in\mathbb{N}$, write $[k]:=\{1,2,\ldots,k\}$, 
 and for each $x$, let $[\![x]\!]$ be the subspace spanned by $x$.
 Let $e$ denote a vector of all ones with dimension known from the context, 
 and let $I_{n\times r}$ represent the $n\times r$ rectangular identity matrix. 
 For each $X\in\mathbb{R}^{n\times r}$, let $\|X\|_1,\|X\|_*$ and $\|X\|$ 
 respectively denote the elementwise $\ell_1$-norm, the nuclear norm, and 
 the spectral norm of $X$, and write ${\rm supp}(X):=\{(i,j)\in[n]\times[r]\ |\ |X_{ij}|\ne 0\}$ 
 and $\overline{{\rm supp}}(X):=\{(i,j)\in[n]\times[r]\ |\ X_{ij}=0\}$.
 For each $X\in\mathbb{R}^{n\times r}$ and each $(i,j)\in[n]\times[r]$, 
 $X_{i\cdot}$ denotes the $i$th row of $X$ and $X_{\cdot j}$ means 
 the $j$th column of $X$, and $\mathbb{B}(X,\delta)$ denotes the closed 
 ball on the Frobenius norm centered at $X$ with radius $\delta$. 
 For a closed set $\Omega\subseteq\mathcal{E}$, 
 ${\rm Proj}_{\Omega}$ denotes the projection operator onto $\Omega$ 
 in the sense of the Euclidean norm of $\mathcal{E}$, 
 ${\rm dist}(x,\Omega)\!:=\inf_{z\in\Omega}\|z\!-\!x\|=\|{\rm Proj}_{\Omega}(x)-x\|$, 
 and $\delta_\Omega$ denotes the indicator function of $\Omega$, 
 i.e., $\delta_\Omega(x)=0$ if $x\in \Omega$ and $+\infty$ otherwise.
 For a given $x\in\Omega$, $\mathcal{T}_{\Omega}(x),\widehat{\mathcal{N}}_{\Omega}(x)$ 
 and $\mathcal{N}_{\Omega}(x)$ denote the tangent cone, the regular normal cone 
 and the limiting (or Mordukhovich) normal cone to $\Omega$ at $x$, respectively. 
 We denote by ${\rm T}_{x}\mathcal{M}$ and ${\rm N}_{x}\mathcal{M}$ respectively
 the tangent space and the normal space to $\mathcal{M}$ at $x\in\mathcal{M}$, 
 and endow the tangent space of $M$ with the Riemannian metric 
 $\langle\cdot,\cdot\rangle_{x}$ induced from the trace inner product, 
 i.e., $\langle G,H\rangle_{x}={\rm tr}(G\tp H)$ for any $G,H\in{\rm T}_{x}M$. 
 
 For a function $h\!:\mathcal{E}\to\overline{\mathbb{R}}\!:=\mathbb{R}\cup\{+\infty\}$ 
 and a compact submanifold $\mathcal{M}$ of $\mathcal{E}$,
 $h_{|\mathcal{M}}\!: \mathcal{M}\to\overline{\mathbb{R}}$ represents
 a restriction of $h$ to $\mathcal{M}$, and if $h$ is differentiable at $x\in\mathcal{M}$, 
 ${\rm grad}\,h(x)$ denotes the Riemannian gradient of $h$ at $x$, 
 the unique tangent vector at $x$ such that for any $\eta\in{\rm T}_{x}\mathcal{M}$,
 $\langle\nabla h(x),\eta\rangle=\langle{\rm grad}\,h(x),\eta\rangle$, 
 which is also the gradient of $h_{|\mathcal{M}}$ at $x$. 
 When $\mathcal{E}=\mathbb{R}^{n\times r}$ 
 and $\mathcal{M}=M$, by our choice of the Riemannian metric, 
 ${\rm grad}\,h(x)={\rm Proj}_{{\rm T}_{x}M}(\nabla h(x))$. For a given $X\in M$, 
 ${\rm N}_{\!X}M=\big\{XS\,|\,S\in\mathbb{S}^r\big\}$ and 
 ${\rm T}_{\!X}M=\big\{H\in\mathbb{R}^{n\times r}\,|\,\mathcal{A}_{X}(H)=0\big\}$
 with $\mathcal{A}_{X}(H)\!:=X\tp H+H\tp X$, and moreover,
 \begin{equation}\label{Proj-St}
   {\rm Proj}_{{\rm T}_{\!X}M}(Z)\!:=Z-\!\frac{1}{2}X\mathcal{A}_{X}(Z)
   \quad\forall Z\in\mathbb{R}^{n\times r}.
 \end{equation} 
 \section{Preliminaries}\label{sec2}
 
 For a proper lsc function $h\!:\mathbb{R}^{n\times r}\to\overline{\mathbb{R}}$
 and a constant $\gamma>0$, $P_{\gamma}h$ and $e_{\gamma}h$ denote the proximal 
 mapping and Moreau envelope of $h$ associated to $\gamma$, respectively, defined by
 \begin{align}\label{prx-mr}
   \mathcal{P}_{\!\gamma}h(x)
   \!:=\mathop{\arg\min}_{z\in\mathbb{R}^{n\times r}}
   \Big\{\frac{1}{2\gamma}\|z-x\|_F^2+h(z)\Big\}\ \ {\rm and}\ \
   e_{\gamma}h(x)\!:=\!\min_{z\in\mathbb{R}^{n\times r}}
   \Big\{\frac{1}{2\gamma}\|z-x\|_F^2+h(z)\Big\}.
 \end{align}
 When $h$ is convex, $\mathcal{P}_{\!\gamma}h$ is a Lipschitz continuous mapping 
 of modulus $1$ from $\mathbb{R}^{n\times r}$ to $\mathbb{R}^{n\times r}$,
 and $e_{\gamma}h$ is a continuously differentiable convex function with
 $\nabla\!e_{\gamma}h(x)=\gamma^{-1}(x\!-\!\mathcal{P}_{\!\gamma}h(x))$.
 \begin{lemma}\label{prox-vartheta}
  Fix any $\gamma>0$. Then, for any $X\in\mathbb{R}^{n\times r}$, 
  $\mathcal{P}_{\!\gamma}\vartheta(X)=\min(X\!+\!\gamma,\max(X,0))$ and 
  \[
    e_{\gamma}\vartheta(X)=\frac{1}{2\gamma}\|\mathcal{P}_{\!\gamma}\vartheta(X)-X\|_F^2
    +\vartheta(\mathcal{P}_{\!\gamma}\vartheta(X))
    =\frac{1}{2\gamma}\sum_{i\in J_1(X)}X_{ij}^2
    +\sum_{i\in J_2(X)}[-X_{ij}-{\gamma}/{2}],
  \]  
  where $J_1(X):=\{(i,j)\in[n]\times[r]\ |\ X_{ij}\in[-\gamma,0]\}$ 
  and $J_2(X):=\{(i,j)\in[n]\times[r]\ |\ X_{ij}<-\gamma\}$. 
 \end{lemma}
 \subsection{Retractions on manifold}\label{sec2.1}
 \begin{definition}\label{Def-retr}
   A retraction on a manifold $\mathcal{M}$ is a $\mathcal{C}^{\infty}$-mapping $R$ from
   the tangent bundle ${\rm T}\!\mathcal{M}:=\bigcup_{x\in\mathcal{M}}{\rm T}_{x}\mathcal{M}$
   onto $\mathcal{M}$, and for any $x\in\mathcal{M}$, the restriction of $R$ to 
   ${\rm T}_{x}\mathcal{M}$, denoted by $R_x$, satisfies the following two conditions:
   \begin{enumerate}[1.]
     \item $R_{x}(0_{x})=x$, where $0_{x}$ denotes
                 the zero element of ${\rm T}_{x}\mathcal{M}$;

     \item the differential of $R_x$ at $0_{x}$, $DR_{x}(0_{x})$, is the identity mapping on ${\rm T}_{x}\mathcal{M}$.
 \end{enumerate}
 \end{definition}
 \begin{remark}\label{remark-retraction}
  For the submanifold $M$ of $\mathbb{R}^{n\times r}$,
  Item 2 in Definition \ref{Def-retr} is equivalent to saying that
  \[
   \lim_{{{\rm T}_{\!x}M}\ni\xi\to 0}\frac{\|R_{x}(\xi)-(x+\xi)\|_F}{\|\xi\|_F}=0.
  \]
  The common retractions for the Stiefel manifold include the exponential mapping,
 the Cayley transformation, the polar decomposition and the QR decomposition.
 For the detailed discussion on their computation complexity, the reader is 
 referred to \cite[Section 3]{Chen20}.
 \end{remark}
 \begin{lemma}\label{property1-manifold}
  (see \cite[Eq. (B.3)$\,\&$(B.4)]{Boumal18})
  Let $\mathcal{M}$ be a compact submanifold of a Euclidean space $\mathcal{E}$, 
  and let $R\!:{\rm T}\mathcal{M}\to \mathcal{M}$ be a retraction on $\mathcal{M}$. 
  Then, there exist constants $c_1>0$ and $c_2>0$ such that for any $X\in\mathcal{M}$ 
  and $\xi\in {\rm T}_{\!X}\mathcal{M}$,
  \begin{equation}\label{key-ineq1}
   \|R_{X}(\xi)-X\|_F\le c_1\|\xi\|_F
   \ \ {\rm and}\ \ \|R_{X}(\xi)-(X+\xi)\|_F\le c_2\|\xi\|_F^2.
  \end{equation}
 \end{lemma}
 \subsection{Subdifferentials of nonsmooth functions on manifold}\label{sec2.1}
   
 First of all, we recall from \cite{Ledyaev07} the subdifferentials of a nonsmooth 
 function on $(\mathcal{M},\langle\cdot,\cdot\rangle_{g})$.
 \begin{definition}\label{Regsubdiff}
  Let $h\!:\mathcal{M}\to\overline{\mathbb{R}}$ be a function
  with ${\rm dom}h\!:=\!\big\{x\in\mathcal{M}\ |\ h(x)<+\infty\big\}\ne\emptyset$.
  We define the Fr\'{e}chet-subdifferential of $h$ at a point $x\in\mathcal{M}$ by
  \[
    \widehat{\partial}h(x)
    \!:=\!\left\{\begin{array}{cl}
    \!\big\{d\theta_{x}\,|\,\exists\theta\in\mathcal{C}^1(\mathcal{M})\ 
    {\rm and}\ h\!-\!\theta\ {\rm attains\ a\ local\ minimum\ at}\ x\big\}
    &{\rm if}\ x\in{\rm dom}h,\\
    \emptyset&{\rm if}\  x\notin{\rm dom}h,
    \end{array}\right.
  \]
  where $d\theta_{x}\in{\rm T}_{\!x}^*\mathcal{M}$, 
  the dual space of ${\rm T}_{\!x}\mathcal{M}$,
  is defined by $d\theta_{x}(\xi)=\langle{\rm grad}\theta(x),\xi\rangle_{g}$ 
  for any $\xi\!\in\!{\rm T}_{\!x}\mathcal{M}$;
  and define its (limiting) subdifferential and singular subdifferential 
  at $x\!\in\!\mathcal{M}$ by
  \begin{align*}
   \partial h(x)
    =\Big\{\lim_{k\to\infty} v^k\,|\,\exists v^k\in\widehat{\partial}h(x^k)
            \ {\rm with}\ (x^k,h(x^k))\to(x,h(x))\Big\},\\
   \partial^{\infty}h(x)
   =\Big\{\lim_{k\to\infty} \lambda^kv^k\,|\,\exists v^k\in\widehat{\partial}h(x^k)
   \ {\rm with}\ (x^k,h(x^k))\to(x,h(x))\ {\rm and}\ \lambda^k\downarrow 0\Big\}.
  \end{align*}
 \end{definition}
 \begin{remark}\label{Remark-subdiff}
  {\bf(a)} By \cite[Remark 3.2]{Ledyaev07}, the support function $\theta$ in the definition of 
  Fr\'{e}chet subdifferential need only be $\mathcal{C}^1$ in a neighborhood of $x$. 
  In addition, by making a translation, the function $\theta$ can be required to satisfy  
  $h\!-\theta\ge\!0$ in a neighborhood of $x$. Thus, if $\mathcal{M}$ is a Euclidean space, 
  the above definitions of subdifferentials coincide with those of 
  \cite[Definition 8.3]{RW98}. For a proper 
  $h\!:\mathcal{M}\to\overline{\mathbb{R}}$ and a point $x\in{\rm dom}h$, 
  obviously, $\widehat{\partial}h(x)\subseteq\partial h(x)$, 
  and the multifunction $\partial h\!:\mathcal{M}
  \rightrightarrows{\rm T}\!\mathcal{M}$ is closed. 
  
  \noindent
  {\bf(b)} For a proper lsc function $h\!:\mathcal{E}\to\overline{\mathbb{R}}$ 
  and a compact submanifold $\mathcal{M}$ of $\mathcal{E}$, 
  by Definition \ref{Regsubdiff}, if $x^*$ is a local minimizer of problem
  $\min_{x\in\mathcal{M}}h(x)$, then 
  $0\in\widehat{\partial}h_{|\mathcal{M}}(x^*)\subseteq\partial h_{|\mathcal{M}}(x^*)$. 
  In view of this, we call $x$ a stationary point of problem 
  $\min_{x\in \mathcal{M}}h(x)$ if $0\in\partial h_{|\mathcal{M}}(x)$. 
 \end{remark}
 \begin{lemma}\label{relation-subdiff}
  For a proper $h\!:\mathcal{E}\to\overline{\mathbb{R}}$ 
  and a compact submanifold $\mathcal{M}$ of $\mathcal{E}$, 
  let $\widetilde{h}:=h+\delta_{\mathcal{M}}$.
  Consider a point $\overline{x}\in{\rm dom}h\cap\mathcal{M}$. 
  Then, $\partial^{s}h_{|\mathcal{M}}(\overline{x})
  ={\rm Proj}_{{\rm T}_{\overline{x}}\mathcal{M}}(\partial^{s}\widetilde{h}(\overline{x}))
  \subseteq\partial^{s}\widetilde{h}(\overline{x})$, 
  where $\partial^{s}$ is one of $\widehat{\partial},\partial$ or $\partial^{\infty}$, 
  so when $h$ is locally Lipschitz at $\overline{x}$, 
  the conclusion also holds for the Clarke subdifferential.
 \end{lemma}
 \begin{proof}
  Pick any $v\in\widehat{\partial}h_{|\mathcal{M}}(\overline{x})$. 
  There exists $\theta\!\in\mathcal{C}^1(\mathcal{M})$ such that 
  $v={\rm grad}\theta(\overline{x})$ 
  and $\overline{x}$ is a local minimizer of $h_{|\mathcal{M}}\!-\theta$ 
  with $h_{|\mathcal{M}}(\overline{x})\!=\theta(\overline{x})$.
  Since $\theta\in\mathcal{C}^1(\mathcal{M})$, there exists a representative function
  $\widehat{\theta}\!:\mathcal{E}\to\mathbb{R}$ that is $\mathcal{C}^1$-smooth
  around $\overline{x}$ with $\widehat{\theta}_{|\mathcal{M}}=\theta$ 
  locally around $\overline{x}$. From \cite[Section 3.6.1]{Absil08},
  \(
    {\rm Proj}_{{\rm T}_{\overline{x}}\mathcal{M}}(\nabla\widehat{\theta}(\overline{x}))
    ={\rm grad}\theta(\overline{x})=v.
  \)
  Then, on some neighborhood of $\overline{x}$, the function
  $\widehat{\theta}$ is smooth and $\widehat{\theta}\le\widetilde{h}_{|\mathcal{M}}
  \!=h_{|\mathcal{M}}+\delta_{\mathcal{M}}=\widetilde{h}$ with 
  $\widehat{\theta}(\overline{x})=\widetilde{h}(\overline{x})$.
  By \cite[Proposition 8.5]{RW98}, 
  $\nabla\widehat{\theta}(\overline{x})\in\widehat{\partial}\widetilde{h}(\overline{x})$,
  which means that $v\in{\rm Proj}_{{\rm T}_{\overline{x}}\mathcal{M}}
  (\widehat{\partial}\widetilde{h}(\overline{x}))$, and hence 
  $\widehat{\partial}h_{|\mathcal{M}}(\overline{x})\subseteq
  {\rm Proj}_{{\rm T}_{\overline{x}}\mathcal{M}}
  (\widehat{\partial}\widetilde{h}(\overline{x}))$.
  To establish the converse inclusion, pick any
  $v\in{\rm Proj}_{{\rm T}_{\overline{x}}\mathcal{M}}
  (\widehat{\partial}\widetilde{h}(\overline{x}))$.
  Then, there exists $\xi\in\widehat{\partial}\widetilde{h}(\overline{x})$ 
  such that $v={\rm Proj}_{{\rm T}_{\overline{x}}\mathcal{M}}(\xi)$. 
  Since $\xi\in\widehat{\partial}\widetilde{h}(\overline{x})$,
  by \cite[Proposition 8.5]{RW98}, on some neighborhood $\mathcal{U}$ of $\overline{x}$,
  there is a smooth function $\theta\le\widetilde{h}$ with $\nabla\theta(\overline{x})=\xi$.
  Clearly, $\theta_{|\mathcal{M}}$ is a smooth function on $\mathcal{U}\cap\mathcal{M}$ 
  such that $\overline{x}$ is a local minimizer of 
  $h_{|\mathcal{M}}-\theta_{|\mathcal{M}}$.
  By Item 1 of Remark \ref{Remark-subdiff}, 
  $v={\rm grad}\theta_{|\mathcal{M}}(\overline{x})
  \in\widehat{\partial}h_{|\mathcal{M}}(\overline{x})$.
  The converse inclusion follows. The equality holds with $\partial^{s}=\widehat{\partial}$. 
  For the inclusion, pick any 
  $v\in{\rm Proj}_{{\rm T}_{\overline{x}}\mathcal{M}}
  (\widehat{\partial}\widetilde{h}(\overline{x}))$.
  Then, there exists $\xi\in\widehat{\partial}\widetilde{h}(\overline{x})$ such that
  $v={\rm Proj}_{{\rm T}_{\overline{x}}\mathcal{M}}(\xi)$. Consequently,
  \begin{align*}
  &\liminf_{\overline{x}\ne x\to\overline{x}}
   \frac{\widetilde{h}(x)-\widetilde{h}(\overline{x})-\langle v,x-\overline{x}\rangle}
   {\|x-\overline{x}\|}
   =\liminf_{\overline{x}\ne x\xrightarrow[\mathcal{M}]{}\overline{x}}
     \frac{\widetilde{h}(x)-\widetilde{h}(\overline{x})-\langle v,x-\overline{x}\rangle}
     {\|x-\overline{x}\|}\\
   &=\liminf_{\overline{x}\ne x\xrightarrow[\mathcal{M}]{}\overline{x}}
     \frac{\widetilde{h}(x)-\widetilde{h}(\overline{x})-\langle\xi,x-\overline{x}\rangle
     +\langle{\rm Proj}_{{\rm N}_{\overline{x}}\mathcal{M}}(\xi),x-\overline{x}\rangle}
     {\|x-\overline{x}\|}\\
   &\ge\liminf_{\overline{x}\ne x\xrightarrow[\mathcal{M}]{}\overline{x}}
    \frac{\widetilde{h}(x)-\widetilde{h}(\overline{x})-\langle\xi,x-\overline{x}\rangle}
    {\|x-\overline{x}\|}
    -\limsup_{\overline{x}\ne x\xrightarrow[\mathcal{M}]{}\overline{x}}
    \frac{\langle-{\rm Proj}_{{\rm N}_{\overline{x}}\mathcal{M}}(\xi),x-\overline{x}\rangle}
    {\|x-\overline{x}\|}\\
   &\ge\liminf_{\overline{x}\ne x\xrightarrow[\mathcal{M}]{}\overline{x}}
   \frac{\widetilde{h}(x)-\widetilde{h}(\overline{x})-\langle\xi,x-\overline{x}\rangle}
   {\|x-\overline{x}\|}\ge 0
  \end{align*}
  where the second inequality holds by $-{\rm Proj}_{{\rm N}_{\overline{x}}\mathcal{M}}(\xi)
  \in{\rm N}_{\overline{x}}\mathcal{M}$ and the last is by 
  $\xi\in\widehat{\partial}\widetilde{h}(\overline{x})$.
  By the definition of the Fr\'{e}chet subgradient, 
  $v\in\widehat{\partial}\widetilde{h}(\overline{x})$.
  The inclusion holds with $\partial^{s}=\widehat{\partial}$.

  To prove that $\partial h_{|\mathcal{M}}(\overline{x})
  ={\rm Proj}_{{\rm T}_{\overline{x}}\mathcal{M}}(\partial\widetilde{h}(\overline{x}))$,
  pick any $v\in\!\partial h_{|\mathcal{M}}(\overline{x})$. 
  Clearly, $v\in{\rm T}_{\overline{x}}\mathcal{M}$.
  Also, there exist $(x^k,h_{|\mathcal{M}}(x^k))
  \to(\overline{x},h_{|\mathcal{M}}(\overline{x}))$
  and $v^k\!\in\widehat{\partial} h_{|\mathcal{M}}(x^k)$ for each $k$ such that $v^k\to v$.
  Since the equality and inclusion hold with $\partial^{s}=\widehat{\partial}$ 
  for all $x\in{\rm dom}h\cap\mathcal{M}$, for each $k$ it holds that
  $v^k\in{\rm Proj}_{{\rm T}_{x^k}\mathcal{M}}(\widehat{\partial}\widetilde{h}(x^k))
  \subseteq\widehat{\partial}\widetilde{h}(x^k)$.
  Since $h_{|\mathcal{M}}(x^k)\to h_{|\mathcal{M}}(\overline{x})$ implies
  $(x^k,\widetilde{h}(x^k))\to(\overline{x},\widetilde{h}(\overline{x}))$,
  we have $v\in \partial\widetilde{h}(\overline{x})$. Together with
  $v\in{\rm T}_{\overline{x}}\mathcal{M}$, 
  $v={\rm Proj}_{{\rm T}_{\overline{x}}\mathcal{M}}(v)
  \in{\rm Proj}_{{\rm T}_{\overline{x}}\mathcal{M}}(\partial\widetilde{h}(\overline{x}))$
  and $\partial h_{|_{\mathcal{M}}}(\overline{x})
  \subseteq{\rm Proj}_{{\rm T}_{\overline{x}}\mathcal{M}}(\partial\widetilde{h}(\overline{x}))$.
  For the converse inclusion, pick any $v\in{\rm Proj}_{{\rm T}_{\overline{x}}\mathcal{M}}(\partial\widetilde{h}(\overline{x}))$. Then,
  there exists $\xi\in\partial\widetilde{h}(\overline{x})$ such that
  $v={\rm Proj}_{{\rm T}_{\overline{x}}\mathcal{M}}(\xi)$. 
  Since $\xi\in\partial\widetilde{h}(\overline{x})$,
  there exists $(x^k,\widetilde{h}(x^k))\to(\overline{x},\widetilde{h}(\overline{x}))$
  with $\xi^k\in\widehat{\partial}\widetilde{h}(x^k)$ for every $k$ such that $\xi^k\to\xi$.
  Since $\widetilde{h}(x^k)\to\widetilde{h}(\overline{x})$ implies that
  $\{x^k\}\subseteq\mathcal{M}$ and $h(x^k)\to h(\overline{x})$, we have
  $h_{|\mathcal{M}}(x^k)\to h_{|\mathcal{M}}(\overline{x})$. 
  Let $v^k:={\rm Proj}_{{\rm T}_{x^k}\mathcal{M}}(\xi^k)$.
  From \eqref{Proj-St}, $v^k\to v$. From the equality with 
  $\partial^{s}=\widehat{\partial}$, for each $k$, 
  $v^k\in\widehat{\partial}h_{|\mathcal{M}}(x^k)$.
  Then, $v\in\partial h_{|\mathcal{M}}(\overline{x})$ and
   $\partial h_{|\mathcal{M}}(\overline{x})\supseteq
   {\rm Proj}_{{\rm T}_{\overline{x}}\mathcal{M}}(\partial\widetilde{h}(\overline{x}))$. 
  Thus, the equality holds with $\partial^{s}=\partial$. For the inclusion,
  pick any $v\in{\rm Proj}_{{\rm T}_{\overline{x}}\mathcal{M}}
  (\partial\widetilde{h}(\overline{x}))$.
  There exists $\xi\in\partial\widetilde{h}(\overline{x})$ such that
   $v={\rm Proj}_{{\rm T}_{\overline{x}}\mathcal{M}}(\xi)$.
   Since $\xi\in\partial\widetilde{h}(\overline{x})$, there exist
   $(x^k,\widetilde{h}(x^k))\to(\overline{x},\widetilde{h}(\overline{x}))$
   and $\xi^k\in\widehat{\partial}\widetilde{h}(x^k)$ for each $k$ such that
   $\xi^k\to\xi$ as $k\to\infty$. From the inclusion with 
   $\partial^{s}=\widehat{\partial}$,
   ${\rm Proj}_{{\rm T}_{x^k}\mathcal{M}}(\xi^k)\in\widehat{\partial}\widetilde{h}(\xi^k)$.
   By the outer semicontinuity of $\partial\widetilde{h}$ and \eqref{Proj-St},
   $v={\rm Proj}_{{\rm T}_{\overline{x}}\mathcal{M}}(\xi)
   \in\partial\widetilde{h}(\overline{x})$.
   The inclusion holds with $\partial^{s}=\partial$. Using the similar arguments 
   yields the equality and inclusion hold with $\partial^{s}=\partial^{\infty}$.
 \end{proof}
  
  Lemma \ref{relation-subdiff} discloses the relation between the subdifferentials 
  of $h$ and those of its restriction on $\mathcal{M}$. By using this lemma, 
  the smoothness of $f$ and the convexity of $\vartheta$, we have the following result.
  \begin{corollary}\label{subdiff-EPfun}
   Fix any $\rho>0$ and $\gamma>0$. For any $X\!\in{\rm St}(n,r)$, 
   the following relations hold
   \begin{subequations}
   \begin{align*}
    \widehat{\partial}(\Theta_{\rho,\gamma})_{|M}(X)
    =\partial(\Theta_{\rho,\gamma})_{|M}(X)=\{{\rm grad}\Theta_{\rho,\gamma}(X)\}
    \!=\big\{{\rm Proj}_{{\rm T}_{\!X}M}(\nabla\Theta_{\rho,\gamma}(X))\big\},\\
     \widehat{\partial}(\Phi_{\!\rho})_{|M}(X)
     =\partial(\Phi_{\!\rho})_{|M}(X)
     \!={\rm Proj}_{{\rm T}_{\!X}M}\big[\nabla\!f(X)+\rho\partial\vartheta(X)\big].
     \qquad\quad
   \end{align*}
   \end{subequations}
  \end{corollary}
 \subsection{KL property on Riemannian manifold}\label{sec2.3}
 \begin{definition}\label{KL-Def1}
  (see \cite[Definition 3.3]{Neto13}) Let $h\!:\mathcal{M}\!\to\overline{\mathbb{R}}$
  be a proper function. The function $h$ is said to have the KL property at
  $\overline{x}\in{\rm dom}\,\partial h$ if there exist $\eta\in(0,+\infty]$,
  a neighborhood $U$ of $\overline{x}$, and a continuous concave function
  $\varphi\!:[0,\eta)\to\mathbb{R}_{+}$ such that
  \begin{enumerate}[1.]
    \item  $\varphi(0)=0$, $\varphi$ is $\mathcal{C}^1$ on $(0,\eta)$,
                and $\varphi'(s)>0$ for all $s\in(0,\eta)$;

    \item for all $x\in U\cap\big\{x\in\mathcal{M}\ |\ 
           h(\overline{x})<h(x)<h(\overline{x})+\eta\big\}$, the KL inequality holds:
           \[
            \varphi'(h(x)-h(\overline{x}))d(0,\partial h(x))\ge 1,
           \]
         where $d(0,\partial h(x)):=\inf_{v\in\partial h(x)}\|v\|_g$
         and $\|\cdot\|_g$ is the norm induced by $\langle\cdot,\cdot\rangle_g$.
  \end{enumerate}
  If $h$ has the KL property at every point of ${\rm dom}\,\partial h$,
  then it is called a KL function.
 \end{definition}
 
 For a function $h\!:\mathcal{E}\!\to\overline{\mathbb{R}}$ and 
 a compact submanifold $\mathcal{M}$ of $\mathcal{E}$,
 the following lemma discloses the relation between the KL property of 
 $h+\delta_{\mathcal{M}}$ and of its restriction $h_{|\mathcal{M}}$. 
 Although this result is not used in this paper, 
 as will be shown in Remark \ref{Remark-KL}, it provides a more convenient
 criterion to identify the KL property of a function defined on $\mathcal{M}$ 
 than \cite[Theorem 6]{Huang21} does.
 
 \begin{lemma}\label{KL-link}
  For a proper $h\!:\mathcal{E}\to\overline{\mathbb{R}}$ 
  and a compact submanifold $\mathcal{M}$ of $\mathcal{E}$, 
  if $h+\delta_{\mathcal{M}}$ is a KL function, then $h_{|{\mathcal{M}}}$ 
  is a KL function; conversely, if $h_{|{\mathcal{M}}}$ is a KL function, 
  then so is $h+\delta_{\mathcal{M}}$.
 \end{lemma}
 \begin{proof}
  Write $\widetilde{h}:=h+\delta_{\mathcal{M}}$. Suppose $\widetilde{h}$ is a KL function.
  Fix any $\overline{x}\in{\rm dom}\partial h_{|\mathcal{M}}$. 
  If $0\notin\partial h_{|\mathcal{M}}(\overline{x})$,
  then $h_{|\mathcal{M}}$ has the KL property at $\overline{x}$ by \cite[Lemma 3.1]{Neto13}.
  So, it suffices to consider the case $0\in\partial h_{|\mathcal{M}}(\overline{x})$.
  Since $\widetilde{h}$ has the KL property at $\overline{x}$,
  by \cite[Definition 3.1]{Attouch10}, there exist $\eta\in(0,+\infty]$,
  a neighborhood $U$ of $\overline{x}$,
  and a continuous concave function $\varphi\!:[0,\eta)\to\mathbb{R}_{+}$
  satisfying Definition \ref{KL-Def1} (i) such that for all
  $x\in U\cap\{z\in{\rm dom}\widetilde{h}\,|\, 
  \widetilde{h}(\overline{x})<\widetilde{h}(z)<\widetilde{h}(\overline{x})+\eta\}$,
  \begin{equation}\label{wth-ineq}
    1\le\varphi'(\widetilde{h}(x)-\widetilde{h}(\overline{x}))
    {\rm dist}(0,\partial\widetilde{h}(x))
    =\varphi'(h(x)-h(\overline{x})){\rm dist}(0,\partial\widetilde{h}(x)).
  \end{equation} 
  Now fix any $x\in U\cap\big\{z\in\mathcal{M}\ |\ h_{|\mathcal{M}}(\overline{x})
  <h_{|\mathcal{M}}(z)\!<h_{|\mathcal{M}}(\overline{x})+\eta\big\}$.
  If $\partial h_{|\mathcal{M}}(x)=\emptyset$, clearly,
  \begin{equation}\label{aim-ineq}
     \varphi'(h_{|\mathcal{M}}(x)-h_{|\mathcal{M}}(\overline{x}))
     {\rm dist}(0,\partial h_{|\mathcal{M}}(x))\ge 1.
  \end{equation}
  We next consider the case that $\partial h_{|\mathcal{M}}(x)\ne\emptyset$.
  By Lemma \ref{relation-subdiff}, $\partial\widetilde{h}(x)\ne\emptyset$ and moreover,
  \[
   {\rm dist}(0,\partial h_{|\mathcal{M}}(x))
   ={\rm dist}(0,{\rm Proj}_{{\rm T}_{x}\mathcal{M}}(\partial\widetilde{h}(x)))
   =\!\inf_{v\in{\rm Proj}_{{\rm T}_{x}\mathcal{M}}(\partial\widetilde{h}(x))}\|v\|
   \ge \inf_{v\in\partial\widetilde{h}(x)}\|v\|.
  \]
  Together with \eqref{wth-ineq}, it follows that
  $\varphi'(h(x)-h(\overline{x})){\rm dist}(0,\partial h_{|\mathcal{M}}(x))\ge 1$, 
  i.e., the inequality \eqref{aim-ineq} also holds.
  This shows that $h_{|\mathcal{M}}$ satisfies the KL property at $\overline{x}$.
  From the arbitrariness of $\overline{x}\in{\rm dom}\partial h_{|\mathcal{M}}$,
  we conclude that $h_{|\mathcal{M}}$ is a KL function.

  Conversely, suppose that $h_{|\mathcal{M}}$ is a KL function.
  Fix any $\overline{x}\in{\rm dom}\partial\widetilde{h}$.
  If $0\notin\partial\widetilde{h}(\overline{x})$, then $\widetilde{h}$ has
  the KL property at $\overline{x}$ by \cite[Lemma 2.1]{Attouch10}. So,
  it suffices to consider the case that $0\in\partial\widetilde{h}(\overline{x})$.
  Since the function $h_{|\mathcal{M}}$ admits the KL property at $\overline{x}$,
  there exist $\eta\in(0,+\infty]$, a neighborhood $U$ of $\overline{x}$,
  and a continuous concave function $\varphi\!:[0,\eta)\to\mathbb{R}_{+}$
  satisfying Definition \ref{KL-Def1} (i) such that for all
  $x\in U\cap\big\{z\in\mathcal{M}\ |\ h_{|\mathcal{M}}(\overline{x})
  <h_{|\mathcal{M}}(z)<h_{|\mathcal{M}}(\overline{x})+\eta\big\}$,
  \[
    1\le\varphi'(h_{|\mathcal{M}}(x)\!-\!h_{|\mathcal{M}}(\overline{x}))
    {\rm dist}(0,\partial h_{|\mathcal{M}}(x))
    =\varphi'(\widetilde{h}(x)-\widetilde{h}(\overline{x}))
    {\rm dist}(0,{\rm Proj}_{{\rm T}_{x}\mathcal{M}}(\partial\widetilde{h}(x)))
  \]
  which the equality is due to Lemma \ref{relation-subdiff} and
  the definition of $\widetilde{h}$. In addition, by Lemma \ref{relation-subdiff},
  ${\rm Proj}_{{\rm T}_{x}\mathcal{M}}(\partial\widetilde{h}(x))
  \subseteq\partial\widetilde{h}(x)$, which implies that 
  $0\in{\rm Proj}_{{\rm N}_{x}\mathcal{M}}(\partial\widetilde{h}(x))$. 
  Then,
  \begin{align*}
   {\rm dist}(0,\partial\widetilde{h}(x))
   &=\inf_{v\in\partial\widetilde{h}(x)}\|v\|
   \ge\inf_{v\in{\rm Proj}_{{\rm T}_{x}\mathcal{M}}(\partial\widetilde{h}(x))
   +{\rm Proj}_{{\rm N}_{x}\mathcal{M}}(\partial\widetilde{h}(x))}\|v\|\\
   &=\inf_{v^1\in{\rm Proj}_{{\rm T}_{x}\mathcal{M}}(\partial\widetilde{h}(x))}\|v^1\|
     +\inf_{v^2\in{\rm Proj}_{{\rm N}_{x}\mathcal{M}}(\partial\widetilde{h}(x))}\|v^2\|\\
  &=\inf_{v^1\in{\rm Proj}_{{\rm T}_{x}\mathcal{M}}(\partial\widetilde{h}(x))}\|v^1\|
   ={\rm dist}(0,{\rm Proj}_{{\rm T}_{x}\mathcal{M}}(\partial\widetilde{h}(x))).
  \end{align*}
  From the last two inequalities, $\widetilde{h}$ has the KL property at $\overline{x}$
  and $\widetilde{h}$ is a KL function.
 \end{proof}
 \begin{remark}\label{Remark-KL}
  Lemma \ref{KL-link} provides a convenient criterion to identify the KL property
  of a function $h\!:\mathcal{M}\to\overline{\mathbb{R}}$ in terms of that of
  its extension 
  \(
   \widetilde{h}(x)=\left\{\begin{array}{cl}
                        h(x) &{\rm if}\ x\in\mathcal{M},\\
                       +\infty &{\rm otherwise}.
                     \end{array}\right.
  \)
 One can easily check the KL property of $\widetilde{h}$ through 
 a collection of criteria included in \cite[Section 4]{Attouch10}.
 From the last part of Lemma \ref{relation-subdiff}
 and the proof of Lemma \ref{KL-link}, the conclusions of Lemma \ref{KL-link} still hold
 if the KL property is defined by the Clarke subdifferential. Recently, Huang and Wei
 \cite{Huang21} provided a criterion to identify such a KL property of $h$ 
 at $x\in\mathcal{M}$, but their criterion involves verifying that the mapping 
 $\phi_{x}^{-1}:\phi_{x}(\mathcal{U})\subset\mathbb{R}^{d}\to\mathbb{R}^{n}$ 
 is semialgebraic, where $(\mathcal{U},\phi_{x})$ is a chart covering $x$. 
 As shown in the proof of \cite[Lemma 8 $\&$9]{Huang21}, this is not an easy task. 
 \end{remark}
 \section{Error bound and global exact penalties}\label{sec3}

  To derive a local error bound for the constraint system 
  $X\in\mathcal{S}_{+}^{n,r}$, for every $X\in\mathbb{R}^{n\times r}$, we denote by $X_r$  
  the submatrix consisting of its first $r$ rows, and by $\sigma(X)$ 
  its singular value vector with entries arranged in a nonincreasing order. First, 
  we show that the distance between a nonnegative matrix and a rectangular identity
  matrix is locally upper bounded by the distance between their singular value vectors.   
 \begin{lemma}\label{lemma1-NonStnr}
  For each $c>0$, define $\mathcal{L}_{c}:=\{X\in\mathbb{R}^{n\times r}\ |\
  X_{r}\ge 0,\,\|X-I_{n\times r}\|_F\le c\|X_r-I_r\|_F\}$. Then, for each $c>0$, 
  there exists $\delta>0$ such that for all $X\in\mathbb{B}(I_{n\times r},\delta)\cap\mathcal{L}_{c}$, 
  $\|X_{r}-I_{r}\|_F\le 2.1\sqrt{r}\|\sigma(X)-e\|$.
 \end{lemma}
 \begin{proof}
  Suppose on the contrary that the conclusion does not hold. There exist $\overline{c}>0$ 
  and a sequence $\{X^k\}_{k\in\mathbb{N}}\subset\mathcal{L}_{\overline{c}}$ with 
  $X^k\to I_{n\times r}$ as $k\to\infty$ such that for each $k\in\mathbb{N}$,
  $\|\sigma(X^k)-e\|<\frac{1}{2.1\sqrt{r}}\|X_{r}^k-I_r\|_F$.
  For each $k$, write $H^k:=X^k-I_{n\times r}$. From \cite[Proposition 6]{DingST14}, 
  for all sufficiently large $k$,
  \[
   \sigma_i(X^k)-1=\lambda_i\Big(\frac{H_r^k+(H_r^k)\tp }{2}\Big)+O(\|H^k\|_F^2)
   \quad{\rm for}\ \ i=1,2,\ldots,r,
  \]
  where $\lambda_i(Z)$ denotes the $i$th largest eigenvalue of a $r\times r$ real symmetric
  matrix $Z$. By taking into account that $\|\sigma(X^k)-e\|\ge\frac{1}{\sqrt{r}}\sum_{i=1}^r|\sigma_i(X^k)-1|$, 
  for all sufficiently large $k$, it holds that
  \begin{align*}
  \|\sigma(X^k)-e\|
  &\ge \frac{1}{\sqrt{r}}\sum_{i=1}^r\Big|\lambda_i\Big(\frac{H_r^k+(H_r^k)\tp }{2}\Big)\Big|+O(\|H^k\|_F^2)\\
  &=\frac{1}{2\sqrt{r}}\big\|H_r^k+(H_r^k)\tp \big\|_*+O(\|H^k\|_F^2)\\
  &\ge\frac{1}{2\sqrt{r}}\big\|H_r^k+(H_r^k)\tp \big\|_F+O(\|H^k\|_F^2)
  >\frac{1}{2.1\sqrt{r}}\|H^k_r\big\|_F,
  \end{align*}
  where the last inequality is using $X_{r}^k\ge 0,H_{r}^k=X_{r}^k-I_{r}$
  and $\|H^k\|_F\le \overline{c}\|H_{r}^k\|_F$. This yields a contradiction to 
  the inequality $\|\sigma(X^k)-e\|<\frac{1}{2.1\sqrt{r}}\|X_{r}^k-I_r\|_F$. 
  The conclusion then follows.
 \end{proof}
 
 By Lemma \ref{lemma1-NonStnr}, we can show that the distance from a nonnegative matrix 
 to $\mathcal{S}_{+}^{n,r}$ is locally upper bounded by the distance between 
 its singular value vector and that of an identity matrix. 
 \begin{lemma}\label{lemma2-NonStnr}
  Consider any $\overline{X}\in\mathcal{S}_{+}^{n,r}$. If the matrix $\overline{X}$ has no zero rows when $n>r>1$, then there exists $\delta>0$ 
  such that for all $X\in\mathbb{B}(\overline{X},\delta)\cap\mathbb{R}_{+}^{n\times r}$,
  \begin{equation*}
   {\rm dist}\big(X,\mathcal{S}_{+}^{n,r}\big)\le \kappa\|\sigma(X)-e\|
   \ \ {\rm with}\ \ 
   \kappa:=\left\{\begin{array}{cl}
      2.1\sqrt{n}&{\rm if}\ n=r,\\
      1 &{\rm if}\ n>r=1,\\
     \frac{2.1\sqrt{r}[1+3r(n\!-\!r)]}{\overline{X}_{\!i^{*}\!j^{*}}}&{\rm if}\ n>r>1
     \end{array}\right. 
  \end{equation*}
  where $\overline{X}_{\!i^{*}\!j^{*}}$ is the smallest nonzero entry of $\overline{X}$.
 \end{lemma}
 \begin{proof}
  We first consider the case that $n=r$. By invoking Lemma \ref{lemma1-NonStnr} for $c=1$, 
  there exists $\delta>0$ such that for all $Z\in\mathbb{B}(I_n,\delta)\cap
  \mathbb{R}_{+}^{n\times n}$, $\|Z-I_n\|_F\le 2.1\sqrt{n}\|\sigma(Z)-e\|$.
  Pick any $X\!\in\mathbb{B}(\overline{X},{\delta})\cap\mathbb{R}_{+}^{n\times n}$. 
  Then, by noting that $\|\overline{X}\tp X\!-I_n\|_F=\!\|X\!-\!\overline{X}\|_F\le\delta$, 
  we have $\overline{X}\tp \!X\in\mathbb{B}(I_n,\delta)\cap\mathbb{R}_{+}^{n\times n}$.
  Consequently,
  \[
   {\rm dist}(X, \mathcal{S}_{+}^{n,n})
   \le\|\overline{X}-X\|_F=\|\overline{X}\tp X-I_n\|_F\\
   \le 2.1\sqrt{n}\|\sigma(\overline{X}\tp X)-e\|=2.1\sqrt{n}\|\sigma(X)-e\|.
  \]
  For the case that $n>r=1$, because $\mathcal{S}_{+}^{n,r}$ reduces to the nonnegative sphere in $\mathbb{R}^n$, from $\overline{X}\in\mathcal{S}_{+}^{n,r}$ and the continuity, there exists $\delta>0$ such that for all $X\in\mathbb{B}(\overline{X},\delta)$, $\|X\|>0$. 
  Fix any $X\in\mathbb{B}(\overline{X},\delta)\cap\mathbb{R}_{+}^{n\times r}$. 
  Then, ${\rm dist}(X, \mathcal{S}_{+}^{n,r})=\|X-\frac{X}{\|X\|}\|=|\|X\|-1|$ 
  and $\sigma(X)=\|X\|$. The desired inequality holds with $\kappa=1$. 
  
  Next we focus on the case that $n>r>1$. Because $\overline{X}\in\mathcal{S}_{+}^{n,r}$ has no zero rows, each row of $\overline{X}$ has only one nonzero entry and each column of $\overline{X}$ has at least one nonzero entry. For each $k\in[n]$, 
  write $l_k:=\arg\max_{j\in[r]}\overline{X}_{kj}$. Clearly, for each $k\in[n]$, $\overline{X}_{\!kl_k}>0$. By the continuity, there exists $\delta_1\in(0,1/3)$ such that for all $X\!\in\mathbb{B}(\overline{X},\delta_1)$ and all
  $k\in[n]$, $X_{\!kl_k}>0$ and $\arg\max_{j\in[r]}X_{kj}=l_k$.
  Take  $c=\frac{\overline{X}_{i^{*}j^{*}}+3(n-r)r}{\overline{X}_{\!i^{*}\!j^{*}}}$ 
  and let $\mathcal{L}_{c}$ be defined as in Lemma \ref{lemma1-NonStnr}. Then, 
  there exists $\delta_2>0$ such that for all 
  $X'\in\mathbb{B}(I_{n\times r},\delta_2)\cap\mathcal{L}_{c}$,
  \begin{equation}\label{Zr-dist}
    \|X_{r}'-I_{r}\|_F\le 2.1\sqrt{r}\|\sigma(X')-e\|.
  \end{equation}
  Set $\delta=\min\big\{\frac{1}{2}\overline{X}_{\!i^{*}\!j^{*}},
  \frac{\delta_2}{9(\sqrt{r}+1)},\delta_1\big\}$.
  Pick any $X\in\mathbb{B}(\overline{X},\delta)\cap\mathbb{R}_{+}^{n\times r}$.
  Define a matrix $\widetilde{X}\!\in\mathbb{R}^{n\times r}$ with
  \[
   \widetilde{X}_{ij}:=\left\{\begin{array}{cl}
     X_{ij}& {\rm if}\ j=l_i,\\
     0 &{\rm if}\ j\ne l_i
   \end{array}\right.\ {\rm for}\ (i,j)\in[n]\times[r].
  \]
  Let $D={\rm Diag}(\|\widetilde{X}_{\cdot1}\|,\|\widetilde{X}_{\cdot2}\|,
  \ldots,\|\widetilde{X}_{\cdot r}\|)$. It is easy to check that 
  $\widetilde{X}D^{-1}\in\mathcal{S}_{+}^{n,r}$.
  Choose a matrix $Z\in\mathbb{R}^{n\times (n-r)}$ such that 
  $[\widetilde{X}\!D^{-1}\ Z]$ is an $n\times n$ orthogonal matrix. 
  Then, it holds that
  \begin{align}\label{temp-ineq30}
   {\rm dist}\big(X,\mathcal{S}_{+}^{n,r}\big)
   \le\|X\!-\!\widetilde{X}D^{-1}\|_F
    =\big\|[\widetilde{X}D^{-1}\ \ Z]\tp (X\!-\!\widetilde{X}D^{-1})\big\|_F
   =\big\|[\widetilde{X}D^{-1}\ \ Z]\tp X\!-\!I_{n\times r}\big\|_F.
  \end{align}
  To invoke \eqref{Zr-dist} with $X'=[\widetilde{X}D^{-1}\ Z]\tp X$, 
  we proceed the arguments by the following two steps.
  \begin{enumerate}
  \renewcommand{\labelenumi}{(\theenumi)}
  \item To prove that $\big\|[\widetilde{X}D^{-1}\ Z]\tp X\!-\!I_{n\times r}\big\|_F
        \le c\| D^{-1}\widetilde{X}\tp X-I_{r}\|_F$. Write $H:=X-\overline{X}$. We have that 
       \begin{align}\label{temp-ineq311}
        \|Z\tp X\|_1
         &=\sum_{i=1}^{n-r}\sum_{j=1}^r\Big|\sum_{k=1}^nZ_{ki}(X_{kj}-\widetilde{X}_{kj})\Big|
          =\sum_{i=1}^{n-r}\sum_{j=1}^r
           \Big|\sum_{k=1}^n(1-\mathbbm{1}_{\{j=l_k\}})Z_{ki}X_{kj}\Big|\nonumber\\
          &=\sum_{i=1}^{n-r}\sum_{j=1}^r\Big|\sum_{k=1}^n(1-\mathbbm{1}_{\{j=l_k\}})Z_{ki}(\overline{X}_{kj}+H_{kj})\Big|\nonumber\\
          &\le\sum_{k=1}^n\sum_{i=1}^{n-r}\sum_{j=1,j\ne l_k}^{r}|Z_{ki}H_{kj}|
           \le (n-r)\sum_{k=1}^n\sum_{j=1,j\ne l_k}^{r}|H_{kj}|,
        \end{align}
        where the first equality is using $Z\tp \widetilde{X}=0$, the first inequality 
        is due to $\overline{X}_{kj}=0$ for $k\in[n]$ and $[r]\ni j\ne l_k$,
        and the last one is from $|Z_{ki}|\le 1$ for each $k\in[n]$ and 
        $i\in[n-r]$. On the other hand, from $\widetilde{X}\in\mathbb{R}_{+}^{n\times r}$ 
        and $X\in\mathbb{R}_{+}^{n\times r}$ it follows that
        \begin{align*}
        &\|D^{-1}\widetilde{X}\tp X\!-I_r\|_1
         \ge\frac{1}{\|D\|}\sum_{i=1}^{r}\sum_{j=1,j\ne i}^{r}|(\widetilde{X}\tp X)_{ij}|
         =\frac{1}{\|D\|}\sum_{i=1}^{r}\sum_{j=1,j\ne i}^{r}\big(\widetilde{X}\tp X\big)_{ij}\\
        &\ge\frac{1}{\|D\|}\sum_{i=1}^{r}\sum_{j=1,j\ne i}^{r}\sum_{k=1}^n\widetilde{X}_{ki}X_{kj}
          =\frac{1}{\|D\|}\sum_{i=1}^{r}\sum_{j=1,j\ne i}^{r}\Big[\sum_{k=1}^n\mathbbm{1}_{\{i=l_k\}}(\overline{X}_{ki}\! 
           +H_{ki})(\overline{X}_{kj}\!+H_{kj})\Big]\\
         &=\frac{1}{\|D\|}\sum_{i=1}^{r}\sum_{j=1,j\ne i}^{r}\sum_{k=1}^n\mathbbm{1}_{\{i=l_k\}}
          (\overline{X}_{ki}H_{kj}+H_{ki}H_{kj})\\
         &=\frac{1}{\|D\|}\sum_{k=1}^n\sum_{j=1,j\ne l_k}^{r}(\overline{X}_{kl_k}H_{kj}+H_{kl_k}H_{kj})
         \ge\frac{1}{\|D\|}\sum_{k=1}^n\sum_{j=1,j\ne l_k}^{r}
         \Big[\frac{1}{2}\overline{X}_{i^{*}j^{*}}|H_{kj}|\Big],
       \end{align*}
       where the second equality is because each row of $\overline{X}$
       has only one nonzero entry, and the last inequality is using $H_{kj}\ge 0$
       when $j\ne l_k$, $|H_{kl_k}|\le\delta\le\frac{1}{2}\overline{X}_{\!i^{*}\!j^{*}}$
       and $\overline{X}_{kl_k}\ge\overline{X}_{\!i^{*}\!j^{*}}$ for $k\in[n]$.
       Note that $D_{ii}\le\|X_{\cdot i}\|\le 1+\delta$ for each $i\in[r]$.
       Combining the last inequality with \eqref{temp-ineq311} yields 
       \[
        \frac{\overline{X}_{\!i^{*}\!j^{*}}}{n-r}\|Z\tp X\|_1
        \le 2\|D\|\|D^{-1}\widetilde{X}\tp X\!-\!I_r\|_1
        \le 3\|D^{-1}\widetilde{X}\tp X\!-\!I_r\|_1.
       \]
      Because $\|[\widetilde{X}D^{-1}\ Z]\tp X\!-\!I_{n\times r}\|_F
      \le \|D^{-1}\widetilde{X}\tp X-I_r\|_F+\|Z\tp X\|_F$, from the last inequality and the definition of the constant $c$, we obtain the stated inequality.
      
 \item To argue that $[\widetilde{X}D^{-1}\ Z]\tp X\in\mathbb{B}(I_{n\times r},\delta_2)$. 
       By the definition of $D$, for each $i\in[r]$, it holds that 
       \[
         D_{ii}=\|\widetilde{X}_{\cdot i}\|\ge\|X_{\cdot i}\|-\|\widetilde{X}_{\cdot i}-X_{\cdot i}\|
         \ge\|X_{\cdot i}\|-\delta\ge 1-2\delta,       
       \] 
       which along with $D_{ii}\le 1+\delta$ for each $i\in[r]$ implies that 
       $-\frac{2\delta}{1-2\delta}\le 1-\frac{1}{D_{ii}}\le\frac{\delta}{1+\delta}$. 
       Then, we have
       \begin{align*}
       \|X-XD^{-1}\|_F\le\|X\|_F\|I-D^{-1}\|\le\frac{\delta}{1+\delta}\big(\|X-\overline{X}\|_F+\|\overline{X}\|_F\big)\le\frac{\delta}{1+\delta}(\sqrt{r}+\delta),\\
       \|(X-\widetilde{X})D^{-1}\|_F\le\|D^{-1}\|\|X-\widetilde{X}\|_F
       \le\delta\|D^{-1}\|\le\frac{\delta}{1-2\delta}.\qquad\qquad
       \end{align*}
       Because $\|[\widetilde{X}D^{-1}\ Z]\tp X-I_{n\times r}\|_F
        =\|X-[\widetilde{X}D^{-1}\ Z]I_{n\times r}\|_F =\|X-\widetilde{X}D^{-1}\|_F$ 
        and $\|X-\widetilde{X}D^{-1}\|_F\le\|X-XD^{-1}\|_F+\|(X-\widetilde{X})D^{-1}\|_F$,
       from the last two inequalities, $\frac{\delta}{1+\delta}\le\frac{2\delta}{1-2\delta}$, 
       $0<\delta\le\frac{\delta_2}{9(\sqrt{r}+1)}$  and $\delta<1/3$, 
       it immediately follows that $\|[\widetilde{X}D^{-1}\ Z]\tp X-I_{n\times r}\|_F\le\delta_2$. 
  \end{enumerate}
  Now by invoking inequality \eqref{Zr-dist} with $X'=[\widetilde{X}D^{-1}\ Z]\tp X$, 
  it is not difficult to obtain that
  \begin{align*}
  2.1\sqrt{r}\|\sigma(X)-e\|=2.1\sqrt{r}\|\sigma([\widetilde{X}D^{-1}\ \ Z]\tp X)-e\|
  \ge\|D^{-1}\widetilde{X}\tp X-I_r\|_F
  \ge\frac{1}{c}\|[\widetilde{X}D^{-1}\ Z]\tp X-I_{n\times r}\|_F,
  \end{align*}
  where the last inequality is due to the result in Step 1. 
  Along with \eqref{temp-ineq30}, we get the desired result.
 \end{proof}

 Now we are ready to derive a local Lipschitzian error bound 
 for the constraint system $X\in\mathcal{S}_{+}^{n,r}$. 
 \begin{proposition}\label{prop-Stnr}
  Fix any $\overline{X}\in\mathcal{S}_{+}^{n,r}$. If $\overline{X}$ has 
  no zero rows when $n>r>1$, then there exists $\delta>0$ such that 
  for all $X\!\in\mathbb{B}(\overline{X},\delta)$, 
  ${\rm dist}\big(X,\mathcal{S}_{+}^{n,r}\big)
   \le(\kappa+1)[{\rm dist}(X,\mathbb{R}_{+}^{n\times r})
                 +{\rm dist}(X,{\rm St}(n,r))]$.
 \end{proposition}
 \begin{proof}
  By Lemma \ref{lemma2-NonStnr}, there exists $\delta'>0$ such that 
  ${\rm dist}(Z,\mathcal{S}_{+}^{n,r})\le \kappa\|\sigma(Z)-e\|$ 
  for all $Z\in\mathbb{B}(\overline{X},\delta')\cap\mathbb{R}_{+}^{n\times r}$.
  Pick any $X\!\in\mathbb{B}(\overline{X},{\delta'}/{2})$. 
  Let $X_{+}\!:=\!{\rm Proj}_{\mathbb{R}_{+}^{n\times r}}(X)$. 
  Clearly, $\|X_{+}-\overline{X}\|_F\le 2\|X-\overline{X}\|_F\le\delta'$. 
  Then, ${\rm dist}(X_{+},\mathcal{S}_{+}^{n,r})\le \kappa\|\sigma(X_{+})-e\|$. 
  Let $X_{+}$ have the SVD as $X_{+}=U{\rm Diag}(\sigma(X_{+}))V\tp $,
  where $U$ and $V$ are respectively an $n\times n$ and $r\times r$ orthogonal matrix. 
  We have
  ${\rm dist}(X_{+},{\rm St}(n,r))=\!\|X_{+}-UV\tp \|_F=\|\sigma(X_{+})-e\|$.
  From the above, we conclude that ${\rm dist}(X_{+},\mathcal{S}_{+}^{n,r})
  \le \kappa{\rm dist}(X_{+},{\rm St}(n,r))$. Consequently, 
  \begin{align*}
   {\rm dist}\big(X,\mathcal{S}_{+}^{n,r}\big)
   &\le\|X-X_{+}\|_F+{\rm dist}\big(X_{+},\mathcal{S}_{+}^{n,r}\big)
   \le\|X-X_{+}\|_F+\kappa{\rm dist}(X_{+},{\rm St}(n,r))\\
   &\le\|X-X_{+}\|_F+\kappa\big[\|X-X_{+}\|_F+{\rm dist}(X,{\rm St}(n,r))\big]\\
   &\le(\kappa+1)\big[{\rm dist}(X,\mathbb{R}_{+}^{n\times r})
       +{\rm dist}(X,{\rm St}(n,r))\big].
  \end{align*}
  This, by the arbitrariness of $X$ in $\mathbb{B}(\overline{X},\delta'/2)$,
  shows that the desired inequality holds.
  \end{proof}
 \begin{remark}
  When $n>r>1$, if $\overline{X}\in\mathcal{S}_{+}^{n,r}$ has a zero row, 
  the conclusion of Proposition \ref{prop-Stnr} does not necessarily hold.
  To see this, consider
  $\overline{X}=\left(\begin{matrix}
              1 & 0 & 0 \\
              0 & 1 & 0\\
            \end{matrix}\right)\tp\in\mathcal{S}_{+}^{3,2}$ and a sequence $\{X^k\}_{k\in\mathbb{N}}$ converging to $\overline{X}$ with 
 $X^k=\left(\begin{matrix}
              1 & 0 & {1}/{k} \\
              0 & 1 & {1}/{k} \\
            \end{matrix}\right)\tp $
 for each $k\in\mathbb{N}$. Because each row of any $X\in\mathcal{S}_{+}^{3,2}$ has at most one nonzero entry, we have ${\rm dist}\big(X^k,\mathcal{S}_{+}^{3,2}\big)\ge{1}/{k}$. However,   
 from the nonnegativity of $X^k$, it follows that 
 ${\rm dist}(X^k,\mathbb{R}_{+}^{n\times r})+{\rm dist}(X^k,{\rm St}(n,r))
 ={\rm dist}(X^k,{\rm St}(n,r))=\|\sigma(X^k)-e\|=\sqrt{1+\frac{2}{k^2}}-1=O(\frac{1}{k^2})$.
 \end{remark}
 
 Due to the compactness of ${\rm St}(n,r)$, when $n=r$ or $n>r=1$, 
 the following global error bound holds. 
 \begin{corollary}\label{gerror}
  When $n=r$ or $n>r=1$, there exists a constant $\kappa'>0$ such that 
  for any $Z\in{\rm St}(n,r)$, ${\rm dist}(Z,\mathcal{S}_{+}^{n,r})
  \le\kappa'{\rm dist}(Z,\mathbb{R}_{+}^{n,r})$.
 \end{corollary} 
 \begin{proof}
  Fix any $Z\in{\rm St}(n,r)$. From Proposition \ref{prop-Stnr}, 
  for each $X\in\mathcal{S}_{+}^{n,r}$, there exists $\delta_{X}$ such that 
  for all $Y\in\mathbb{B}(X,\delta_{X})$, 
  ${\rm dist}(Y,\mathcal{S}_{+}^{n,r})\le\kappa[{\rm dist}(Y,\mathbb{R}_{+}^{n\times r})
  +{\rm dist}(Y,{\rm St}(n,r))]$. Because $\bigcup_{X\in\mathcal{S}_{+}^{n,r}}\mathbb{B}^{\circ}(X,\delta_{X})$ is 
  an open covering of the compact set $\mathcal{S}_{+}^{n,r}$, according to
  Heine-Borel covering theorem, there exist $X^1,X^2,\ldots,X^p\in\mathcal{S}_{+}^{n,r}$ 
  such that $\mathcal{S}_{+}^{n,r}\subset\bigcup_{i=1}^p\mathbb{B}^{\circ}(X^i,\delta_{X^i}):=D$.
  When $Z\in D$, it holds that
  ${\rm dist}(Z,\mathcal{S}_{+}^{n,r})\le\kappa{\rm dist}(Z,\mathbb{R}_{+}^{n\times r})$, 
  so we only need to consider the case that $Z\notin D$. 
  Let $\overline{D}={\rm cl}[{\rm St}(n,r)\backslash D]$. 
  There must exist a constant $\widetilde{\kappa}>0$ 
  such that $\min_{Y\in\overline{D}}{\rm dist}(Y,\mathbb{R}_{+}^{n\times r})
  \ge\widetilde{\kappa}$. If not, there exists a sequence 
  $\{Y^k\}_{k\in\mathbb{N}}\subset\overline{D}$ such that 
   ${\rm dist}(Y^k,\mathbb{R}_{+}^{n\times r})\le 1/k$, 
   which by the compactness of the set $\overline{D}$ and the continuity of 
   the distance function means that there is a cluster point, 
   say $\overline{Y}\in\overline{D}$, of $\{Y^k\}_{k\in\mathbb{N}}$ 
   such that $\overline{Y}\in\mathbb{R}_{+}^{n\times r}$. 
   Then $\overline{Y}\in\mathcal{S}_{+}^{n,r}\subset D$, 
   a contradiction to the fact that $\overline{Y}\in\overline{D}$.
   In addition, because both $\overline{D}$ and $\mathcal{S}_{+}^{n,r}$ 
   are compact, there exists a constant $c_0>0$ such that for all 
   $Y\in\overline{D}$, ${\rm dist}(Y,\mathcal{S}_{+}^{n,r})\le c_0$.
   Together with $\min_{Y\in\overline{D}}{\rm dist}(Y,\mathbb{R}_{+}^{n\times r})
  \ge\widetilde{\kappa}$ and $Z\in\overline{D}$, it holds that 
  ${\rm dist}(Z,\mathcal{S}_{+}^{n,r})\le (c_0/\widetilde{\kappa})
  {\rm dist}(Z,\mathbb{R}_{+}^{n\times r})$. Thus,
  ${\rm dist}(Z,\mathcal{S}_{+}^{n,r})\le\kappa'{\rm dist}(Z,\mathbb{R}_{+}^{n\times r})$ 
  with $\kappa'=\max\{c_0/\widetilde{\kappa},\kappa\}$. 
  The desired result then follows.
 \end{proof}
  
  Now we are in a position to establish the global exact penalty result 
  for problem \eqref{Nsmooth-pen}.
 \begin{theorem}\label{theorem1-epenalty}
  If every global minimizer of problem \eqref{orth-prob} has no zero rows when $n>r>1$, then 
  \begin{enumerate}[1.]
  \item  when $n=r$ or $n>r=1$, for all $X\in{\rm St}(n,r)$, 
         $f(X)-f^*+\kappa'L_{\!f}\vartheta(X)\ge 0$;
         
  \item  when $n>r>1$, for every global minimizer $X^*$ of \eqref{orth-prob},
         there exists $\delta>0$ such that for all $\epsilon\ge0$
         and all $X\in\mathbb{B}(X^*,\delta)\cap\mathcal{F}_{\epsilon}$ 
         with $\mathcal{F}_{\epsilon}:=\big\{X\in{\rm St}(n,r)\,|\, 
         \vartheta(X)=\epsilon\big\}$,
         $f(X)-f^*+\kappa L_{\!f}\,\vartheta(X)\ge 0$;
  \end{enumerate}
  where $f^*$ is the optimal value of problem \eqref{orth-prob}, and consequently there exists a threshold $\overline{\rho}>0$ such that problem \eqref{Nsmooth-pen} associated to each $\rho\ge\overline{\rho}$
  has the same global optimal solution set as \eqref{orth-prob} does.
 \end{theorem}
 \begin{proof}
  We first consider the case $n>r>1$. By Proposition \ref{prop-Stnr}, 
  there exists $0<\delta_0<1$ such that
  \begin{equation}\label{temp-ineq31}
   {\rm dist}\big(Z,\mathcal{S}_{+}^{n,r}\big)
   \le\kappa\big[{\rm dist}(Z,\mathbb{R}_{+}^{n\times r})+{\rm dist}(Z,{\rm St}(n,r))\big]
   \ \ {\rm for\ all}\ Z\in\mathbb{B}(X^*,\delta_0).
  \end{equation}
  Fix any $\epsilon\ge 0$ and any $X\in\mathbb{B}(X^*,\delta_0)\cap\mathcal{F}_{\epsilon}$. 
  Then, there exists $\overline{X}\!\in\mathcal{S}_{+}^{n,r}$ such that 
  \[
    \|X\!-\!\overline{X}\|_F
    ={\rm dist}(X,\mathcal{S}_{+}^{n,r})
    \le \kappa\big[{\rm dist}(X,\mathbb{R}_{+}^{n\times r})+{\rm dist}(X,{\rm St}(n,r))\big]
    \le\kappa\vartheta(X),
  \]
  where the last inequality is due to $X\in{\rm St}(n,r)$ implied by 
  $X\in\mathcal{F}_{\epsilon}$. In addition, from $\overline{X}\in\mathcal{S}_{+}^{n,r}$, 
  \[
    f(X)-f^*\ge f(X)-f(\overline{X})\ge-L_{\!f}\|X-\overline{X}\|_F. 
  \]
  The last two inequalities imply that $f(X)-f^*+\kappa L_{\!f}\vartheta(X)\ge 0$. 
  By the arbitrariness of $\epsilon\ge 0$, the desired result holds.
  For the case $n>r=1$ or $n=r$, by combining the last inequality and Corollary \ref{gerror}, 
  it is easy to obtain the conclusion. Thus, the first part of the conclusions follows.
  The second part holds by the first part, \citep[Proposition 2.1(b)]{LiuBiPan18}, 
  and the compactness of $\mathcal{S}_{+}^{n,r}$. 
 \end{proof}
 
 By following the proof of Theorem \ref{theorem1-epenalty}, 
 it is not difficult to achieve the following conclusion.
 \begin{corollary}\label{local-Epenalty}
  Suppose that each local minimizer of problem \eqref{orth-prob} has no zero rows when $n>r>1$.  
  Then, for every local minimizer $X^*$ of \eqref{orth-prob},
  \begin{enumerate}[1.]
  \item  when $n=r$ or $n>r=1$, there exists $\delta>0$ such that for all 
         $X\in\mathbb{B}(X^*,\delta)\cap{\rm St}(n,r)$,
         \[
            f(X)-f(X^*)+\kappa'L_{\!f}\vartheta(X)\ge 0;
         \]
         
  \item  when $n>r>1$, there exists $\delta>0$ such that for all $\epsilon\ge0$
         and all $X\in\mathbb{B}(X^*,\delta)\cap\mathcal{F}_{\epsilon}$, 
         where $\mathcal{F}_{\epsilon}$ is defined as above,
         $f(X)-f(X^*)+\kappa L_{\!f}\,\vartheta(X)\ge 0$.
  \end{enumerate}
  Hence, each local minimizer of \eqref{orth-prob} is locally optimal 
  to problem \eqref{Nsmooth-pen} with $\rho\ge\max(\kappa',\kappa)L_{\!f}$, 
  and conversely, each nonnegative local minimizer of \eqref{Nsmooth-pen}
  with $\rho>0$ is also locally optimal to \eqref{orth-prob}. 
 \end{corollary} 
 \begin{remark}\label{remark-Epenalty}
 {\bf(a)} When the distance in ${\rm dist}(X,\mathbb{R}_{+}^{n\times r})$ is defined by other norms of $\mathbb{R}^{n\times r}$, the conclusion of Proposition \ref{prop-Stnr} still holds (but with a different $\kappa$), and from the proof of Theorem \ref{theorem1-epenalty}, problem \eqref{Nsmooth-pen} is still a global exact penalty of \eqref{orth-prob} if the penalty term is replaced by a distance on other norms to $\mathbb{R}_{+}^{n\times r}$.
 
 \noindent
 {\bf(b)} In \cite[Theorem 4.1]{Han79}, Han and Mangasarian showed that if $\overline{X}$
  is globally optimal to all penalty problems associated to $\rho>\overline{\rho}$,
  then it is globally optimal to \eqref{orth-prob}. Such a condition 
  is very strong and unverifiable. Our global exact penalty result requires 
  that each global optimal solution has no zero rows when $n>r>1$, 
  and as will be shown by Example \ref{example1} below, 
  this requirement cannot be removed.     
 \end{remark}
 \begin{example}\label{example1}
  Consider $\min_{X\in\mathcal{S}_{+}^{3,2}}f(X)$ with $f(X)=-2X_{11}-2X_{22}-\min\{X_{31},X_{32}\}$ for $X\in\mathbb{R}^{3\times 2}$. Note that $-4$ is a lower bound of $f$ on the nonnegative orthogonal set $\mathcal{S}_{+}^{3,2}$. Hence, 
    $X^*=\left(\begin{matrix}
        1 & 0 & 0\\ 0 & 1 & 0
    \end{matrix}\right)^{\top}$ and $f^*=-4$. In fact, using GloptiPoly 3 (see \cite{Henrion09}) can vadilate that $X^*$ is indeed a global optimal solution of the problem. 
    For each $k$, let  
    $X^k=\left(\begin{matrix}
        \sqrt{1-\frac{1}{k^2}} & 0 & 1/k\\ 
        -\frac{1}{k^2} & \sqrt{1-\frac{1}{k^2}} & \frac{1}{k}\sqrt{1-\frac{1}{k^2}}
    \end{matrix}\right)^{\top}$. 
    It is easy to check that 
    $X^k\in\big\{X\in{\rm St}(3,2)\ |\ \vartheta(X)=\frac{1}{k^2}\big\}$. 
    Note that $f(X^k)-f^*=4-4\sqrt{1-\frac{1}{k^2}}-\frac{1}{k}\sqrt{1-\frac{1}{k^2}}=-\frac{1}{k}+O(\frac{1}{k^2})$ but $\vartheta(X^k)=\frac{1}{k^2}$. Hence, 
    $\min_{X\in{\rm St}(3,2)}\big\{f(X)+\rho\vartheta(X)\big\}$ 
    is not a global exact penalty for $\min_{X\in\mathcal{S}_{+}^{3,2}}f(X)$. 
 \end{example}

  Next we establish the global exact penalty result for problem \eqref{smooth-pen} 
  under a mild condition for $f$.
 \begin{theorem}\label{theorem2-epenalty}
  Fix any $\gamma>0$. Suppose that every global (or local) minimizer of \eqref{orth-prob} has no zero rows when $n>r>1$, and that for every global (or local) minimizer $X^*$ of \eqref{orth-prob}, there exist $\delta'>0$ and $L'>0$ such that 
  for all $X\in\mathbb{B}(X^*,\delta')\cap{\rm St}(n,r)$ and all
  $\overline{X}\in{\rm Proj}_{\mathcal{S}_{+}^{n,r}}(X)$,
  \begin{equation}\label{growth}
    f(X)-f(\overline{X})\ge -L'\|X-\overline{X}\|_F^2.
  \end{equation}
  Then, for every global (or local) minimizer $X^*$ of \eqref{orth-prob}, 
  \begin{enumerate}[1.]
  \item when $n=r$ or $n>r=1$, there exists $\delta>0$ such that for all 
        $X\in\mathbb{B}(X^*,\delta)\cap{\rm St}(n,r)$, 
        \[
         f(X)-f(X^*)+2\gamma(\kappa')^2 L'e_{\gamma}\vartheta(X)\ge 0;
       \] 
        
  \item when $n>r>1$, there exists $\delta>0$ such that for all $\epsilon\ge0$ 
       and all $X\!\in\mathbb{B}(X^*,\delta)\cap\mathcal{G}_{\epsilon}$ 
       where $\mathcal{G}_{\epsilon}:=\{X\in{\rm St}(n,r)\,|\, 
        e_{\gamma}\vartheta(X)=\epsilon\}$,
       $f(X)-f(X^*)+2\gamma\kappa^2 L'e_{\gamma}\vartheta(X)\ge 0$,
  \end{enumerate}
  and consequently, there exists a threshold $\widehat{\rho}>0$ such that
  the penalty problem \eqref{smooth-pen} associated to every $\rho\ge\widehat{\rho}$
  has the same global optimal solution set as problem \eqref{orth-prob} does.
 \end{theorem}
 \begin{proof}
  It suffices to consider the local minimizer case, and for the global minimizer case, 
  the proof is similar. Fix any local minimizer $X^*$ of problem \eqref{orth-prob}. 
  There exists $0<\epsilon'\le\gamma$ such that
  \begin{equation}\label{local-optimal}
    f(X)\ge f(X^*)\quad{\rm for\ all}\ X\in\mathbb{B}(X^*,\epsilon')\cap\mathcal{S}_{+}^{n,r}.
  \end{equation}
  By continuity, there exists $0<\delta_1\le\epsilon'$ such that 
  for all $X\in\mathbb{B}(X^*,\delta_1)$ and all $(i,j)\in{\rm supp}(X^*)$, $X_{ij}>0$. 
  We first consider Item 2 of the conclusions. Now inequality \eqref{temp-ineq31} 
  still holds. Set $\delta:=\min(\epsilon',\delta',\delta_0)/2$ where $\delta_0$ is 
  same as the one in \eqref{temp-ineq31}. Pick any $\epsilon\ge 0$
  and fix any $X\!\in\mathbb{B}(X^*,\delta)\cap\mathcal{G}_{\epsilon}$.
  From \eqref{temp-ineq31} and $X\in{\rm St}(n,r)$, there exists 
  $\overline{X}\!\in{\rm Proj}_{\mathcal{S}_{+}^{n,r}}(X)$ such that 
  \begin{equation}\label{Xbar-ineq}
   \|\overline{X}-X^*\|_F\le\epsilon'\ \ {\rm and}\ \ 
   \|X-\overline{X}\|_F^2\le\kappa^2{\rm dist}^2(X,\mathbb{R}_{+}^{n\times r}).
  \end{equation}
  Let $\overline{J}_1=\big\{(i,j)\in\overline{{\rm supp}}(X^*)\,|\,X_{ij}\in[-\gamma,0]\big\}$
  and $\overline{J}_2=\big\{(i,j)\in\overline{{\rm supp}}(X^*)\,|\,X_{ij}>0\big\}$.
  Because $\|X-X^*\|_{F}\le\delta\le\epsilon'\le\gamma$, $X_{ij}\ge-\gamma$ 
  for all $(i,j)\in\overline{{\rm supp}}(X^*)$, which means that $\overline{{\rm supp}}(X^*)
  =\overline{J}_1\cup\overline{J}_2$. By Lemma \ref{prox-vartheta}, 
  \begin{equation}\label{Pgamma-equa}
    [\mathcal{P}_{\!\gamma}\vartheta(X)]_{ij}
    =\left\{\begin{array}{cl}
      X_{ij}&{\rm if}\ (i,j)\in{\rm supp}(X^*)\cup\overline{J}_2,\\
      0 &{\rm if}\ (i,j)\in\overline{J}_1.
     \end{array}\right. 
  \end{equation}
  Together with the expression of $e_{\gamma}\vartheta(X)$ in Lemma \ref{prox-vartheta},  
  it is immediate to check that 
  \[
    e_{\gamma}\vartheta(X)=\frac{1}{2\gamma}\sum_{(i,j)\in\overline{J}_1}X_{ij}^2
    \ \ {\rm and}\ \
    {\rm dist}(X,\mathbb{R}_{+}^{n\times r})
    =\sqrt{\sum_{(i,j)\in\overline{J}_1}|X_{ij}|^2},
  \]
  which implies that $e_{\gamma}\vartheta(X)=\frac{1}{2\gamma}
  {\rm dist}^2(X,\mathbb{R}_{+}^{n\times r})$.
  Together with \eqref{growth}-\eqref{Xbar-ineq}, it then follows that
  \[
    f(X)\!-\!f(X^*)\ge f(X)\!-\!f(\overline{X})\ge -L'\|X\!-\!\overline{X}\|_F^2
    \ge -2\gamma\kappa^2L'e_{\gamma}\vartheta(X).
  \]
  This shows that Item 2 of the conclusions holds.
  By using Corollary \ref{gerror} and the same arguments as above, 
  we can prove that Item 1 of the conclusions holds with 
  $\delta =\min(\epsilon,\delta')/2$. 
 \end{proof}
 \begin{remark}\label{remark-Epenalty1}
  Note that for all $X\in\mathbb{B}(X^*,\delta')\cap{\rm St}(n,r)$ 
  and all $\overline{X}\in{\rm Proj}_{\mathcal{S}_{+}^{n,r}}(X)$, 
  $\|\overline{X}-X^*\|_F\le 2\delta'$. Hence, the assumption in \eqref{growth}
  is actually the lower second-order calmness of $f$ relative to a neighborhood of the global (or local) minimizer set of \eqref{orth-prob}, which is not necessarily stronger than its locally Lipschitz 
 continuity relative to this neighborhood. 
 \end{remark}
 
 Because ${\rm dist}^2(X,\mathbb{R}_{+}^{n\times r})=\|\max(0,-X)\|_F^2$ 
 for any $X\in\mathbb{R}^{n\times r}$, the proof of Theorem \ref{theorem2-epenalty} 
 implies that under the lower second-order calmness assumption of $f$ in \eqref{growth},  problem \eqref{Quad-Pen} is also a global exact penalty for \eqref{orth-prob}, i.e., the following conclusion holds.
\begin{corollary}\label{corollary2-epenalty}
  Under the assumptions of Theorem \ref{theorem2-epenalty}, 
  for every global (or local) minimizer $X^*$ of \eqref{orth-prob}, 
  \begin{enumerate}[1.]
  \item when $n=r$ or $n>r=1$, there exists $\delta>0$ such that for all 
        $X\in\mathbb{B}(X^*,\delta)\cap{\rm St}(n,r)$, 
        \[
         f(X)-f(X^*)+(\kappa')^2L'\|\max(0,-X)\|_F^2\ge 0;
        \] 
        
  \item when $n>r>1$, there exists $\delta>0$ such that for all $\epsilon\ge0$ 
       and all $X\!\in\mathbb{B}(X^*,\delta)\cap\mathcal{G}_{\epsilon}$ 
       where $\mathcal{G}_{\epsilon}:=\{X\in{\rm St}(n,r)\,|\, 
        \|\max(0,-X)\|_F^2=\epsilon\}$,
       $f(X)-f(X^*)+\kappa^2 L'\|\max(0,-X)\|_F^2\ge 0$,
  \end{enumerate}
  and consequently there exists a threshold $\widehat{\rho}>0$ such that
  the penalty problem \eqref{Quad-Pen} associated to every $\rho\ge\widehat{\rho}$
  has the same global optimal solution set as problem \eqref{orth-prob} does.
 \end{corollary} 

 To close this section, we provide another application of Proposition \ref{prop-Stnr} in verifying the Clarke regularity of the set $\mathcal{S}_{+}^{n,r}$ and achieving its normal cone under a mild condition. It is worth pointing out that the regular normal cone to $\mathcal{S}_{+}^{n,r}$ was characterized in \citep[Section 2.1]{JiangM22}. 
 \begin{corollary}\label{regularity}
  Consider any $X\in\mathcal{S}_{+}^{n,r}$. If the matrix $X$ has no zero rows when $n>r>1$, 
  then      
  \begin{align*}
	\widehat{\mathcal{N}}_{\mathcal{S}_{+}^{n,r}}(X)
	&=\mathcal{N}_{\mathcal{S}_{+}^{n,r}}(X)
	=\mathcal{N}_{\mathbb{R}_{+}^{n\times r}}(X)
	+\mathcal{N}_{{\rm St}(n,r)}(X)\\
	&=\big\{V\in\mathbb{R}^{n\times r}\ |\ V_{ij}\in\lambda_j X_{ij}\ \ {\rm with}\ \ \lambda_j\in\mathbb{R},\ \ \forall(i,j)\in{\rm supp}(X)\big\}.
 \end{align*}
 \end{corollary}
 \begin{proof}
  By combining Proposition \ref{prop-Stnr} and \cite[Section 3.1]{Ioffe08}, 
  we have the following inclusion:
  \[
	\mathcal{N}_{\mathcal{S}_{+}^{n,r}}(X)\subset
	\mathcal{N}_{\mathbb{R}_{+}^{n\times r}}(X)
	+\mathcal{N}_{{\rm St}(n,r)}(X).
  \]
  On the other hand, from \citep[Theorem 6.42]{RW98} and the regularity of 
	${\rm St}(n,r)$, it follows that 
	\[
	\widehat{\mathcal{N}}_{\mathcal{S}_{+}^{n,r}}(X)\supset
	\mathcal{N}_{\mathbb{R}_{+}^{n\times r}}(X)
	+\mathcal{N}_{{\rm St}(n,r)}(X).
	\]
	Combining the last two inclusions with $\widehat{\mathcal{N}}_{\mathcal{S}_{+}^{n,r}}(X)
	\subset\mathcal{N}_{\mathcal{S}_{+}^{n,r}}(X)$ yields the first two equalities. 
	The last equality follows by the first equality and \citep[Eq. (2.1)\&(2.2)]{JiangM22}.
 \end{proof} 
 \section{Exact penalty method for problem \eqref{orth-prob}}\label{sec4}
 
 From Section \ref{sec3}, problem \eqref{Nsmooth-pen} associated to each $\rho\ge\overline{\rho}$ is equivalent to problem \eqref{orth-prob} in a global sense if every global minimizer of \eqref{orth-prob} has no rows when $n>r>1$; and so are problems \eqref{smooth-pen} and \eqref{Quad-Pen} associated to each $\rho\ge\widehat{\rho}$ provided that $f$ satisfies the assumptions of Theorem \ref{theorem2-epenalty}. 
 When solving a single penalty problem with a certain $\rho\ge\overline{\rho}$ 
 or $\widehat{\rho}$, one cannot expect a solution of high quality due to the nonconvexity 
 of \eqref{orth-prob} and its penalty problems. Hence, it is necessary to solve 
 a series of penalty problems with increasing $\rho$. Solving 
 the nonsmooth penalty problem \eqref{Nsmooth-pen} is much more difficult than 
 solving the smooth penalty problems \eqref{smooth-pen} and \eqref{Quad-Pen}, though some active 
 exploration were given to the solution methods of the former (see \citep{Chen20,Huang21}). In particular, the asymptotic convergence of the penalty method 
 based on \eqref{Nsmooth-pen} requires that every penalty problem is solved to 
 an exact stationary point, which is impractical in computation. Inspired by 
 the feasible method proposed in \citep{Wen13} for optimization with orthogonality constraints, 
 we develop the following penalty method by seeking the approximate stationary points 
 of the penalty problem \eqref{smooth-pen} or \eqref{Quad-Pen}. Although the global exact penalty of \eqref{smooth-pen} and \eqref{Quad-Pen} requires that every global minimizer of \eqref{orth-prob} has no rows when $n>r>1$, as will be shown in Theorem \ref{convergence} below, the convergence results of this penalty method do not need this requirement.
 \renewcommand{\thealgorithm}{\arabic{algorithm}}
 \begin{algorithm}[h]
  \setcounter{algorithm}{0}
 \caption{(Penalty methods based on $\Theta_{\rho,\gamma}$)}\label{PenAlg}
 \begin{algorithmic}[1]
 \State Initialization:
  Input $\gamma\ge 0,l_{\rm max}\in\mathbb{N},\epsilon>0$. Choose $\rho_{\rm max}>0,
  \tau_{\rm min}\in(0,1),\sigma_{\!\rho}>1,\sigma_{\tau}\in(0,1),
  \rho_0>0$ and $\tau_0>0$.
  Choose $X^0\in{\rm St}(n,r)$. Set $X^{0,0}=X^0$ and $\upsilon^0=\Theta_{\rho_0,\gamma}(X^0)$.

 \For{$l=0,1,2,\ldots,l_{\rm max}$}
 \State\label{step3} Starting from $X^{l,0}$, to seek a point $X^{l+1}\in{\rm St}(n,r)$ 
      such that $\|{\rm Proj}_{{\rm T}_{\!X^{l\!+\!1}}M}(\nabla\Theta_{\rho_l,\gamma}(X^{l\!+\!1}))\|_F
      \!\le\!\tau_l$ \hspace*{0.2in} and $\Theta_{\!\rho_l,\gamma}(X^{l\!+\!1})\!\le\!\upsilon^l$ by solving the following penalized problem
        \begin{equation}\label{EP-subprob}
          \min_{X\in{\rm St}(n,r)}\Theta_{\!\rho_l,\gamma}(X).
         \end{equation}
        
  \State  If $\vartheta(X^{l+1})\le\epsilon$, then stop; otherwise 
          let $\rho_{l+1}\!=\min\{\sigma_{\!\rho}\rho_l,\rho_{\rm max}\}$ 
          and $\tau_{l+1}\!=\max\{\sigma_{\!\tau}\tau_l,\tau_{\rm min}\}$.
          
  \State\label{step5} Construct $\widetilde{X}^{l+1}\in\mathcal{S}_{+}^{n,r}$ with $X^{l+1}$. If $\Theta_{\!\rho_{l+1},\gamma}(\widetilde{X}^{l+1})
          <\Theta_{\!\rho_{l+1},\gamma}(X^{l+1})$, then set $X^{l+1,0}=\widetilde{X}^{l+1}$ and  \hspace*{0.2in}$
          \upsilon^{l+1}=$ $\Theta_{\!\rho_{l+1},\gamma}(\widetilde{X}^{l+1})$; else set 
          $X^{l+1,0}=X^{l+1}$ and $\upsilon^{l+1}=\Theta_{\!\rho_{l+1},\gamma}(X^{l+1})$.
 \EndFor 
 \end{algorithmic}
 \end{algorithm}
 \begin{remark}\label{remark-Alg1}
  {\bf(a)}\ From Remark \ref{Remark-subdiff} (a) and Corollary \ref{subdiff-EPfun}, 
  $X^*\in{\rm St}(n,r)$ is a stationary point of problem \eqref{EP-subprob} if and only if ${\rm Proj}_{{\rm T}_{\!X^*}M}(\nabla\Theta_{\rho_{l},\gamma}(X^*))= 0$. This means that Step \ref{step3} of Algorithm \ref{PenAlg} is seeking an approximate stationary point $X^{l+1}$ of \eqref{EP-subprob} with $\Theta_{\!\rho_l,\gamma}(X^{l+1})\le\upsilon^l$. As will be shown by Item 2 of Theorem \ref{theorem1-ManPG1} and Remark \ref{remark-algA} (a), when Algorithm \ref{ManPG1} is applied to solving \eqref{EP-subprob}, it can return such $X^{l+1}$.

  \noindent
  {\bf(b)} The condition $\Theta_{\!\rho_l,\gamma}(X^{l+1})\le\upsilon^l$ requires that the objective value of $X^{l+1}$ is not more than $\Theta_{\rho_{l},\gamma}(\widetilde{X}^{l})$, where $\widetilde{X}^{l}$ is a feasible point 
  constructed from $X^{l}$. As will be shown in Section \ref{sec4.1}, the introduction of $\widetilde{X}^{l}$ is crucial to guarantee that each cluster point of the sequence $\{X^l\}_{l\in\mathbb{N}}$ is feasible to \eqref{orth-prob}. In practical computation, one may adopt \cite[Procedure 1]{JiangM22} to yield such a feasible point $\widetilde{X}^{l}$. 
 \end{remark} 
 \subsection{Convergence analysis of Algorithm \ref{PenAlg}}\label{sec4.1}
  
  Before analyzing the convergence of Algorithm \ref{PenAlg}, we introduce the notion of stationary point for \eqref{orth-prob}. 
  \begin{definition}\label{spoint}
  	We call $X\in\mathcal{S}_{+}^{n,r}$ a stationary point of problem \eqref{orth-prob}
  	if $0\in\nabla\!f(X)+ {\rm N}_{X}M+\mathcal{N}_{\mathbb{R}_{+}^{n\times r}}(X)$
  	or equivalently 
  	$0\in{\rm Proj}_{{\rm T}_{\!X}M}\big(\nabla\!f(X)
  	+\mathcal{N}_{\mathbb{R}_{+}^{n\times r}}(X)\big)$ holds.
  \end{definition}
  
  The following lemma states that the stationary point in Definition \ref{spoint} 
  is the same as the one in \citep{JiangM22}, and Proposition \ref{SSOC} in Appendix \ref{appendix_A} provides a second-order sufficient condition for a stationary point of problem \eqref{orth-prob} to be a strong local optimal solution by leveraging the second subderivative.
  \begin{lemma}\label{spoint-equiv}
  	A matrix $X\in \mathcal{S}_{+}^{n,r}$ is a stationary point of 
  	problem \eqref{orth-prob} iff
  	it is a stationary point of \eqref{orth-prob} in terms of 
  	$-\nabla\!f(X)\in[\mathcal{T}_{\mathcal{S}_{+}^{n,r}}(X)]^{\circ}$, 
  	where $K^{\circ}$ denotes the negative polar of a cone $K$.
  \end{lemma}
  \begin{proof}
  	When $n>r=1$ or the matrix $X$ has no zero rows when $n>r>1$, from Corollary \ref{regularity} 
  	we have
  	\[
  	{\rm N}_{X}M+\mathcal{N}_{\mathbb{R}_{+}^{n\times r}}(X)
  	=\widehat{\mathcal{N}}_{\mathcal{S}_{+}^{n,r}}(X)
  	=[\mathcal{T}_{\mathcal{S}_{+}^{n,r}}(X)]^{\circ},
  	\] 
  	which together with \cite[Theorem 2.1]{JiangM22} implies the desired conclusion. 
  	Thus, it suffices to consider the case when $X$ has zero rows when $n>r>1$. 
  	Note that 
  	${\rm N}_{X}M+\mathcal{N}_{\mathbb{R}_{+}^{n\times r}}(X)
  	\subset\widehat{\mathcal{N}}_{\mathcal{S}_{+}^{n,r}}(X)
  	=[\mathcal{T}_{\mathcal{S}_{+}^{n,r}}(X)]^{\circ}$.
  	We only need to prove that the converse inclusion holds. 
  	Pick any $H\in[\mathcal{T}_{\mathcal{S}_{+}^{n,r}}(X)]^{\circ}$.   
  	Define the index sets
  	$\overline{J}_1(X):=\{(i,j)\in\overline{{\rm supp}}(X)\,|\,\|X_{i\cdot}\|>0\}$ 
  	and $\overline{J}_2(X)=\overline{{\rm supp}}(X)\backslash\overline{J}_1(X)$.  
  	By the expression of $[\mathcal{T}_{\mathcal{S}_{+}^{n,r}}(X)]^{\circ}$
  	in \citep[Lemma 2.2]{JiangM22}, there exists $\lambda_j$ such that $H_{ij}=\lambda_j X_{ij}$ for each $(i,j)\in{\rm supp}(X)$, and 
  	for each $(i,j)\in\overline{J}_2(X)$, $H_{ij}\le 0$. 
  	For each $(i,j)\in[n]\times[r]$, let $S_{ij}=\lambda_j$ if $i=j$, 
  	$S_{ij}=H_{sj}/{X_{sl_s}}$ if $(s,j)\in\overline{J}_1(X)$ and $i=l_s$, 
  	otherwise $S_{ij}=0$, where $l_s:=\mathop{\arg\max}_{j\in[r]}X_{sj}$. Then, after an elementary calculation, it holds that
  	\[
  	(XS)_{ij}=\left\{\begin{array}{cl}
  		X_{ij}S_{jj} &{\rm if}\ (i,j)\in{\rm supp}(X), \\
  		X_{il_i}S_{l_ij} &{\rm if}\  (i,j)\in\overline{J}_1(X), \\
  		0 &{\rm if}\  (i,j)\in\overline{J}_2(X)
  	\end{array}\right.\ {\rm and}\ \  
  	(H-XS)_{ij}=\left\{\begin{array}{cl}
  		0 &{\rm if}\ (i,j)\in{\rm supp}(X), \\
  		0 &{\rm if}\ (i,j)\in\overline{J}_1(X), \\
  		H_{ij} & {\rm if}\  (i,j)\in\overline{J}_2(X).
  	\end{array}\right.
  	\]
  	Clearly, $H-XS\in\mathcal{N}_{\mathbb{R}_+^{n\times r}}(X)$. 
  	Then $H\in{\rm N}_{X}M+\mathcal{N}_{\mathbb{R}_+^{n\times r}}(X)$.
  	The converse inclusion follows.
  \end{proof}
      
  Now we focus on the asymptotic convergence of Algorithm \ref{PenAlg}. 
  The main result is stated as follows. 
  \begin{theorem}\label{convergence}
   Let $\{X^l\}_{l\in\mathbb{N}}$ be the sequence yielded by Algorithm \ref{PenAlg} 
   with $l_{\rm max}=+\infty,\,\epsilon=0,\,\tau_{\rm min}=0$ and $\rho_{\rm max}=+\infty$.
   Then, $\{X^l\}_{l\in\mathbb{N}}$ is bounded and its every 
   cluster point is a stationary point of \eqref{orth-prob}.  
  \end{theorem} 
  \begin{proof}
   The boundedness of $\{X^l\}_{l\in\mathbb{N}}$ is immediate because 
   $\{X^l\}_{l\in\mathbb{N}}\subset{\rm St}(n,r)$. Let $X^{\infty}$ be an arbitrary 
   cluster point of $\{X^l\}_{l\in\mathbb{N}}$. Then there exists an index set 
   $\mathcal{L}\subset\mathbb{N}$ such that $\lim_{\mathcal{L}\ni l\to\infty}X^{l}=X^{\infty}$. 
   
   \noindent
   {\bf Case 1: $\gamma>0$.} In this case, based on Step \ref{step5} of Algorithm \ref{PenAlg}, it follows that for all $l\in\mathbb{N}$, 
   \[
     f(X^{l})+\rho_{l-1}e_{\gamma}\vartheta(X^{l})
     =\Theta_{\!\rho_{l-1},\gamma}(X^{l})\le\Theta_{\!\rho_{l-1},\gamma}(\widetilde{X}^{l-1})
     =f(\widetilde{X}^{l-1}).
  \]
  Because the sequence $\{f(\widetilde{X}^{l-1})\}_{l\in\mathbb{N}}$ is bounded 
  and $\liminf_{l\to\infty}f(X^{l})>-\infty$, the last inequality and the nonnegativity of $e_{\gamma}\theta$ implies that 
  $e_{\gamma}\vartheta(X^{\infty})=0$. From $e_{\gamma}\vartheta(X^{\infty})=0$ 
  and the expression of $e_{\gamma}\vartheta$, we have  
  \[
    0=e_{\gamma}\vartheta(X^{\infty})
    =\frac{1}{2\gamma}\|\mathcal{P}_{\gamma}\vartheta(X^{\infty})-X^{\infty}\|_F^2
     +\vartheta(\mathcal{P}_{\gamma}\vartheta(X^{\infty}))
     \ge\frac{1}{2\gamma}\|\mathcal{P}_{\gamma}\vartheta(X^{\infty})-X^{\infty}\|_F^2,
  \]
  which implies that $\mathcal{P}_{\gamma}\vartheta(X^{\infty})=X^{\infty}$ 
  and $\vartheta(X^{\infty})=0$. Hence, $X^{\infty}$ is a feasible point of problem \eqref{orth-prob}. From Step \ref{step3} of Algorithm \ref{PenAlg}, for each $l\in\mathbb{N}$, there exists $H^{l}\in\mathbb{R}^{n\times r}$ with $\|H^{l}\|_F\le\tau_{l-1}$ such that 
  \begin{align*}
    H^{l}={\rm Proj}_{{\rm T}_{\!X^{l}}M}(\nabla\Theta_{\rho_{l-1},\gamma}(X^{l}))
    &={\rm Proj}_{{\rm T}_{\!X^{l}}M}\big(\nabla\!f(X^{l})+\rho_{l-1}\gamma^{-1}(X^{l}-
    \mathcal{P}_{\gamma}\vartheta(X^{l}))\big)\\
    &={\rm grad}f(X^{l})+{\rm Proj}_{{\rm T}_{\!X^{l}}M}\big[\rho_{l-1}\gamma^{-1}(X^{l}-
    \mathcal{P}_{\gamma}\vartheta(X^{l}))\big].
  \end{align*}
  Since $\lim_{\mathcal{L}\ni l\to\infty}X^{l}=X^{\infty}$, 
  by the expression of $\mathcal{P}_{\gamma}\vartheta$ 
  in Lemma \ref{prox-vartheta} (see also \eqref{Pgamma-equa}), there exists $\widehat{l}\in\mathbb{N}$ 
  such that for all $\mathcal{L}\ni l\ge\widehat{l}$, 
  $X^{l}-\mathcal{P}_{\gamma}\vartheta(X^{l})\in\mathcal{N}_{\mathbb{R}_{+}^{n\times r}}
  (X^{\infty})$. Along with the last inequality, for each $\mathcal{L}\ni l\ge\widehat{l}$,
  \[
    H^{l}-{\rm grad}f(X^{l})\in{\rm Proj}_{{\rm T}_{\!X^{l}}M}
    \big[\mathcal{N}_{\mathbb{R}_{+}^{n\times r}}(X^{\infty})\big].
  \]
  Passing the limit $\mathcal{L}\ni l\to\infty$ to this inclusion 
  and using the continuity of ${\rm Proj}_{{\rm T}_{\!X}M}$ on the variable $X$ and \cite[Proposition 5.52]{RW98}, 
  we obtain $-{\rm grad}f(X^{\infty})\in{\rm Proj}_{{\rm T}_{\!X^{\infty}}M}
  \big[\mathcal{N}_{\mathbb{R}_{+}^{n\times r}}(X^{\infty})\big]$. 
  This, together with Definition \ref{spoint} and the feasibility of $X^{\infty}$, 
  shows that $X^{\infty}$ is a stationary point of \eqref{orth-prob}.
  
  \noindent
  {\bf Case 2: $\gamma=0$.} In this case, based on Step \ref{step5} of Algorithm \ref{PenAlg}, it follows that for all $l\in\mathbb{N}$, 
  \[
    f(X^{l})+\rho_{l-1}\|\max(0,-X^{l})\|_F^2
    =\Theta_{\!\rho_{l-1},\gamma}(X^{l})\le\Theta_{\!\rho_{l-1},\gamma}(\widetilde{X}^{l-1})
    =f(\widetilde{X}^{l-1}).
  \]
  Then, following the same arguments as those of Case 1 yileds the desired conclusion. 
  \end{proof}
  
  From Lemma \ref{spoint-equiv}, our stationary point is same as the one defined 
  in \citep[Eq.\,(2.7a)-(2.7b)]{JiangM22}. Then, by comparing Theorem \ref{convergence} 
  with the convergence results in \cite[Section 4.2]{JiangM22}, we see that any cluster 
  point yielded by the penalty method based on \eqref{smooth-pen} or \eqref{Quad-Pen} is a stationary point defined in \citep[Eq.\,(2.7a)-(2.7b)]{JiangM22} without requiring the cluster point to have no zero rows when $n>r>1$. This partly confirms the superiority of the penalty problems \eqref{smooth-pen} and \eqref{Quad-Pen}. In addition, the following corollary shows that Algorithm \ref{PenAlg} returns a local minimizer of \eqref{orth-prob} within a finite number of steps, provided that a local minimizer $X^{l+1}$ of some penalty problem \eqref{EP-subprob} associated to $\rho_{l}>\widehat{\rho}$ lies in a neighborhood of a point $X^*\in\mathbb{R}_{+}^{n\times r}$ which is locally optimal to \eqref{smooth-pen} or \eqref{Quad-Pen} associated to $\widehat{\rho}$. A similar result was achieved in \citep[Corollary 4.2]{JiangM22} by requiring that every penalty problem is solved to a local optimal solution. 
 \begin{corollary}\label{finite-stop} 
  Suppose that for some $\rho_{l}>\widehat{\rho}$ and $X^*\in\mathbb{R}_{+}^{n\times r}$, there exist $\varepsilon>0$ and a local minimizer $X^{l+1}$ of \eqref{EP-subprob} with  $X^{l+1}\in\mathbb{B}(X^*,\varepsilon)$ such that  $\Theta_{\rho_{l},\gamma}(X^{l+1})\le\Theta_{\rho_{l},\gamma}(X^*)$ and $  \Theta_{\widehat{\rho},\gamma}(X^*)\le\Theta_{\widehat{\rho},\gamma}(Z)$
  for all $Z\in\mathbb{B}(X^*,\varepsilon)\cap{\rm St}(n,r)$. Then, $X^{l+1}$ is locally optimal to problem \eqref{orth-prob}.
  \end{corollary} 
 \begin{proof}
 From the given assumption, $X^{l+1}\in{\rm St}(n,r)$ and the expression of $\Theta_{\rho,\gamma}$, it follows that   
 \[
  \Theta_{\rho_{l},\gamma}(X^{l+1})\le\Theta_{\rho_{l},\gamma}(X^*)
  =f(X^*)=\Theta_{\widehat{\rho},\gamma}(X^*)
  \le\Theta_{\widehat{\rho},\gamma}(X^{l+1}), 
 \]
 where the equalities are due to $X^*\in\mathbb{R}_{+}^{n\times r}$.
  When $\gamma=0$, from the last inequality, it follows that
  \[
    f(X^{l+1})+\rho_{l}\big\|\max(-X^{l+1},0)\big\|_F^2\le f(X^{l+1})+\widehat{\rho}\big\|\max(-X^{l+1},0)\big\|_F^2,
  \]
  which along with $\rho_{l}>\widehat{\rho}$ implies that $X^{l+1}\in\mathbb{R}_{+}^{n\times r}$, and $X^{l+1}$ is a feasible point of \eqref{orth-prob}. When $\gamma>0$, 
  \[
  f(X^{l+1})+\rho_{l}e_{\gamma}\vartheta(X^{l+1})\le f(X^{l+1})+\widehat{\rho}e_{\gamma}\vartheta(X^{l+1}),
  \]
  which along with $\rho_{l}>\widehat{\rho}$ implies that $e_{\gamma}\vartheta(X^{l+1})=0$, and then $X^{l+1}$ is a feasible point of \eqref{orth-prob}. Together with the local optimality of $X^{l+1}$ to \eqref{EP-subprob}, we conclude that $X^{l+1}$ is a local minimizer of  \eqref{orth-prob}. 
  \end{proof} 
  
 \subsection{PGM on manifold for solving subproblem \eqref{EP-subprob}}\label{sec4.2}
 
  Fix $\rho=\rho_l$. Inspired by the work \citep{Chen20}, we develop a PGM
  on manifold to find an approximate stationary point of \eqref{EP-subprob}. 
  This PGM first computes a tangent direction $V^k$ to $M$ given by 
  \begin{equation}\label{subprob}
   V^k:=\mathop{\arg\min}_{V\in{\rm T}_{\!X^k}M}
   \Big\{\langle\nabla\Theta_{\rho,\gamma}(X^k),V\rangle+\frac{1}{2t_k}\|V\|_F^2\Big\}
   =-t_k\,{\rm grad}\Theta_{\rho,\gamma}(X^k),
  \end{equation}
  and then pullbacks it to the manifold $M$ by leveraging the retraction mapping $R_{X^k}$. 
  It is well known that the efficiency of the PGM depends much on the choice of 
  the step-size $t_k$, which is determined by the Lipschitz constant $L_{\rho,\gamma}$
  of $\nabla\Theta_{\rho,\gamma}$. Taking into account that $L_{\rho,\gamma}$ is unknown 
  if $L_{\nabla\!f}$ is unachievable, we develop a nonmonotone line-search PGM 
  on manifold with $V^k$ in \eqref{subprob}.
 \begin{algorithm}[H]
 \renewcommand{\thealgorithm}{A}
 \caption{(Nonmonotone line-search PGM on manifold for \eqref{EP-subprob})}\label{ManPG1}
 \begin{algorithmic}[1] 
 \State Initialization: Fix $\rho=\rho_l$. Choose $\eta\in(0,1),\alpha>0, 
         m\in\mathbb{N},0<t_{\rm min}\le t_{\rm max}$ and $X^0\in{\rm St}(n,r)$. 
        
 \For{$k=0,1,2,\ldots$}
 
 \State\label{Astep3} Choose $t_k\in[t_{\rm min},t_{\rm max}]$ by the BB rule;
 
 \State Let $V^k=-t_k\,{\rm grad}\Theta_{\!\rho,\gamma}(X^k)$ 
        and set $X^{k+1}=R_{X^k}(V^k)$;
 	
	    \While{$\Theta_{\rho,\gamma}(X^{k+1})>\max_{j=\max(0,k-m),\ldots,k}
	    \Theta_{\rho,\gamma}(X^j)-\frac{\alpha}{2t_k}\|V^k\|_F^2$}\label{Astep5}
	    
	    \State Let $t_k\leftarrow\eta t_k$ and  
	           $V^k\leftarrow -t_k\,{\rm grad}\Theta_{\!\rho,\gamma}(X^k)$;
 	    
	    \State\label{Astep7} Set $X^{k+1}=R_{X^k}(V^k)$;
	    \EndWhile \label{Astep8}

 \EndFor  
 \end{algorithmic}
 \end{algorithm}
  A good step-size initialization at each outer iteration can greatly reduce 
  the line-search cost. Inspired by the work \cite{Wright09}, we initialize 
  the step-size $t_k$ in Step \ref{Astep3} by the Barzilai-Borwein (BB) rule \cite{Barzilai88}: 
  \begin{equation}\label{BB-stepsize}
    t_k=\max\Big\{\min\Big\{\frac{\|\Delta X^k\|_F^2}{|\langle\Delta X^k,\Delta Y^k\rangle|},
    \frac{|\langle\Delta X^k,\Delta Y^k\rangle|}{\|\Delta Y^k\|_F^2},t_{\rm max}\Big\},
    t_{\rm min}\Big\},
  \end{equation}
  where $\Delta X^k\!:=X^k\!-\!X^{k-1}$ and $\Delta Y^k\!:={\rm grad}\Theta_{\rho,\gamma}(X^k)
  \!-\!{\rm grad}\Theta_{\rho,\gamma}(X^{k-1})$. Such a step-size is accepted when 
  the nonmonotone line-search criterion (see \cite{Grippo86,Grippo02}) in Step \ref{Astep5} is satisfied. The following lemma shows that the nonmonotone line-search criterion 
  in Step \ref{Astep5} is well defined.
 \begin{lemma}\label{lemma1-ManPG1}
  For each $k\in\mathbb{N}$, the nonmonotone line-search criterion 
  in Step \ref{Astep5} is satisfied whenever $t_k\le\frac{1-\alpha}{2c_2c_{\rho,\gamma}
  +c_1^2L_{\rho,\gamma}}$, where $c_{\rho,\gamma}\!:={\displaystyle\max_{Z\in{\rm St}(n,r)}}
  \|\nabla\Theta_{\rho,\gamma}(Z)\|_F$ and $c_1,c_2$ are the constants from Lemma \ref{property1-manifold}.
 \end{lemma}
 \begin{proof}
  Recall that $\nabla\!f$ is Lipschitz continuous on ${\rm St}(n,r)$ 
  with constant $L_{\nabla\!f}$. Hence, $\nabla\Theta_{\rho,\gamma}$ is Lipschitz 
  continuous on ${\rm St}(n,r)$ with constant $L_{\rho,\gamma}$. From the descent lemma, 
  \begin{equation}\label{Thetarho-ineq1}
    \Theta_{\rho,\gamma}(X^{k+1})-\Theta_{\rho,\gamma}(X^k)
    \le\langle\nabla\Theta_{\rho,\gamma}(X^k),X^{k+1}-X^k\rangle
    +(L_{\rho,\gamma}/2)\|X^{k+1}-X^k\|_F^2.
  \end{equation}
  Notice that $X^{k+1}=R_{X^k}(V^k)$. By Lemma \ref{property1-manifold}, 
  $\|X^{k+1}-X^k-V^k\|_F\le c_2\|V^k\|_F^2$, and moreover,
  \begin{align}\label{Thetarho-ineq2}
    \langle\nabla\Theta_{\rho,\gamma}(X^k),X^{k+1}-X^k\rangle
    &=\langle\nabla\Theta_{\rho,\gamma}(X^k),X^{k+1}-X^k-V^k\rangle
     +\langle\nabla\Theta_{\rho,\gamma}(X^k),V^k\rangle \nonumber\\
    &\le c_2\|\nabla\Theta_{\rho,\gamma}(X^k)\|_F\|V^k\|_F^2 
     +\langle{\rm grad}\Theta_{\rho,\gamma}(X^k),V^k\rangle \nonumber\\
    &\le c_2\|\nabla\Theta_{\rho,\gamma}(X^k)\|_F\|V^k\|_F^2-\frac{1}{2t_k}\|V^k\|_F^2
    \nonumber\\
    &\le c_2c_{\rho,\gamma}\|V^k\|_F^2-\frac{1}{2t_k}\|V^k\|_F^2
  \end{align}
  where the second inequality is also using \eqref{subprob}.
  From \eqref{Thetarho-ineq1}-\eqref{Thetarho-ineq2} 
  and $t_k\le\frac{1-\alpha}{2c_2c_{\rho,\gamma}
  	+c_1^2L_{\rho,\gamma}}$, we have
  \begin{equation*}
   \Theta_{\rho,\gamma}(X^{k+1})
   \le\Theta_{\rho,\gamma}(X^k)-\Big(\frac{1}{2t_k}-c_2c_{\rho,\gamma}
    -\frac{1}{2}c_1^2L_{\rho,\gamma}\Big)\|V^k\|_F^2
    \le-\frac{\alpha}{2t_k}\|V^k\|_F^2.
  \end{equation*}
  Since $\Theta_{\rho,\gamma}(X^k)
  \le{\displaystyle\max_{j=\max(0,k-m),\ldots,k}}\Theta_{\rho,\gamma}(X^j)$, 
  the desired result follows.
 \end{proof}
 
 The following theorem states the convergence results of Algorithm \ref{ManPG1}, including 
 the convergence of $\{\Theta_{\rho,\gamma}(X^k)\}_{k\in\mathbb{N}}, 
 \{V^k\}_{k\in\mathbb{N}}$ and the whole iterate sequence, whose proof is put in Appendix \ref{appendix_B}.
 \begin{theorem}\label{theorem1-ManPG1}
  Let $\{(X^k,V^k)\}_{k=0}^{\infty}$ be the sequence generated by Algorithm \ref{ManPG1}.
  For each $k$, write $\ell(k):=\mathop{\arg\max}_{j=\max(0,k-m),\ldots,k}
  \Theta_{\rho,\gamma}(X^j)$. Then, the following assertions hold.
  \begin{enumerate}[1.]
   \item The sequences $\{\Theta_{\rho,\gamma}(X^k)\}_{k\in\mathbb{N}}$ 
          and $\{V^k\}_{k\in\mathbb{N}}$ are convergent, and $\lim_{k\to\infty}V^k=0$.

   \item The sequence $\{X^k\}_{k\in\mathbb{N}}$ is bounded and every cluster point
          is a stationary point of \eqref{EP-subprob}.
               
   \item Suppose that $f$ is definable in an o-minimal structure $\mathscr{O}$ 
         over $\mathbb{R}$, and that 
         \begin{equation}\label{mild-cond}
          \sum_{\mathcal{K}\ni k=0}^{\infty}
         \sqrt{\Theta_{\rho,\gamma}(X^{\ell(k+1)})-\Theta_{\rho,\gamma}(X^{k+1})}<\infty
         \ \ {\rm when}\ 
         \liminf_{\mathcal{K}\ni k\to\infty}
         \frac{\Theta_{\rho,\gamma}(X^{\ell(k)})-\Theta_{\rho,\gamma}(X^{\ell(k+1)})}
         {\|X^{k+1}-X^k\|_F^2}=0
         \end{equation}         
        with $\mathcal{K}\!:=\!\big\{k\in\mathbb{N}\ |\ 
          \Theta_{\rho,\gamma}(X^{\ell(k+1)})-\Theta_{\rho,\gamma}(X^{k+1})
          \ge\frac{\alpha\|X^{k+1}-X^k\|_F^2}{6t_{\max}}\big\}$. 
        Then $\sum_{k=0}^{\infty}\|X^{k+1}-X^k\|_F<\infty$. 
  \end{enumerate}
  \end{theorem}
  \begin{remark}\label{remark-algA}
  {\bf(a)}	From the proof of Item 1 of Theorem \ref{theorem1-ManPG1}, the sequence $\{\Theta_{\rho,\gamma}(X^{\ell(k)})\}_{k\in\mathbb{N}}$ 
  	is monotonically decreasing. Together with the definition of $l(k)$, 
  	it follows that for each $k\in\mathbb{N}$, 
  	\[
  	\Theta_{\rho,\gamma}(X^{k})\le\Theta_{\rho,\gamma}(X^{\ell(k)})\le\Theta_{\rho,\gamma}(X^{\ell(0)})=\Theta_{\rho,\gamma}(X^{0}).
  	\]
  
  \noindent
  {\bf(b)} Item 3 of Theorem \ref{theorem1-ManPG1} establishes the convergence of the whole 
  iterate sequence generated by a nonmonotone descent method under the KL assumption 
  on $f$ and condition \eqref{mild-cond} on the nonmonotonicity of 
  the objective values. When $m=0$, condition \eqref{mild-cond} is superfluous. 
 \end{remark} 
\section{Numerical experiments}\label{sec5}

We employ Algorithm \ref{PenAlg} armed with Algorithm \ref{ManPG1} to solve QAP problems, graph matching problems, projection problems onto $\mathcal{S}_{+}^{n,k}$, and orthogonal nonnegative 
matrix factorization. Among others, Algorithm \ref{PenAlg} based on $\Theta_{\rho,\gamma}$ and $\Theta_{\rho,0}$ are respectively called SEPPG+ and SEPPG0. 
We compare their performance with that of the ALM \cite{Wen13} and the exact penalty method EP4Orth+ \citep{JiangM22}. All experiments are performed in MATLAB on a workstation running 
on 64-bit Windows System with an Intel(R) Xeon(R) W-2245 CPU 3.90GHz and 128 GB RAM.
\subsection{Implementation details of SEPPG+ and SEPPG0}\label{sec5.1}

First of all, we take a closer look at the choice of the parameters in Algorithm \ref{PenAlg}. Our preliminary numerical tests indicated that when $\rho_l$ is small, if $X^{l+1,0}$ is feasible, Algorithm \ref{PenAlg} has poor numerical performance. In view of this, in the subsequent tests, we set $X^{l+1,0}=X^{l+1}$ and $\upsilon^{l+1}=\Theta_{\!\rho_{l+1},\gamma}(X^{l+1})$ when $\rho_{l}<10^3$. The parameter $\rho_0$ and the initial $X^0$ are specified in the experiments. Recall that $\vartheta(X)=0$ if and only if $e_{\gamma}\vartheta(X)=0$. The parameter $\gamma>0$ is unnecessary too small. Consider that a smaller $\gamma$ leads to a worse smoothness of $\Theta_{\rho_{l},\gamma}$ and then requires more computation cost, so we choose $\gamma=0.05$ for SEPPG+. In addition, the increasing ratio $\sigma_{\rho}$ of the penalty parameter $\rho$ has a major influence on the quality of solutions. A smaller $\sigma_{\rho}$ usually leads to better solutions, but more penalty problems need to be solved. We take a trade off by choosing $\sigma_{\rho}=1.05$ if $\rho_l\le 1$, otherwise $\sigma_{\rho}=1.1$. 
We choose $\tau_0=1.0$ for SEPPG+ and $\tau_0=0.005$ for SEPPG0. Other parameters in Algorithm \ref{PenAlg} are chosen as follows:
\[
l_{\rm max}=2000,\,\epsilon=10^{-6},\,\rho_{\rm max}=10^{10},\,
\tau_{\rm min}=10^{-5},\,\sigma_{\tau}=0.95.
\]
We terminate Algorithm \ref{PenAlg} at $X^{l}$ whenever $\vartheta(X^{l})\le\epsilon$, or $\vartheta(X^{l})\le 5\epsilon$ and 
$\frac{|f(X^{l})-f(X^{l-9})|}{1+|f(X^{l})|}\le 10^{-8}$.

For the implementation of Algorithm \ref{ManPG1}, we choose 
$\eta=0.1,\alpha=10^{-4},m=5,t_{\rm min}=10^{-12}$ and $t_{\rm max}=10^{12}$, 
and adopt the QR decomposition for the retraction in Step \ref{Astep7}, 
which is computed with \textrm{myQR} from the toolbox \textbf{OptM}; 
see \url{https://github.com/optsuite/OptM}. 
\subsection{ALM and its implementation}\label{sec6.2}

For a given $\mu>0$, let $L_{\mu}(X,\Lambda)$ be the augmented Lagrangian function of \eqref{orth-prob} defined as
\[
L_{\mu}(X,\Lambda):=f(X)+\frac{\mu}{2}\|\min(0,X-\mu^{-1}\Lambda)\|_F^2-\frac{1}{2\mu}\|\Lambda\|_F^2.
\]
The ALM proposed in \cite{Wen13} for solving problem \eqref{orth-prob} consists of the following iterates:
\begin{subnumcases}{}\label{XY1-subprob}
	X^{k+1}\in\mathop{\arg\min}_{X\in{\rm St}(n,r)}L_{\mu_k}(X,\Lambda^k);\\
	\label{Lambda1-subprob}
	\Lambda^{k+1}=\max(\Lambda^k\!-\!\mu_kX^{k+1},0);\\
	\label{muk1-subprob}
	\mu_{k+1}=1.2\mu_k.
\end{subnumcases}
The iterates \eqref{XY1-subprob}-\eqref{muk1-subprob} are actually the ALM for the equivalent reformulation
of \eqref{orth-prob}:
\begin{equation}\label{orth-prob1}
	\min_{X,Y\in\mathbb{R}^{n\times r}}\Big\{f(X)+\delta_{\mathbb{R}_{+}^{n\times r}}(Y)\ \ {\rm s.t.}\ \
	X-Y=0,X\!\in{\rm St}(n,r)\Big\}.
\end{equation}
Indeed, for a given $\mu>0$, the augmented Lagrangian function for \eqref{orth-prob1} takes the form of
\[
L_{\mu}(X,Y;\Lambda)=f(X)+\delta_{{\rm St}(n,r)}(X)+\delta_{\mathbb{R}_{+}^{n\times r}}(Y)-\langle\Lambda,X-Y\rangle
+\frac{\mu}{2}\|X-Y\|_F^2,
\]
and the ALM for problem \eqref{orth-prob1} consists of the following iterates:
\begin{subnumcases}{}\label{XY-subprob}
	(X^{k+1},Y^{k+1})\in\mathop{\arg\min}_{X,Y\in\mathbb{R}^{n\times r}}L_{\mu_k}(X,Y;\Lambda^k);\\
	\label{Lambda-subprob}
	\Lambda^{k+1}=\Lambda^k-\mu_k(X^{k+1}-Y^{k+1});\\
	\label{muk-subprob}
	\mu_{k+1}=\varrho\mu_k\quad{\rm for\ some}\ \varrho>1.
\end{subnumcases}
Note that $\min_{X,Y\in\mathbb{R}^{n\times r}}L_{\mu_k}(X,Y;\Lambda^k)
=\min_{X\in\mathbb{R}^{n\times r}}[\min_{Y\in\mathbb{R}^{n\times r}}L_{\mu_k}(X,Y;\Lambda^k)]$.
An elementary calculation shows that
$\min_{Y\in\mathbb{R}^{n\times r}}L_{\mu_k}(X,Y;\Lambda^k)=L_{\mu_k}(X,\Lambda^k)+\delta_{{\rm St}(n,r)}(X)$ with the unique optimal solution $Y_{k}^*(X)=\max(0,X\!-\!\mu_k^{-1}\Lambda^k)$ and $Y^{k+1}=Y_{k}^*(X^{k+1})$. Thus, the iterates \eqref{XY1-subprob}-\eqref{muk1-subprob}
are precisely those \eqref{XY-subprob}-\eqref{muk-subprob} with $\varrho=1.2$.
For the implementation of the above ALM, the reader is referred to
the toolbox \textbf{OptM}; see \url{https://github.com/optsuite/OptM}. 

In the subsequent testing, we run OptM with $\textrm{opts.tol}=10^{-6},
\textrm{opts.omxitr}=10^3$, $\textrm{opts.tolsub}=10^{-6},\textrm{opts.mxitr}=100,
\textrm{opts.xtol}=10^{-5}, \textrm{opts.gtol}=10^{-5},
\textrm{opts.ftol}=10^{-8}$ and $\textrm{opts.tau}=10^{-3}$. 
The initial penalty parameter $\mu_0$ is specified in the experiment.  
\subsection{Quadratic assignment problems}\label{sec6.3}

Because all entries of every feasible point of a QAP are binary, 
problem \eqref{QAP} is equivalent to
\begin{equation}\label{QAP1}
	\min_{X\in\mathbb{R}^{n\times n}}
	\Big\{\langle A, (X\circ X)B(X\circ X)^{\top}\rangle\ \ {\rm s.t.}\ \
	X^{\top}X=I_n,\ X\ge 0\Big\},
\end{equation}
where ``$\circ $'' denotes the Hadamard product. When $A$ and $B$ are 
nonnegative matrices, say, those instances from QAPLIB, the lifted reformulation 
\eqref{QAP1} is superior to the original model \eqref{QAP} because now the objective 
function of \eqref{QAP1} is nonnegative at any $X\in\mathbb{R}^{n\times n}$ 
and will account for a larger proportion, compared with the objective function of  
\eqref{orth-prob}, in the penalty function $\Theta_{\rho,\gamma}$. 
The lifted reformulation \eqref{QAP1} was first adopted by Wen and Yin \cite{Wen13} 
when applying the ALM to the QAPs. In this section, we test SEPPG+ and SEPPG0 on the \textbf{134} instances from QAPLIB by solving \eqref{QAP1} except ``\textbf{esc16f}'' and 
``\textbf{tai10b}'' because the matrix $A$ for \textbf{esc16f} is zero and 
the best bound for \textbf{tai10b} is not provided. We also compare the performance of
SEPPG+ and SEPPG0 with that of ALM \cite{Wen13} and EP4Orth+ \citep{JiangM22} for these instances by solving \eqref{QAP1}. We measure the quality of a solution yielded by a solver 
in terms of the relative gap and the nonnegative infeasibility, which are respectively defined by 
\[
\textbf{relgap}:=\Big[\frac{f(X^{\!f})-\textrm{Best}}{\textrm{Best}}\times 100\Big]\%
\ \ {\rm and}\ \ 
\textbf{Ninf}:=\vartheta(X^{\!f}). 
\]
 Here, $X^{\!f}$ denotes the output of a solver and \textrm{Best} is the best upper bound 
of an instance. In the tests, we find that the orthogonal infeasibility yielded 
by four solvers are all less than $10^{-10}$ except that EP4Orth+ fails to ``tai256c'' 
(see Table \ref{table4}), so we do not report this index.

For SEPPG+ and SEPPG0, we choose the initial $\rho_0$ to be $c_0|f(X^0)|/\vartheta(X^0)$ for some $c_0>0$, where $X^0$ is the starting point specified later, and for ALM, we choose 
$\mu_0=10$, more stable than $\mu_0=1$ and $\mu_0=0.1$. 
For EP4Orth+, all parameters are default except $\textrm{sigma\_max}=10^{10}$ 
and $\textrm{maxiter}=2000$. For the four solvers, because the quality of the returned solutions depends on the starting point, we run each QAPLIB instance from \textbf{100} different starting points and take the average of the \textbf{100} results. Among others, SEPPG+, SEPPG0 and ALM are using the same starting point, generated randomly by the MATLAB command $\textrm{orth}(\widetilde{X}^0)$ with $\widetilde{X}^0=\textrm{randn}(n,r)$, while EP4Orth+ is using the starting point $|\widetilde{X}^0|{\rm diag}(\frac{1}{\|\widetilde{X}_{\cdot 1}^0\|},\ldots,\frac{1}{\|\widetilde{X}_{\cdot r}^0\|})$. 

Table \ref{table1} reports the average results of four solvers on the 21 instances 
with $n\ge 80$, including the average median gaps, the average CPU time in second, 
and the average nonnegative infeasibility for each instance. Consider that a rounding technique is imposed on the solutions returned by EP4Orth+ so that they always belong to 
$\mathcal{S}_{+}^{n,r}$. For SEPPG+, SEPPG0 and ALM, we also reported the average median gaps of 
the solutions imposed on the same rounding technique, denoted by ``rmed''. We see that 
the rounding technique has no effect on the quality of the solutions yielded by 
SEPPG+, SEPPG0 and ALM, and moreover, SEPPG+ and SEPPG0 yield the best average median gap 
respectively for \textbf{6} and \textbf{13} instances, while ALM yields the best average median gap for \textbf{2} instances. SEPPG+ and SEPPG0 take the comparable CPU time as EP4Orth+ does, but ALM requires the least CPU time.
\begin{table}[h]
	\renewcommand\arraystretch{1.3}
	\setlength{\abovecaptionskip}{2pt}
	\setlength{\belowcaptionskip}{0pt}
	\centering
	\captionsetup{font={small}}
	\caption{Numerical results of four solvers on 21 QAPLIB instances with $n\ge80$}\label{table1}
	\tiny
	\begin{tabular*}{\textwidth}{@{\extracolsep{\fill}}cllllllllllllllll@{\extracolsep{\fill}}}
		\hline
		Name& \multicolumn{4}{l}{SEPPG+}&\multicolumn{4}{l}{SEPPG0}&\multicolumn{3}{l}{EP4Orth+}&\multicolumn{4}{l}{ALM}\\
		\cmidrule(lr){2-5} \cmidrule(lr){6-9} \cmidrule(lr){10-12} \cmidrule(lr){13-16}
		&  med     &  rmed    &time  &Ninf  &  med     &  rmed    &time  &Ninf   &  med     &  rmed    &time   &  med     &  rmed    &time  &Ninf  \\
		&  gap(\%) &  gap(\%) &      &      & gap(\%)  & gap(\%)  &   &     &  gap(\%) & gap(\%)   &      &  gap(\%) &  gap(\%) &      &          \\
		\hline
		esc128      & {\bf 3.281}  & {\bf 3.281}  & 102.3   &  3.6e-7  & 4.156  & 4.156  & 73.9  & 3.4e-6  & 194.000  & 194.000   & 56.0  &25.938  &25.938 &0.7  &3.0e-5\\
		lipa80a     & 0.821  & 0.821  & 23.3   &  1.4e-7  & {\bf 0.778}  & {\bf 0.778}  & 17.7  & 1.6e-6   & 2.028  &  2.028 & 0.7   & 0.791  & 0.791  &1.3  &6.4e-6\\
		lipa80b     & {\bf 11.676}  & {\bf 11.676}  &3.0   & 3.1e-7  & 14.891  & 14.891  & 1.6  & 0   & 27.207  & 27.207    & 22.9  & 20.226 &20.226 &1.6  &2.7e-5\\
		lipa90a     & 0.745  & 0.745  & 31.3   &  1.7e-7  & {\bf 0.703}  & {\bf 0.703}  & 18.6  & 6.1e-7   & 1.820  & 1.820   &  1.0  & 0.718 & 0.718  &1.6  &8.0e-6\\
		lipa90b     & {\bf 11.429}  & {\bf 11.429}  & 4.2   & 4.6e-8  & 13.541  & 13.541  & 1.8  & 0   & 27.700  & 27.700    & 23.8  & 21.165 &21.165 &2.1  &2.6e-5\\
		sko81       & 1.623  & 1.623  & 5.1   &  3.2e-6  & {\bf 1.618}  & {\bf 1.618}  & 5.1  & 0   & 13.240  & 13.240   &  33.0  & 1.747  &1.747  &0.9  &1.8e-5\\
		sko90       & {\bf 1.608}  & {\bf 1.608}  & 6.3   &  2.9e-6  & { 1.637}  & { 1.637}  & 5.7  & 0   & 12.870  & 12.870   &  38.7  & 1.634 &1.634  &1.0 &1.5e-5\\
		sko100a     & { 1.480}  & { 1.480}  & 8.5   &  3.7e-6  & {\bf 1.426}  & {\bf 1.426}  & 7.2  &  0  & 12.773  & 12.773    & 51.7 & 1.620 &1.620  &1.4  &2.0e-5\\
		sko100b     & { 1.460}  & { 1.460} & 8.4   &  3.8e-6  & {\bf 1.451}  & {\bf 1.451}  & 7.1  &  0  & 12.290  & 12.290    & 51.0  & 1.547 &1.547  &1.4  &2.0e-5\\
		sko100c     & { 1.678}  & { 1.678}  & 8.8   &  3.5e-6 & {\bf 1.622}  & {\bf 1.622}  & 7.5  &  0  & 13.234 & 13.234    & 51.6   & 1.731&1.731  &1.4  &2.0e-5\\
		sko100d     & { 1.507}  & { 1.507}  & 8.5   &  3.8e-6 & {\bf 1.471}  & {\bf 1.471}  & 7.0  & 0   & 12.540  & 12.540    & 51.7  & 1.649&1.649  &1.4  &2.1e-5\\
		sko100e     & { 1.681}  & { 1.681}  & 8.5   &  3.9e-6 & {\bf 1.660}  & {\bf 1.660}  & 7.2  & 0   & 13.478  & 13.478     & 52.5  & 1.790&1.790  &1.4  &2.2e-5\\
		sko100f     & { 1.489}  & { 1.489}  & 8.8   &  3.7e-6 & {\bf 1.464}  & {\bf 1.464}  & 7.3  & 0   & 12.250  & 2.250   & 51.7   & 1.542&1.542  &1.4  &1.8e-5\\
		tai80a      & 3.009  & 3.009  & 3.1   &  2.3e-9 & 3.020  & 3.020  & 1.3  &  0  & 9.481  & 9.481   &  22.3  & {\bf 2.876} &{\bf 2.876}  &1.5  &2.1e-5\\
		tai80b      & { 4.792}  & { 4.792}  & 47.6  & 0  & {\bf 4.580}  & {\bf 4.580}  & 48.4  & 0   & 33.349  & 33.349   &  34.7   & 4.975&4.975  &4.0  &1.2e-5\\
		tai100a     & {\bf 2.734}  & {\bf 2.734}  & 5.2   &  3.3e-8 & 2.737  & 2.737  & 2.1  &  0  & 8.720  & 8.720    & 29.9   & 2.743&2.743  &2.5  &2.2e-5\\
		tai100b     & 4.271  & 4.271 & 84.1   &  2.1e-12 & {\bf 3.908}  & {\bf 3.908}  & 76.9  &  0  & 39.591  & 39.591   & 58.5   & { 4.208}&{ 4.208}  &6.7  &8.9e-6\\
		tai150b     & 2.889 & 2.889 & 170.9   &0 & 2.991  & 2.991  & 35.0  &  0  & 22.580  & 22.580    &  154.6  & {\bf 2.849} &{\bf 2.849}  
		&10.8  &1.8e-5\\
		tai256c     & { 1.219}  & { 1.219} & 181.1   &  6.2e-7 & {\bf 1.216}  & {\bf 1.216}  & 32.7  &0    & 203.637*  & --    &  317.5  & 1.272&1.272  &16.9  &7.9e-5\\
		tho150      & {\bf 1.872}  & {\bf 1.872}  & 44.8   &  1.7e-7 & 1.885  & 1.885  & 12.2  &  0  & 15.948  & 15.948    & 147.2   & 1.933&1.933  &6.1  &5.6e-5\\
		wil100      & { 0.694}  & { 0.694}  & 26.9   &  3.8e-6 & {\bf 0.686}  & {\bf 0.686}  & 8.3  &  0  & 8.162  & 8.162  & 52.6  & 0.768 &0.768  &1.5  &2.2e-5\\
		\hline
	\end{tabular*}
	{\small where, for ``tai256c'', EP4Orth+ cannot return a solution satisfying the orthogonality.}  
\end{table}

Table \ref{table2} and \ref{table3} report the levels of the average minimum and 
median gaps of four solvers for the 134 instances. We see that SEPPG+ and SEPPG0 yield 
better minimum and median gaps than ALM and EP4Orth+ do, and for half of instances 
their average minimum and median gaps are respectively less than $\textbf{0.5\%}$ 
and $\textbf{5\%}$. We also make a pairwise comparison for four solvers as in 
\cite{Dai05,JiangLW16}. Given a QAPLIB instance, solver $i$ is called the winner 
if the average minimum gap (or median gap) obtained by this solver is the smallest 
one among them. The second line of Table \ref{table4} reports the winners of each 
solver for the $134$ QAPLIB instances. From the third line, we perform pairwise 
comparisons of the solvers. The third and fourth lines indicate that SEPPG+ and SEPPG0 
are superior to ALM and EP4Orth+ whether in terms of the minimum gap or the median one. 
\begin{table}[h]
	\renewcommand\arraystretch{1.3}
	\setlength{\abovecaptionskip}{2pt}
	\setlength{\belowcaptionskip}{0pt}
	\centering
	\captionsetup{font={small}}
	\caption{Levels of the average minimum gaps on the 134 QAPLIB instances}\label{table2}
	\small
	\begin{tabular*}{\columnwidth}{@{\extracolsep\fill}cccccccccccccccc@{\extracolsep\fill}}
		\hline
		Min Gap$\le$\%&0.0&0.1& 0.2&0.3&0.4&0.5&0.6&0.7&0.8&0.9&1.0&2.0&3.0&4.0\\
		\hline	
		SEPPG+&{ 45}&{\bf 58}&{\bf 62}&{\bf 66}&{ 69}&{ 71}&{ 75}&{ 84}
		&{ 90}&{\bf 98}&{\bf 99}&{ 111}&{\bf 121}&{\bf 122}\\
		SEPPG0&{\bf 46}&{ 55}&{\bf 62}&{ 64}&{\bf 70}&{\bf 75}&{\bf 79}&{\bf 89}
		&{\bf 93}&{ 97}&{ 98}&{\bf 112}&{ 119}&{ 121}\\
		EP4Orth+&6&6&6&7&7&9&9&11&14&14&14&18&22&25\\
		ALM& 42&51&55&57&62&64&68&76&80&86&91&102&114&117\\
		\hline
	\end{tabular*}
\end{table}

\begin{table}[h]
	\renewcommand\arraystretch{1.3}
	\setlength{\abovecaptionskip}{2pt}
	\setlength{\belowcaptionskip}{0pt}
	\centering
	\captionsetup{font={small}}
	\caption{Levels of the average median gaps of on the 134 QAPLIB instances}\label{table3}
	\small
	\begin{tabular*}{\columnwidth}{@{\extracolsep\fill}cccccccccccccccc@{\extracolsep\fill}}
		\hline
		Med Gap$\le$\%&0.3&0.5& 0.7&1.0&3.0&5.0&7.0&10.0&15.0&20.0&25.0&30.0&40.0\\
		\hline
		SEPPG+&{\bf 4}&{ 8}&{ 15}&{ 23}&{ 55}&{\bf 86}
		&{ 93}&{ 101}&{\bf 114}&{\bf 117}&{ 119}
		&{\bf 122}&{\bf 122}\\
		SEPPG0&{\bf 4}&{\bf 9}&{\bf 17}&{\bf 25}&{\bf 58}&{\bf 86}
		&{\bf 96}&{\bf 102}&{\bf 114}&{\bf 117}&{\bf 121}
		&{\bf 122}&{\bf 122}\\
		EP4Orth+&0&0&0&0&10&15&17&23&44&53&73&84&91\\
		ALM&1&5&12&20&51&79&88&94&103&109&115&118&{\bf 122}\\
		\hline
	\end{tabular*}
\end{table}

\begin{table}[h]
	\renewcommand\arraystretch{1.3}
	\setlength{\abovecaptionskip}{2pt}
	\setlength{\belowcaptionskip}{0pt}
	\centering
	\captionsetup{font={small}}
	\caption{Pairwise comparison of four solvers on the 134 QAPLIB instances}\label{table4}
	\scriptsize
	\begin{tabular*}{\columnwidth}{@{\extracolsep\fill}rcccc|rcccc@{\extracolsep\fill}}
		\hline
		Min Gap&SEPPG+&SEPPG0&EP4Orth+&ALM&Med Gap&SEPPG+&SEPPG0&EP4Orth+&ALM\\
		\hline
		winners&73 &{\bf 75} &6 &66 & winners&{ 25} &{\bf 69} &0 &41 \\
		SEPPG+  &--&{\bf 93\,:\,90} &{\bf 134\,:\,6}& { 91\,:\,80} &SEPPG+  &--&39\,:\,96 & {\bf 134\,:\,0}&{ 70\,:\,64} \\
		SEPPG0  &90\,:\,93&--&133\,:\,7&  {\bf97\,:\,73} &SEPPG0  &{\bf96\,:\,39} &-- & {\bf 134\,:\,0}&{\bf 87\,:\,47} \\
		EP4Orth+ &6\,:\,134&7\,:\,133 & --&6\,:\,134 &EP4Orth+ &0\,:\,134&0\,:\,134 & --&0\,:\,134 \\
		ALM  &80\,:\,91& 73\,:\,97 &  {\bf134\,:\,6} &--&ALM    &64\,:\,70&47\,:\,87 &  {\bf134\,:\,0}&-- \\
		\hline
	\end{tabular*}
\end{table} 

\subsection{Graph matching problems}\label{sec6.4}

The CMU house image sequence (see \url{http://www.cs.cmu.edu/\textasciitilde cil/v-images.html})
is commonly used to test the performance of graph matching algorithms \cite{Anita20, Zhou15}.
This dataset consists of 111 frames of a house, each of which has been manually labeled with
30 landmarks. We connected the landmarks via Delaunay triangulation. Given an image pair
$(\mathcal{ G}_1,\mathcal{ G}_2)$, the edge-affinity matrix $K^q$ is computed by
$K_{c_ic_j}^{q}=\exp\big(-\!\frac{(q_{c_i}^{1}-q_{c_j}^{2})^2}{2500}\big)$ and
the node-affinity $K^p$ was set to zero, where $K_{c_ic_j}^{q}$ measures the similarity
between the $c_i$th edge of $\mathcal{ G}_1$ and the $c_j$th edge of $\mathcal{ G}_2$,
and the edge feature $q_{c}$ was computed as the pairwise distance between the connected nodes.
The graph matching problem is modelled as 
\begin{equation}\label{eq:gm_obj}
 \max_{X\in\mathcal{S}_{+}^{n,n}}{\rm vec}(X)^{\top}K{\rm vec}(X).
\end{equation}
Since the data matrix $K$ is nonnegative, we solve problem \eqref{eq:gm_obj},
rather than its lifted reformulation as in \eqref{QAP1}, with the four solvers
and compare their performance with that of FGM-D, a graph matching algorithm \cite{Zhou15}. 
The FGM-D is a path-following algorithm with a heuristic strategy, designed by 
the convex and concave relaxations for \eqref{orth-prob}. 

The starting points of the four solvers are chosen in the same way 
as in Section \ref{sec6.3}. For SEPPG+ and SEPPG0, we choose the initial $\rho_0$ in the same way as in Section \ref{sec6.3}, and for ALM, choose $\mu_0=0.1$, much better 
than $\mu_0=1$. For EP4Orth+, all parameters are default except 
$\textrm{sigma\_max}=10^{10}$ and $\textrm{maxiter}=2000$. 

\begin{figure}[h]
	\centering
	\includegraphics[width=\textwidth]{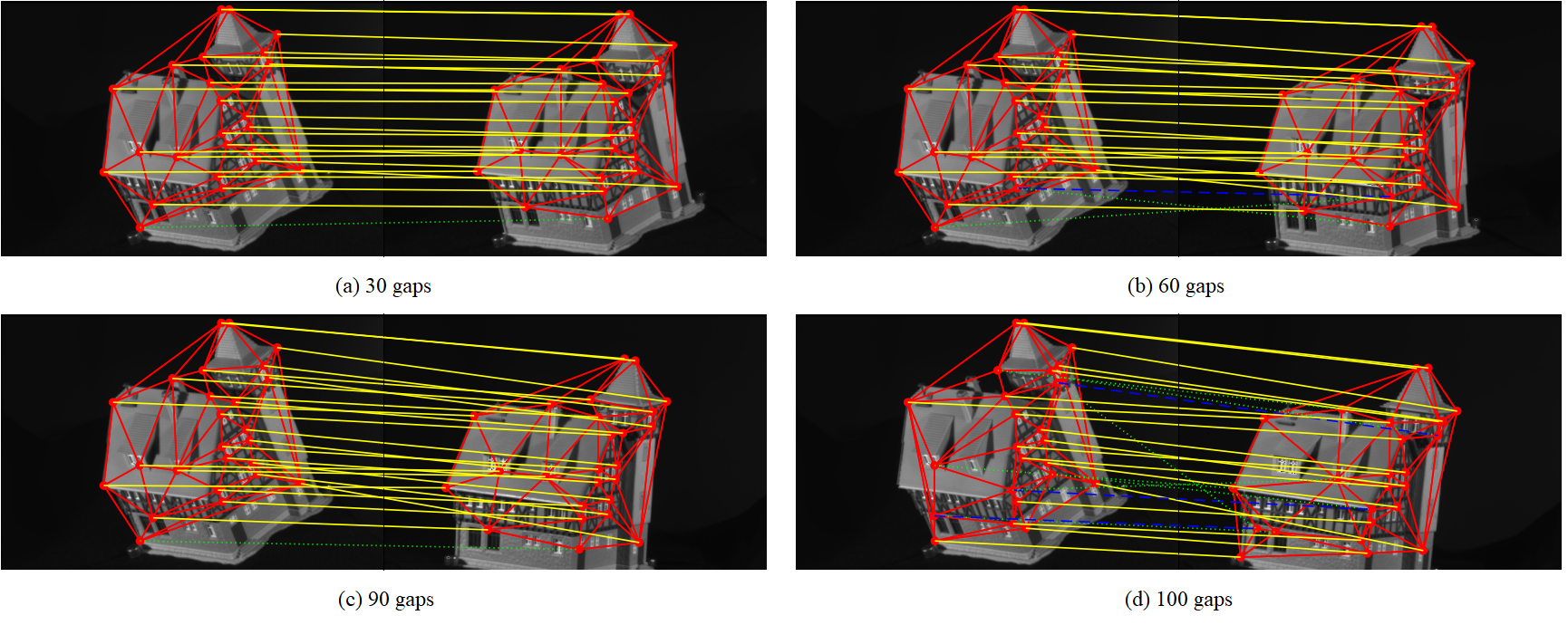}
	\caption{\small Example frame pairs for CMU house dataset generated by SEPPG+}
	\label{fig1}
\end{figure}

The sequence gap of a pair of images is defined as the number of images lying
between them, and the graph matching problems become more difficult in
the presence of a larger sequence gap. In the first experiment, we use image $\sharp 1$ to
match image $\sharp 30$, $\sharp 60$, $\sharp 90$, and $\sharp 100$, respectively,
and remove 5 out of the 30 landmarks for the purpose of noise simulation. 
The illustration (via SEPPG+) is shown in Figure \ref{fig1}, where the correct 
matches of the landmarks are drawn in yellow line, the incorrect ones are colored in blue, 
and the green lines label the missing landmarks.
The accuracy is defined by $\frac{n_{c}}{30-n_{m}}$ where $n_{c}$ denotes the number of
landmarks that are correctly matched, and $n_m$ denotes the number of missing landmarks.
Figure \ref{fig1} displays the matching effect of SEPPG+ in terms of the accuracy corresponding
to the best objective value yielded by running $10$ times from different starting points.
We see that SEPPG+ can provide a satisfactory matching by running $10$ times randomly 
(with accuracy attaining $92\%$).

\begin{figure} 
\centering
    \subfloat{%
        \includegraphics[width=\textwidth]{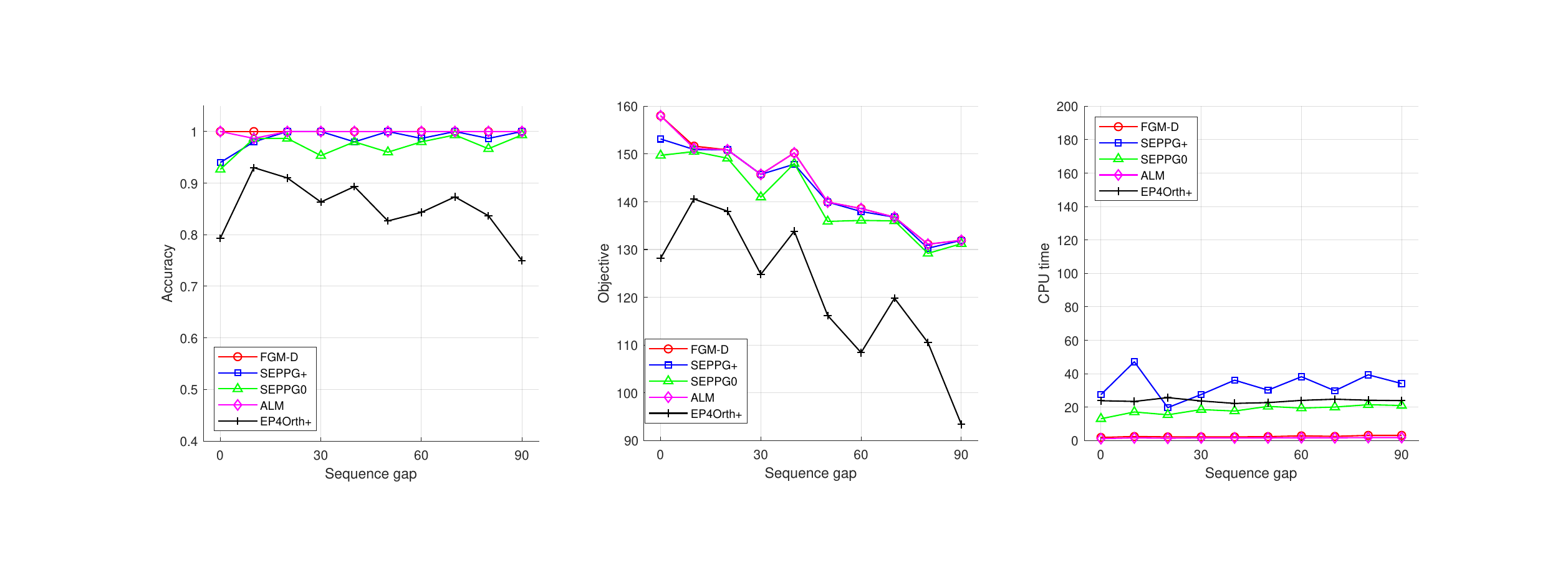}
    }

    \vspace{-1cm}
    \subfloat{%
        \includegraphics[width=\textwidth]{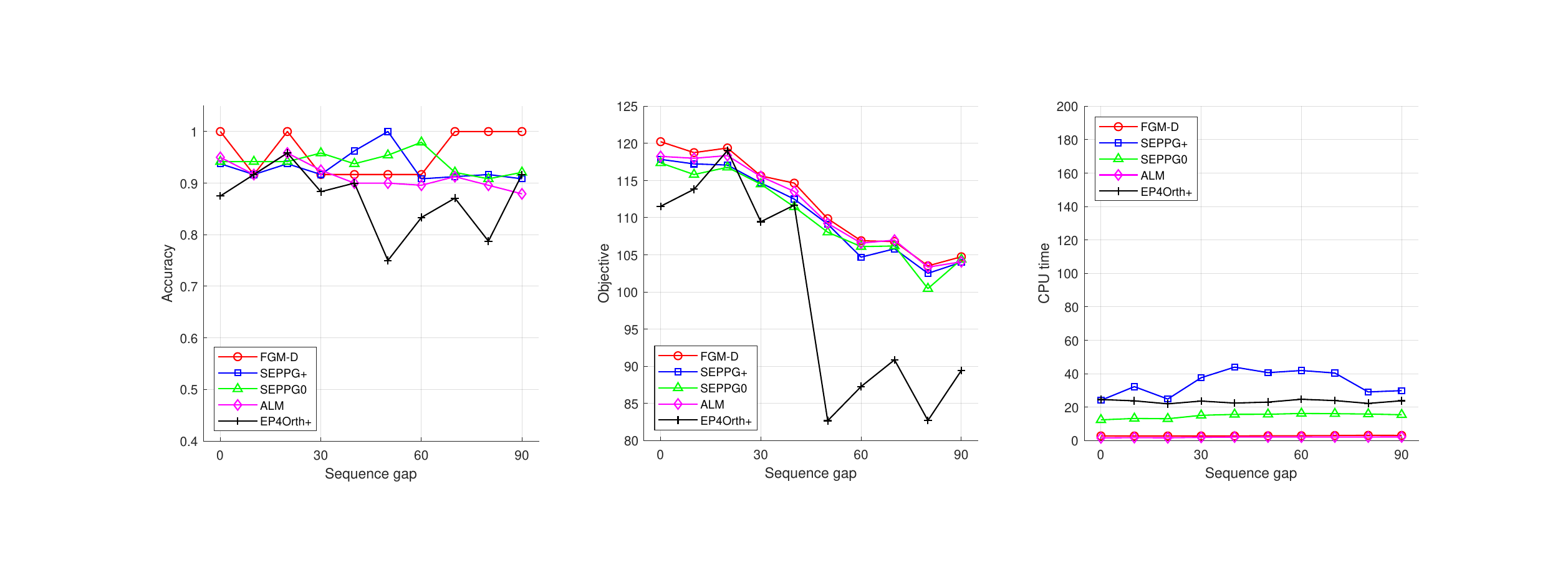}
    }
	\caption{\small Performance of solvers on (a) graph pairs with perfect 30 nodes (b) graph pairs with 5 nodes removed randomly}
	\label{fig2}
\end{figure}

Next we undertake the experiments in two scenarios. In the first scenario,
we match a sequence of original image pairs with different sequence gap varying from 0 to 90.
Figure \ref{fig2} (a) displays the performance of FGM-D and four solvers in terms of
the accuracy and the objective values of \eqref{eq:gm_obj}, where the accuracy and
the objective values take the average for $10$ times running from the starting points 
generated randomly. We see that SEPPG+, SEPPG0 and ALM yield the comparable accuracy 
and objective values, which are close to those yielded by FGM-D. 
In the second scenario, we simulate noisy images by randomly removing $5$ landmarks 
from each image and test the performance on matching the sequence
of noisy image pairs. Figure \ref{fig2} (b) displays the performance of FGM-D and
four solvers in terms of the average matching accuracy and objective value of \eqref{eq:gm_obj} for $10$ times running. We see that in terms of the accuracy 
and objective value, SEPPG+, SEPPG0 and ALM have the comparable performance which are 
remarkably superior to EP4Orth+. In addition, the objective values yielded by 
SEPPG+, SEPPG0 and ALM are very close to the one yielded by FGM-D though their accuracy 
is a little less favorable than the one yielded by FGM-D.  
\subsection{Projection onto nonnegative orthogonal set}\label{sec6.5}

Given a matrix $C\in\mathbb{R}^{n\times r}$, we consider the problem of computing its projection onto $\mathcal{S}_{+}^{n,r}$:
\begin{align}\label{projprob}
	\min_{X\in\mathcal{S}_{+}^{n,r}}\|X-C\|_F^2.
\end{align}
We evaluate a solver in terms of the objective value (obj), the infeasibility (infeas) 
and the CPU time in second, where the infeasibility is measured 
by $\|(X^{f})\tp X^{f}-I_{r}\|_F+\|\min(0,X^{f})\|_F$. 

The starting points of the four solvers are obtained by imposing the rounding technique from EP4Orth+ on the matrix $C$. For SEPPG+, SEPPG0 and ALM, we choose $\rho_0=1/\|C\|$ and $\mu_0=1/\|C\|$. For EP4Orth+, all parameters are chosen in the same way as in \citep[Section 5.1]{JiangM22} except $\gamma_2=1.05$, which is much better than $\gamma_2=5$ for the first experiment below.

In the first experiment, for each $(n,r)$, we generate $50$ binary matrices 
$C\in\{0,1\}^{n\times r}$ by the MATLAB command $\textrm{randi}(2,n,r)-1$ and 
apply the four solvers to solving the $50$ instances. Table \ref{table5} 
reports the averaged results of the objective values, the infeasibility
and the CPU time. We see that for such $C$, the four solvers yield 
the comparable objective values, and EP4Orth+ yields a little worse objective 
values but low infeasibility among the four solvers.
\begin{table}[h]
	\renewcommand\arraystretch{1.3}
	\setlength{\abovecaptionskip}{2pt}
	\setlength{\belowcaptionskip}{0pt}
	\centering
	\captionsetup{font={small}}
	\caption{Numerical results for the projection problem of 0-1 matrix $C$ 
		onto $\mathcal{S}_{+}^{n,r}$}\label{table5}
	\small
	\begin{tabular*}{\textwidth}{@{\extracolsep{\fill}}lllllllllllllll@{\extracolsep{\fill}}}
		\hline
		&\multicolumn{3}{l}{SEPPG+} &\multicolumn{3}{l}{SEPPG0} & \multicolumn{3}{l}{EP4Orth+}&\multicolumn{3}{l}{ALM}\\
		\cmidrule(lr){2-4} \cmidrule(lr){5-7} \cmidrule(lr){8-10} \cmidrule(lr){11-13}
		$(n,r)$ & obj   &infeas  &time & obj   &infeas  &time  & obj   &infeas  &time  &obj   &infeas  &time    \\
		\hline
		(300,5)  & 677.5  & 3.6e-8 & 0.9 & {\bf 677.3} & 9.4e-16& 2.5& 677.5  & 1.6e-15  & 0.0 & {\bf 677.3}  & 1.2e-7 & 0.1 \\
		(300,10) & 1399.1 & 4.6e-8 & 2.1 & {\bf 1398.8}& 5.4e-16& 5.4& 1400.2 & 1.3e-15  & 0.1 & {\bf 1398.8} & 1.8e-7 & 0.2 \\
		(500,5)  & 1152.2 & 2.7e-8 & 1.2 & 1152.0& 6.9e-16& 3.2& 1152.3 & 1.2e-15 & 0.1 & {\bf 1151.9} & 9.6e-8 & 0.1 \\
		(500,10) & 2367.9 & 4.5e-8 & 3.0 & {\bf 2367.7}& 4.5e-16& 7.2& 2369.4 & 1.1e-15  & 0.2 & {\bf 2367.7} & 1.4e-7 & 0.4 \\
		(1000,5) & 2365.5 & 2.1e-8 & 2.5 & {\bf 2365.0}& 1.0e-15& 6.1& 2365.4 & 1.7e-15  & 0.1 & {\bf 2365.0} & 6.4e-8 & 0.2 \\
		(1000,10)& 4806.0 & 2.8e-8 & 7.6 & 4805.9& 5.2e-16& 12.0& 4808.5 & 1.4e-15  & 0.4 & {\bf 4805.8} & 9.2e-8 & 0.9 \\
		\hline
	\end{tabular*}
\end{table}

In the second experiment, we generate the matrix $C$ in the same way as in \cite{JiangM22}. 
The parameter $\xi\in[0,1]$ is used to control the magnitude of noise level 
and a larger $\xi$ makes it harder to find the ground truth $X^*$. For such $C$, 
the projection problem \eqref{projprob} has a unique global optimal solution. 
For each $(\xi,n,r)$, we generate $50$ such matrices $C\in\mathbb{R}^{n\times r}$ 
and then use the four solvers to solve them. Table \ref{table6} reports the averaged 
results of the gap, the infeasibility and the CPU time, 
where $\textrm{gap}=\frac{\|X^{f}-C\|_F}{\|X^*-C\|_F}-1$. We see that 
for such $C$, EP4Orth+ yields the best gap and infeasibility, while 
SEPPG+, SEPPG0 and ALM yield the comparable infeasibility.  
\begin{table}[h]
	\renewcommand\arraystretch{1.3}
	\setlength{\abovecaptionskip}{2pt}
	\setlength{\belowcaptionskip}{0pt}
	\centering
	\captionsetup{font={small}}
	\caption{Numerical results for the projection problems onto $\mathcal{S}_{+}^{n,r}$ generated in \cite{JiangM22}}\label{table6}
	\small
	\begin{tabular*}{\textwidth}{@{\extracolsep{\fill}}lllllllllllllll@{\extracolsep{\fill}}}
		\hline
		&\multicolumn{3}{l}{SEPPG+} &\multicolumn{3}{l}{SEPPG0} & \multicolumn{3}{l}{EP4Orth+}&\multicolumn{3}{l}{ALM}\\
		\cmidrule(lr){2-4} \cmidrule(lr){5-7} \cmidrule(lr){8-10} \cmidrule(lr){11-13}
		$(\xi,n,r)$ & gap   &infeas  &time & gap   &infeas  &time  & gap   &infeas  &time  &gap   &infeas  &time    \\
		\hline
		(0.5,500,5)   & 6.2e-9 & 6.5e-8 & 0.1 & 3.8e-9& 4.1e-8& 0.1 & {\bf 4.4e-18} & 3.5e-16 & 0.0 & 2.2e-9 & 4.9e-8 & 0.0 \\
		(0.8,500,5)   & 2.1e-5 & 5.3e-8 & 0.2 & 4.1e-6& 6.3e-8& 0.2 & {\bf 4.4e-18} & 3.7e-16 & 0.1 & 2.6e-9 & 5.7e-8 & 0.0 \\
		(1.0,500,5)   & 5.4e-4 & 4.3e-8 & 0.4 & 2.6e-5& 5.1e-8& 0.6 & 2.8e-4 & 3.6e-16 & 0.1 & {\bf 9.7e-6} & 9.7e-8 & 0.1 \\
		(0.5,1000,10) & 9.4e-9 & 5.7e-8 & 0.5 & 3.5e-9& 5.6e-8& 0.6 & {\bf 1.3e-17} & 4.9e-16 & 0.3 & 9.3e-10 & 4.1e-8 & 0.2 \\
		(0.8,1000,10) & 3.5e-6 & 3.5e-8 & 2.7 & 3.2e-6& 5.3e-8& 3.4 & {\bf 8.9e-18}  & 4.7e-16 & 0.5 & 2.0e-6 & 1.2e-7 & 0.6 \\
		(1.0,1000,10) & 1.3e-3 & 3.0e-8 & 3.8 & 1.1e-3& 5.1e-16& 10.1 & 4.1e-4  & 4.4e-16 & 0.5 & {\bf 6.8e-5} & 7.7e-8 & 1.4 \\
		\hline
	\end{tabular*}
\end{table}
\subsection{Orthogonal nonnegative matrix factorization}\label{sec6.6}

Given a matrix $A\in\mathbb{R}_+^{n\times p}$, the orthogonal 
nonnegative matrix factorization (ONMF) problem is 
\begin{align}\label{ONMF}
	\min_{X\in\mathcal{S}^{n,r}_+,Y\in\mathbb{R}_+^{p\times r}}
	\|A-XY^\top\|_F^2.
\end{align}
Similar to \citep{JiangM22}, we use the alternating minimization method to 
solve \eqref{ONMF}, and when the variable $X$ is fixed, take 
${\rm Proj}_{\mathbb{R}_{+}^{p\times r}}(A\tp X(X\tp X)^{-1})$ as 
an approximate solution to the convex problem on $Y$. 
We evaluate the four solvers on text and image datasets from \cite{Cai08}, 
which are available at \url{http://www.cad.zju.edu.cn/home/dengcai/Data/data.html}. 
The original data is preprocessed by using the same way as in \cite{JiangM22}. 
We evaluate the performance of clustering results in terms of three criteria: 
purity \cite{Ding06}, entropy \cite{ZhaoK04} and normalized mutual 
information (NMI) \cite{XuLG03}. Purity gives a measure of the predominance 
of the largest category per cluster, and larger purity means better clustering results. 
In addition, smaller entropy and larger NMI also imply better clustering results.
Let $r$ be the number of clusters and let $n$ be the total number of data points. 
Suppose that $\mathcal{C}=\cup_{i=1}^r\mathcal{C}_i$ 
and $\mathcal{C}'=\cup_{i=1}^r\mathcal{C}_i'$ are the clustering results  
yielded by the ground truth and a certain solver, respectively.
The purity, entropy and NMI are respectively computed by 
\[
\textrm{Pidx}:=\frac{\sum_{i=1}^r\max_j\{n_{ji}\}}{n},\
\textrm{Eidx}:=-\frac{\sum_{j=1}^r\sum_{i=1}^rn_{ij}
	\log_2\frac{n_{ij}}{n_j'}}{n\log_2 r},\ 
\textrm{NMI}:=\frac{\sum_{i=1}^r\sum_{j=1}^r\frac{n_{ij}}{n}
	\log_2\frac{nn_{ij}}{n_in_j'}}{\max(H(\mathcal{C}),H(\mathcal{C}'))}
\] 
where $H(\mathcal{C})=-\sum_{i=1}^r\frac{n_i}{n}\log_2\frac{n_i}{n}$ and 
$H(\mathcal{C'})=-\sum_{i=1}^r\frac{n_i'}{n}\log_2\frac{n_i'}{n}$ with 
$n_i=|\mathcal{C}_i|,n_i'=|\mathcal{C}_i'|$ and $n_{ij}=|\mathcal{C}_i\cap\mathcal{C}_j'|$. 
The starting points of SEPPG+, SEPPG0 and ALM are generated in the same way as 
for EP4Orth+, which depends on the singular value vectors corresponding to 
the first $r$ largest singular values of $A$. For SEPPG+, SEPPG0 and ALM, we choose $\rho_0=1/\|A\|$ and $\mu_0=1/\|A\|$. For EP4Orth+, all parameters are set to be the same as those in \citep[Section 5.2.1]{JiangM22}.
\begin{table}[h]
	\renewcommand\arraystretch{1.3}
	\setlength{\abovecaptionskip}{2pt}
	\setlength{\belowcaptionskip}{0pt}
	\centering
	\captionsetup{font={small}}
	\caption{Numerical results of four solvers on real datasets}\label{table7}
	\tiny
		\begin{tabular*}{\textwidth}{@{\extracolsep{\fill}}cllllllllllllllllllllll@{\extracolsep{\fill}}}
			\hline
			Name &\multicolumn{5}{l}{SEPPG+}&\multicolumn{5}{l}{SEPPG0} & \multicolumn{5}{l}{EP4Orth+}&\multicolumn{5}{l}{ALM}\\
			\cmidrule(lr){2-6} \cmidrule(lr){7-11} \cmidrule(lr){12-16} \cmidrule(lr){17-21}
			& Pidx  & NMI  & Eidx  &infeasi  &time  &  Pidx    &NMI   & Eidx & infeasi  & time  
			& Pidx  & NMI  & Eidx  &infeasi  &time  & Pidx  & NMI  & Eidx  &infeasi  &time   \\
			& (\%)  &(\%)   &(\%)  &  &  & (\%)  &(\%)   &(\%)  &  &  &  (\%)  &(\%)   &(\%)  &  &  &  (\%)  &(\%)   &(\%)  &  &    \\
			\hline
			MNIST       & {\bf 60.7}  & {\bf 49.4}  & {\bf 50.6}  & 8.7e-16  & 23.3 & 60.3&48.7 &51.3 & 1.1e-15& 34.0 & 60.3  & 48.9  & 51.0  & 9.6e-16 & 19.7  & 55.5  & 47.0  & 52.9  & 7.2e-16  & 9.5   \\
			Yale        & {\bf 46.1}  & { 46.9}  & { 53.1}  & 9.1e-16  & 2.6 &44.2 &{\bf47.1} &{\bf52.9} & 3.7e-16&  3.9 & 44.8  & 46.1  & 53.9  & 4.6e-16 & 1.3  & 39.4  & 42.4  & 57.6  & 8.5e-16  & 1.5   \\
			TDT2-l10    & 82.7  & 77.0  & 22.9  & 7.2e-16  & 7.2  &81.0 &75.2 &24.8 & 3.5e-16& 27.3 & {\bf 84.5}  & {\bf 79.9}  
			& {\bf 20.0}  & 7.7e-16 & 4.8  & 82.8  & 77.4  & 22.5  & 1.0e-15  & 12.9   \\
			TDT2-l20    & {\bf 83.1}  & 83.8  & 15.9  & 1.1e-15  & 36.9 &82.8 &83.2 &16.4 & 1.0e-15&  124.9& {\bf 83.1}  & {\bf 84.2}  & {\bf 15.5}  & 9.9e-16 & 22.5  & {\bf 83.1}  &  83.9  & 15.8  & 1.4e-15  & 56.2   \\
			TDT2-t10    & {\bf 85.8}  & {\bf 70.2}  & 20.8  & 8.9e-16  & 20.4 &83.5 &69.3 &{\bf19.6} & 9.9e-16&  81.8& 85.6  & 69.9  & 20.9  & 1.1e-15 & 14.0  & {\bf 85.8}  & 70.1  & { 20.7}  & 1.0e-15  & 29.3   \\
			TDT2-t20    & 81.6  & 67.6  & 19.7  & 1.7e-15  & 40.5 &81.5 &66.8 &19.7 & 7.2e-16& 124.6 & 82.3  & {\bf 69.5}  & 18.1  & 1.0e-15 & 27.8  & {\bf 82.7}  &  69.1  & {\bf 17.9}  & 1.4e-15  & 56.3   \\
			Ret-t10     & {\bf 74.3}  & {\bf 62.3}  & {\bf 36.4}  & 1.2e-15  & 18.5 &72.3 &59.9 &38.7 & 1.1e-15& 68.6 & 71.1  & 58.9  & 39.7  & 6.5e-16 & 12.7  & 74.0  &  61.8  & 36.9  & 1.0e-15  & 16.3   \\
			Ret-t20     & { 65.6}  & 56.6  & 38.1  & 1.5e-15  & 35.9 &{\bf65.8} &56.6 &{\bf37.7} & 1.3e-15& 113.6 & 65.4  & 56.2  & 38.5  & 7.6e-16 & 27.1  & 65.5  & {\bf 56.7}  & { 38.0}  & 1.5e-15  & 54.8   \\
			NewsG-t5    & { 41.6}  & {\bf 22.9}  & {\bf 77.0}  & 1.4e-15  & 16.4 &{\bf42.6} &22.5 &77.4 & 2.7e-16& 75.6 & 41.4  & 22.7  & 77.1  & 7.0e-16 & 10.9  & { 41.6}  & 22.7  & 77.2  & 3.8e-16  & 17.1   \\
			\hline
		\end{tabular*}
	\end{table}
	
	Table \ref{table7} reports the text and images clustering results yielded by 
	the four solvers. We see that the purity, entropy and NMI yielded by SEPPG+ 
	are a little better than those yielded by SEPPG0 and EP4Orth+, which are better than 
	those yielded by ALM. 
	
	From the previous numerical comparisons, we conclude that SEPPG+ and SEPPG0 have better 
	performance than EP4Orth+ does in terms of the quality of solutions 
	especially for QAPs and graph matching problems, and SEPPG+ and SEPPG0 have better 
	performance than ALM does in terms of the quality of solutions for 
	QAPs and orthogonal nonnegative matrix factorization. Among the four solvers, 
	ALM requires the least time for QAPs and graph matching, while EP4Orth+ requires 
	the least time for projection problems and the orthogonal nonnegative matrix 
	factorization.
 \section{Conclusions}\label{sec7}

 For the nonnegative orthogonal constraint system, we have derived a locally 
 Lipschitzian error bound under the no zero assumption for $n>r>1$, 
 and the global Lipschitzian error bound for $n>r=1$ or $n=r$. Furthermore, under the assumption that every global minimizer has no zero rows when $n>r>1$,
 we established that problem \eqref{Nsmooth-pen} is a global exact 
 penalty of \eqref{orth-prob}, and so are problems \eqref{smooth-pen} and \eqref{Quad-Pen}
 for those $f$ satisfying the lower second-order calmness in \eqref{growth}. 
 Example \ref{example1} illustrates that the no zero row assumption
 for every global minimizer when $n>r>1$ cannot be removed. 
 A practical penalty algorithm was developed by solving approximately 
 a series of smooth penalty problems \eqref{smooth-pen} or \eqref{Quad-Pen}, and any cluster point of the generated sequence is shown to be a stationary point of \eqref{orth-prob}. Numerical comparisons with ALM \cite{Wen13} and the exact penalty method \citep{JiangM22} on several classes of test problems demonstrate that
 our penalty method has an advantage in terms of the quality of solutions 
 especially for those difficult QAPs, though it requires a little more time than the others do. 
 
 \appendix
 
 \section{To bound $\|X\tp X-I\|_F$ by ${\rm dist}(X,{\rm St}(n,r))$}
 \begin{lemma}\label{Bound-relation}
  For any $X\in\mathbb{R}^{n\times r}$, 
 	${\rm dist}(X,{\rm St}(n,r))\le\|X\tp X-I\|_F
 	\le(1+\|X\|){\rm dist}(X,{\rm St}(n,r))$.
 \end{lemma}
 \begin{proof}
  Fix any $X\in\mathbb{R}^{n\times r}$. Let $X$ have the SVD given by 
  $X=P\left[\begin{matrix}
  	 \Sigma \\ 0
  	 \end{matrix}\right]Q^{\top}$, where $\Sigma={\rm diag}(\sigma_1,\ldots,\sigma_r)$ with $\sigma_1\ge\sigma_2\ge\cdots\ge\sigma_r$, 
  $P$ is an $n\times n$ orthogonal matrix, and $Q$ is a $r\times r$ orthogonal matrix. 
  Let $P_1$ be the matrix consisting of the first $r$ columns of $P$. Then, 
  ${\rm dist}(X,{\rm St}(n,r))=\|X-P_1Q^{\top}\|_F=\|P_1\Sigma Q^{\top}-P_1Q^{\top}\|_F=\sqrt{\sum_{i=1}^r(\sigma_i-1)^2}$. On the other hand,
  from the SVD of $X$, $\|X\tp X-I\|_F=\sqrt{\sum_{i=1}^r(\sigma_i^2-1)^2}$.
  The two sides imply that the desired result holds.
 \end{proof}

 \section{A second-order sufficient condition of problem \eqref{orth-prob}}\label{appendix_A}
 
 \begin{proposition}\label{SSOC}
 	Let $\widetilde{f}=f+\delta_{\mathcal{S}_{+}^{n,r}}$ be the extended valued objective   
 	function of \eqref{orth-prob}. Suppose that $f$ is twice continuously differentiable 
 	on $\mathcal{O}$. Then, $\widetilde{f}$ is properly twice epi-differentiable at 
 	any $X\in\mathcal{S}_{+}^{n,r}$ for $V\in\partial\!\widetilde{f}(X)$, 
 	and any stationary point $\overline{X}\in\mathcal{S}_{+}^{n,r}$ satisfying 
 	the metric subregularity constraint qualification (MSCQ) is a strong local optimal 
 	solution of \eqref{orth-prob} if for any nonzero $H\in\mathcal{C}(\overline{X})$ 
 	with $-\overline{X}H\tp H\in\mathcal{T}_{\mathbb{R}_{+}^{n\times r}}(\overline{X})$, 
 	where $\mathcal{C}(\overline{X}):=\big\{H\in\mathbb{R}^{n\times r}\ |\ 
 	\mathcal{A}_{\overline{X}}(H)=0,\langle H,\nabla\!f(\overline{X})\rangle=0\big\}$, 
 	\begin{equation}\label{SOSC-ineq}
 		\langle H,\nabla^2\!f(\overline{X})H\rangle
 		-\langle H\tp H,\overline{X}\tp\nabla\!f(\overline{X})\rangle>0.
 	\end{equation}
 \end{proposition} 
 \begin{proof}
 	Note that $\mathcal{S}_{+}^{n,r}=F^{-1}(\mathbb{R}_{+}^{n\times r}\times\{0\})$ 
 	with $F(X)=(X;X^{\top}X-I_{r})$ for $X\in\mathbb{R}^{n\times r}$.  
 	From the given assumption, by invoking \citep[Theorem 5.6]{Mordu21} 
 	with $\Theta=\mathbb{R}_{+}^{n\times r}\times\{0\}$ and $\Omega=\mathcal{S}_{+}^{n,r}$, 
 	it follows that
 	\begin{equation}\label{subderive-equa0}
 		d^2\widetilde{f}(\overline{X}|0)
 		=\left\{\begin{array}{cl}
 			\langle H,\nabla^2\!f(\overline{X})H\rangle
 			+2{\displaystyle\max_{Y\in\mathcal{Y}(\overline{X},-\!\nabla\!f(\overline{X}))}}\langle H,HY\rangle, &{\rm if}\ H\in\mathcal{C}(\overline{X}),\\
 			+\infty, & {\rm otherwise,}
 		\end{array}\right.
 	\end{equation}
 	where $\mathcal{Y}(\overline{X},-\!\nabla\!f(\overline{X}))
 	=\big\{Y\in\mathbb{S}^r\,|\,-\!\nabla\!f(\overline{X})-\mathcal{A}_{\overline{X}}^*(Y)
 	\in\mathcal{N}_{\mathbb{R}_{+}^{n\times r}}(\overline{X})\big\}$. 
 	For each $Y\in\mathcal{Y}(\overline{X},-\!\nabla\!f(\overline{X}))$, 
 	obviously, there exists $U\in\mathcal{N}_{\mathbb{R}_{+}^{n\times r}}(\overline{X})$ 
 	such that $U=-\!\nabla\!f(\overline{X})-\mathcal{A}_{\overline{X}}^*(Y)
 	=-\!\nabla\!f(\overline{X})-2\overline{X}Y$. Together with 
 	$\langle H,HY\rangle=\langle H,H\overline{X}\tp\overline{X}Y\rangle$, 
 	it immediately follows that 
 	\begin{align*}
 		2\max_{Y\in\mathcal{Y}(\overline{X},-\!\nabla\!f(\overline{X}))}\langle H,HY\rangle
 		&=\max_{U\in\mathcal{N}_{\mathbb{R}_{+}^{n\times r}}(\overline{X})}
 		-\langle H,H\overline{X}\tp U\rangle-\langle H,H\overline{X}\tp\nabla\!f(\overline{X})\rangle\\
 		&=\delta_{\mathcal{T}_{\mathbb{R}_{+}^{n\times r}}(\overline{X})}
 		(-\overline{X}H\tp H)-\langle H,H\overline{X}\tp\!\nabla\!f(\overline{X})\rangle.
 	\end{align*} 
 	Combining this with \eqref{subderive-equa0} and invoking \cite[Theorem 13.24]{RW98}, 
 	we obtain the second part. The first part is immediate by \citep[Corollary 5.11]{Mordu21} 
 	and the polyhedral convexity of $\mathbb{R}_{+}^{n\times r}\times\{0\}$.
 \end{proof}
 \begin{remark}\label{remark-SOSC}
 	\begin{enumerate}[1.]
 		\item By Proposition \ref{prop-Stnr} and \cite[Section 3.1]{Ioffe08}, if a stationary point 
 		$\overline{X}\in\mathcal{S}_{+}^{n,r}$ has no zero rows when $n>r$, 
 		then it necessarily satisfies the MSCQ.
 		
 		\item We state that $\mathcal{C}(\overline{X})$ coincides with the critical cone 
 		defined in \citep[Eq.\,(2.12)]{JiangM22}. From \citep[Lemma 2.2]{JiangM22}, 
 		$\mathcal{T}_{\mathbb{R}_{+}^{n\times r}}\cap{\rm Ker}\mathcal{A}_{\overline{X}}$ 
 		is the linearized cone in \citep[Eq.\,(2.12)]{JiangM22}, so it suffices to argue 
 		that for every $H\in\mathcal{C}(\overline{X})$, 
 		$\langle H,\nabla\!f(\overline{X})\rangle=0$ if and only if $H_{ij}=0$ when 
 		$[\nabla\!f(\overline{X})]_{ij}>0$ for all $(i,j)\in\overline{J}_2(\overline{X})
 		:=\{(i,j)\in\overline{{\rm supp}}(\overline{X})\,|\,\|\overline{X}_{i\cdot}\|=0\}$.
 		Fix any $H\in\mathcal{C}(\overline{X})$. Since $\overline{X}$ is a stationary point 
 		of \eqref{orth-prob}, by using Definition \ref{spoint} and noting that 
 		$H\in{\rm T}_{\overline{X}}M$ and $\langle H,\nabla\!f(\overline{X})\rangle=0$,
 		we deduce that $\langle H,G\rangle=0$ for every 
 		$G\in\mathcal{N}_{\mathbb{R}_{+}^{n\times r}}(\overline{X})$, 
 		which along with $H\in\mathcal{T}_{\mathbb{R}_{+}^{n\times r}}(\overline{X})$ 
 		means that $H_{ij}G_{ij}=0$ for every 
 		$G\in\mathcal{N}_{\mathbb{R}_{+}^{n\times r}}(\overline{X})$ and every $(i,j)$.
 		Note that $0\in\nabla\!f(\overline{X})+\overline{X}S+G$ for some $S\in\mathbb{S}^r$
 		and $G\in\mathcal{N}_{\mathbb{R}_{+}^{n\times r}}(\overline{X})$ and 
 		$[XS]_{ij}=0$ for $(i,j)\in\overline{J}_2(\overline{X})$. Then the stated fact holds.

 		\item Comparing with the second-order sufficient condition in \citep[Theorem 2.3]{JiangM22},   
 		we see that the term
 		$-\langle H\tp H,{\rm off}(\overline{X}\tp\!\nabla\!f(\overline{X}))\rangle$ in 
 		\eqref{SOSC-ineq} corresponds to 
 		$-\langle H\tp H,{\rm off}(\overline{\lambda}VV\tp)\rangle$ in 
 		\citep[Eq.\,(2.16)]{JiangM22}, where 
 		$\overline{\lambda}\in\mathbb{R}$ is the Lagrange multiplier of 
 		the nonconvex equality constraint in \citep[Eq.\,(1.8)]{JiangM22}.
 		Clearly, our second-order sufficient condition does not involve such a multiplier.
 		Instead, our second-order sufficient condition involves those nonzero directions 
 		from $\mathcal{C}(\overline{X})\cap\{H\in\mathbb{R}^{n\times r}\ |\
 		-\overline{X}H\tp H\in\mathcal{T}_{\mathbb{R}_{+}^{n\times r}}(\overline{X})\}$,
 		which is a subset of the critical cone $\mathcal{C}(\overline{X})$.        
 	\end{enumerate} 
 \end{remark}

\section{The proof of Theorem \ref{theorem1-ManPG1}}\label{appendix_B}
\begin{proof}
	{\bf Item 1.} By the definition of $\ell(k)$, it is immediate to check that 
	for each $k\in\mathbb{N}$, 
	\begin{align*}
		\Theta_{\rho,\gamma}(X^{\ell(k+1)})
		&=\max\Big\{\max_{j=1,\ldots,\min(m,k+1)}\Theta_{\rho,\gamma}(X^{k+1-j}),\
		\Theta_{\rho,\gamma}(X^{k+1})\Big\}\\
		&\le\max\Big\{\Theta_{\rho,\gamma}(X^{\ell(k)}),\ 
		\Theta_{\rho,\gamma}(X^{\ell(k)})-\frac{\alpha}{2t_k}\|V^k\|_F^2\Big\}
		\le\Theta_{\rho,\gamma}(X^{\ell(k)}),
	\end{align*}
	where the first inequality is due to the criterion in Step \ref{Astep5}.
	This means that the sequence $\{\Theta_{\rho,\gamma}(X^{\ell(k)})\}_{k\in\mathbb{N}}$
	is monotonically decreasing. Together with the lower boundedness of $\Theta_{\rho,\gamma}$ 
	on ${\rm St}(n,r)$, the sequence $\{\Theta_{\rho,\gamma}(X^{\ell(k)})\}_{k\in\mathbb{N}}$ 
	is convergent, say, $\lim_{k\to\infty}\Theta_{\rho,\gamma}(X^{\ell(k)})
	=\varpi_{\!\rho,\gamma}^*$. In addition, by following the same arguments as those 
	for \cite[Lemma 4]{Wright09}, we can obtain $\lim_{k\to\infty}V^{\ell(k)-j}=0$ 
	for all $j\ge 1$. Notice that $\ell(k)$ is one of the indices $k\!-\!m,\ldots,k$. 
	Hence, $k-(m+1)=\ell(k)-j$ for some $j\in\{1,2,\ldots,m+1\}$. Consequently, $0=\lim_{k\to\infty}V^{k-(m+1)}=\lim_{k\to\infty}V^k$.
	Notice that $X^{\ell(k)}=R_{X^{\ell(k)-1}}(V^{\ell(k)-1})$ for each $k\in\mathbb{N}$.
	From $\lim_{k\to\infty}V^{\ell(k)-j}=0$ and Remark \ref{remark-retraction},
	for all sufficiently large $k$, 
	$X^{\ell(k)}=X^{\ell(k)-1}+V^{\ell(k)-1}+o(\|V^{\ell(k)-1}\|_F)$, 
	which implies that
	$ X^{\ell(k)}=X^{k-(m+1)}+{\textstyle\sum_{j=1}^{\ell(k)-(k-m-1)}}
	\!\big[V^{\ell(k)-j}+o(\|V^{\ell(k)-j}\|_F)\big]$.
	Thus, we obtain $\lim_{k\to\infty}(X^{\ell(k)}-X^{k-(m+1)})=0$. Along with the continuity 
	of $\Theta_{\rho,\gamma}$, $\lim_{k\to\infty}\Theta_{\rho,\gamma}(X^k)
	=\lim_{k\to\infty}\Theta_{\rho,\gamma}(X^{k-m-1})
	=\lim_{k\to\infty}\Theta_{\rho,\gamma}(X^{\ell(k)})=\varpi_{\!\rho,\gamma}^*$.
	
	\medskip
	
	{\bf Item 2.} The boundedness of $\{X^k\}_{k\in\mathbb{N}}$ is immediate due to
	$\{X^k\}_{k\in\mathbb{N}}\subset{\rm St}(n,r)$. Let $\overline{X}$ be an arbitrary
	cluster point of $\{X^k\}_{k\in\mathbb{N}}$. Then, there exists an index set 
	$\mathcal{K}\subset\mathbb{N}$ such that 
	$\lim_{\mathcal{K}\ni k\to\infty}X^k=\overline{X}$. Together with Item 1 
	and $t_k\ge t_{\rm min}>0$, from \eqref{Proj-St} it follows that
	\(
	0=\!\lim_{\mathcal{K}\ni k\to\infty}{\rm grad}\Theta_{\rho,\gamma}(X^k)
	={\rm grad}\Theta_{\rho,\gamma}(\overline{X}).
	\)
	Recall that $\rho=\rho_l$. This shows that $\overline{X}$ is a stationary 
	point of subproblem \eqref{EP-subprob}.
	
	\medskip
	
	{\bf Item 3.} Since $\vartheta(X)=\min_{Z\in\mathbb{R}_{+}^{n\times r}}\|Z\!-\!X\|_1$,
	it is a semi-algebraic function on $\mathbb{R}^{n\times r}$, which means that
	$e_{\gamma}\vartheta$ is a semi-algebraic function on $\mathbb{R}^{n\times r}$.
	Since $\delta_{M}$ is semi-algebraic, $e_{\gamma}\vartheta+\delta_{M}$
	is a semi-algebraic function on $\mathbb{R}^{n\times r}$ and hence is definable 
	in the o-minimal structure $\mathscr{O}$. Since $f$ is definable in the o-minimal 
	structure $\mathscr{O}$, $\widetilde{\Theta}_{\rho,\gamma}\!
	=\!f+\rho e_{\gamma}\vartheta+\delta_{M}$ is definable in the o-minimal 
	structure $\mathscr{O}$ by \cite[Section 4.3]{Attouch10}.
	By \cite[Theorem 4.1]{Attouch10}, $\widetilde{\Theta}_{\rho,\gamma}$ is a KL function. 
	Recall that the sequence $\{X^k\}_{k\in\mathbb{N}}$ is bounded.
	By the expression of ${\rm Proj}_{{\rm T}_{\!X}M}$ in \eqref{Proj-St},
	there exists a constant $c'>0$ such that for all $k\in\mathbb{N}$, 
	$\|[{\rm Proj}_{{\rm T}_{\!X^{k+1}}M}-{\rm Proj}_{{\rm T}_{\!X^{k}}M}]
	(\nabla\Theta_{\rho,\gamma}(X^{k}))\|_F\le c'\|X^{k+1}-X^k\|_F$.    
	Next we argue that there exists $k_0\in\mathbb{N}$ such that for all $k\ge k_0$, 
	the following two conditions hold:
	\begin{itemize}[(H1)]
		\item[(H1)] $\widetilde{\Theta}_{\rho,\gamma}(X^{k+1})
		+\frac{\alpha}{3t_{\rm max}}\|X^{k+1}-X^k\|_F^2
		\le \widetilde{\Theta}_{\rho,\gamma}(X^{\ell(k)})$;
		
		\item[(H2)] there exists $W^{k}\in\partial\widetilde{\Theta}_{\rho,\gamma}(X^{k})$ 
		such that $\|W^{k+1}\|_F\le 
		(L_{\nabla f}+\frac{2\rho}{\gamma}+c'+\frac{2}{t_{\rm min}})\|X^{k+1}-X^{k}\|_F$.
	\end{itemize}     
	Indeed, from the iterates of Algorithm \ref{ManPG1},  
	$\|X^{k+1}-X^k\|_F=\|R_{X^k}(V^k)-(X^k+V^k)+V^k\|_F$ for each $k\in\mathbb{N}$.
	Since $\lim_{k\to\infty}V^k=0$, by Remark \ref{remark-retraction}
	we have $\|R_{X^k}(V^k)-(X^k+V^k)\|_F=o(\|V^k\|_F)$. Then, 
	there exist $k_0\in\mathbb{N}$ such that for all $k\ge k_0$, 
	$\frac{1}{2}\|V^k\|_F\le\|X^{k+1}-X^k\|_F\le\frac{3}{2}\|V^k\|_F$.
	Together with Lemma \ref{lemma1-ManPG1} and  
	$\{X^k\}_{k\in\mathbb{N}}\subset M$, for each $k\ge k_0$ we have 
	$\widetilde{\Theta}_{\rho,\gamma}(X^{k+1})+\frac{\alpha}{3t_{\rm max}}\|X^{k+1}-X^k\|_F^2
	\le \widetilde{\Theta}_{\rho,\gamma}(X^k)$.
	This shows that condition (H1) is satisfied. 
	For each $k\in\mathbb{N}$, let $W^k={\rm grad}\Theta_{\rho,\gamma}(X^{k})$. 
	Noting that $\partial\widetilde{\Theta}_{\rho,\gamma}(X^k)
	=\nabla\Theta_{\rho,\gamma}(X^k)+{\rm N}_{X^k}M$,
	we have $W^k\in{\rm Proj}_{{\rm T}_{\!X^k}M}(\partial\widetilde{\Theta}_{\rho,\gamma}(X^k))
	\subseteq\partial\widetilde{\Theta}_{\rho,\gamma}(X^k)$, 
	where the inclusion is due to Lemma \ref{relation-subdiff}.
	Recall that $\|X^{k+1}-X^k\|_F\ge \frac{1}{2}\|V^k\|_F$ for each $k\ge k_0$.  
	Together with $V^k=-t_k{\rm grad}\Theta_{\rho,\gamma}(X^{k})=-t_kW^k$,
	we have $\|W^k\|_F\le\frac{2}{t_{\rm min}}\|X^{k+1}-X^k\|_F$. 
	Then, for each $k\ge k_0$,
	\begin{align*}
		\|W^{k+1}\|_F&\le\|W^{k+1}-W^{k}\|_F+\|W^k\|_F
		=\|{\rm Proj}_{{\rm T}_{\!X^{k+1}}M}(\nabla\Theta_{\rho,\gamma}(X^{k+1}))
		-{\rm Proj}_{{\rm T}_{\!X^{k}}M}(\nabla\Theta_{\rho,\gamma}(X^{k}))\|_F+\|W^k\|_F\\
		&\le\|\nabla\Theta_{\rho,\gamma}(X^{k+1})-\nabla\Theta_{\rho,\gamma}(X^{k})\|_F
		+\|[{\rm Proj}_{{\rm T}_{\!X^{k+1}}M}-{\rm Proj}_{{\rm T}_{\!X^{k}}M}]
		(\nabla\Theta_{\rho,\gamma}(X^{k}))\|_F+\|W^k\|_F\\
		&\le [L_{\nabla f}+ 2\rho/\gamma]\|X^{k+1}-X^k\|_F+c'\|X^{k+1}-X^k\|_F
		+\frac{2}{t_{\rm min}}\|X^{k+1}-X^k\|_F.
	\end{align*}
	This shows that condition (H2) holds. Thus, the sequence 
	$\{X^k\}_{k\in\mathbb{N}}$ satisfies conditions H1-H2 in \citep{QianPan22} 
	with $\Phi=\widetilde{\Theta}_{\rho,\gamma}$. 
	Since $\widetilde{\Theta}_{\rho,\gamma}$ is continuous relative to $M$ 
	and $\{X^k\}\subseteq M$, the assumption of \citep[Theorem 3.3]{QianPan22} 
	holds. By invoking \citep[Theorem 3.3]{QianPan22}, 
	we obtain the conclusion.
\end{proof}

\section*{Acknowledgments}
 The authors are deeply indebted to Professor Wen Zaiwen from Peking University 
 and Professor Jiang Bo from Nanjing Normal University for their codes sharing 
 and the helpful discussions about optimization with nonnegative orthogonal constraint. 
 The authors would like to thank the anonymous referees for their valuable suggestions and comments, which helped improve this paper. 

\section*{Funding}

 This work is supported by the National Natural Science Foundation of China
 under projects No.11971177,  the Guangdong Basic and Applied Basic Research 
 Foundation No. 2021A1515010210, and the Science and Technology Projects in Guangzhou No. 202201010566.

 \end{document}